\documentclass[reqno,a4paper]{amsart}
\usepackage{amsthm}
\usepackage[T1]{fontenc}
\usepackage[utf8]{inputenc}
\usepackage[english]{babel}
\usepackage{microtype}

\usepackage{amssymb, amsmath}
\usepackage{enumitem}
\usepackage{tikz-cd}
\usepackage{orcidlink}
\usepackage{doi}
\renewcommand{\doi}[1]{DOI:~\href{https://doi.org/#1}{#1}}
\usepackage{url}

\usepackage[all]{xy}
\usepackage{calrsfs}
\usepackage{mathabx}
\usepackage{xcolor}
\usepackage{mathtools}
\usepackage{xurl}
\usepackage{nameref,hyperref,cleveref}
\hypersetup{hidelinks}

\theoremstyle{plain}
\newtheorem{theorem}[subsection]{Theorem}
\newtheorem{proposition}[subsection]{Proposition}
\newtheorem{corollary}[subsection]{Corollary}
\newtheorem{lemma}[subsection]{Lemma}

\theoremstyle{definition}
\newtheorem{definition}[subsection]{Definition}

\theoremstyle{remark}
\newtheorem{remark}[subsection]{Remark}
\newtheorem{example}[subsection]{Example}

\crefname{examples}{Examples}{Examples}
\Crefname{examples}{Examples}{Examples}
\crefname{problem}{Open problem}{Open problems}
\Crefname{problem}{Open problem}{Open problems}

\numberwithin{equation}{section}

\newenvironment{tfae}
{
\begin{enumerate}}
{\end{enumerate}}

\newcommand{\noproof}{\hfil\qed}

\newcommand{\cC}{\ensuremath{\mathcal{C}}}

\newcommand{\F}{\mathbb{F}}         
\newcommand{\E}{\mathcal{E}}

\newcommand{\N}{\mathbb{N}}
\newcommand{\V}{\mathcal{V}}
\newcommand{\W}{\mathcal{W}}

\DeclareMathOperator{\bchar}{char}
\DeclareMathOperator{\Ker}{Ker}
\DeclareMathOperator{\Aut}{Aut}

\DeclareMathOperator{\op}{op}
\DeclareMathOperator{\Act}{Act}

\DeclareMathOperator{\Hom}{Hom}
\DeclareMathOperator{\Bim}{Bim}
\DeclareMathOperator{\SplExt}{SplExt}
\DeclareMathOperator{\Imm}{Im}
\DeclareMathOperator{\End}{End}
\DeclareMathOperator{\Der}{Der}

\DeclareMathOperator{\Bider}{Bider}
\DeclareMathOperator{\Inn}{Inn}
\DeclareMathOperator{\M}{M}

\DeclareMathOperator{\USGA}{USGA}
\DeclareMathOperator{\Ann}{Ann}

\DeclareMathOperator{\id}{id}

\newcommand{\Set}{\ensuremath{\mathsf{Set}}}

\newcommand{\Assoc}{\ensuremath{\mathsf{Assoc}}}

\newcommand{\Lie}{\ensuremath{\mathsf{Lie}}}
\newcommand{\PreLie}{\ensuremath{\mathsf{PreLie}}}
\newcommand{\Leib}{\ensuremath{\mathsf{Leib}}}
\newcommand{\Pois}{\ensuremath{\mathsf{Pois}}}
\newcommand{\CPois}{\ensuremath{\mathsf{CPois}}}

\newcommand{\CAssoc}{\ensuremath{\mathsf{CAssoc}}}

\newcommand{\Com}{\ensuremath{\mathsf{Com}}}
\newcommand{\ACom}{\ensuremath{\mathsf{ACom}}}
\newcommand{\JJord}{\ensuremath{\mathsf{JJord}}}

\newcommand{\AbAlg}{\ensuremath{\mathsf{AbAlg}}}
\newcommand{\Alg}{\ensuremath{\mathsf{Alg}}}
\newcommand{\PAlg}{\ensuremath{\mathsf{PAlg}}}
\newcommand{\Nil}{\ensuremath{\mathsf{Nil}}}
\newcommand{\Sol}{\ensuremath{\mathsf{Sol}}}
\newcommand{\Alt}{\ensuremath{\mathsf{Alt}}}

\newcommand{\Grp}{\ensuremath{\mathsf{Grp}}}

\newcommand{\Vect}{\ensuremath{\mathsf{Vec}}}

\newcommand{\NP}{\ensuremath{\mathsf{NPois}}}

\newcommand{\SplExtX}{\ensuremath{\mathsf{SplExt}(X)}}

\newcommand{\Match}{\ensuremath{\mathsf{MDA}}}
\newcommand{\Dias}{\ensuremath{\mathsf{Dias}}}
\newcommand{\CLie}{\ensuremath{\mathsf{CLie}}}
\newcommand{\MatchLie}{\ensuremath{\mathsf{MLie}}}
\newcommand{\PostLie}{\operatorname{\mathsf{PostLie}}}
\newcommand{\Tetra}{\operatorname{Tetra}}

\date{}

\begin{document}

\title[Actions of algebras with several bilinear multiplications]{On the representability of actions of algebras with several bilinear multiplications}

\author[M.~Mancini]{Manuel Mancini~\orcidlink{0000-0003-2142-6193}}

\email{manuel.mancini@uclouvain.be; manuel.mancini@unipa.it}

\address{Institut de Recherche en Mathématique et Physique, Université catholique de Louvain, chemin du cyclotron 2 bte L7.01.02, B--1348 Louvain-la-Neuve, Belgium.}
\address{Dipartimento di Matematica e Informatica, Università degli Studi di Palermo, Via Archirafi 34, 90123 Palermo, Italy.}

\begin{abstract}
We study the representability of actions in varieties of algebras with $n \geq 1$ not necessarily associative bilinear multiplications. We first prove that, for homogeneous varieties, action accessibility is equivalent to being an Orzech category of interest. For an action accessible, operadic variety we then construct the \emph{external weak actor} $\E(X)$ of an algebra $X$, an algebra with $n$ bilinear partial operations together with a natural monomorphism encoding the actions on~$X$. We characterise the algebras admitting an actor in terms of $\E(X)$, and we deduce an obstruction to action representability: if, for some algebra $X$ with trivial multiplications, $\E(X)$ does not satisfy a multilinear identity of the variety, then the variety is not action representable. Computing the external weak actor explicitly, we obtain that associative and matching dialgebras are weakly action representable but not action representable, and that Poisson algebras, left-right non-commutative Poisson algebras and matching Lie algebras are not action representable. Finally, we introduce the \emph{reduced type} $(p,q)$ of such a variety and prove that, over a field of characteristic different from~$2$, an operadic variety of reduced type $(p,0)$, $(0,q)$ or $(1,1)$ is action representable if and only if it is term equivalent to the variety of Lie algebras or to that of abelian algebras.
For the reduced types $(p,0)$ and $(0,q)$, the arguments require no assumption on the cardinality of the field, whereas the case $(1,1)$ requires the field to be infinite.
\end{abstract}


\maketitle

\section{Introduction}\label{SecIntro}

Semi-abelian categories were introduced by G.~Janelidze, L.~Márki, and W.~Tholen in~\cite{Semi-Ab} as a categorical framework encompassing many algebraic categories, including groups, rings, and varieties of algebras. Within this setting, internal actions are naturally described by split extensions~\cite{BJK}. More precisely, for every object $X$ of a semi-abelian category $\cC$, there is a natural isomorphism
\[
\Act(-,X)\cong\SplExt(-,X),
\]
where $\Act(B,X)$ denotes the set of internal actions of $B$ on $X$, and $\SplExt(B,X)$ denotes the set of isomorphism classes of split extensions of $B$ by $X$.

A semi-abelian category $\cC$ is said to be \emph{action representable} if, for every object~$X$ of $\cC$, the functor $\SplExt(-,X)$ is representable; in this case, the representing object~$[X]$ is called the \emph{actor} of $X$. Classical examples of action representable categories are the category of groups and the category of Lie algebras over a unital commutative ring. However, action representability is a restrictive condition: over an infinite field, the only action representable operadic varieties of non-associative algebras with a single bilinear multiplication are the varieties of abelian algebras, of Lie algebras and of quasi-Lie algebras, the last two coinciding in characteristic different from~$2$; see~\cite[Theorem~6.1]{Tim}.

Several notions weaker than action representability have been introduced in the literature. D.~Bourn and G.~Janelidze defined the notion of \emph{action accessible category} in~\cite{act_accessible}, a class that includes relevant examples such as rings, associative algebras, and Leibniz algebras. More recently, G.~Janelidze introduced in~\cite{WRA} the notion of weakly representable action: a category $\cC$ is \emph{weakly action representable} if, for every object $X$ of $\cC$, there exist an object $T=T(X)$ of $\cC$ and a natural monomorphism of functors
\[
\SplExt(-,X)\rightarrowtail\Hom_{\cC}(-,T).
\]
G.~Janelidze also proved that every weakly action representable semi-abelian category is action accessible, so that one has the implications
\[
\textit{action representability}
\Longrightarrow
\textit{weak action representability}
\Longrightarrow
\textit{action accessibility}.
\]
Orzech categories of interest~\cite{Orzech} provide an important algebraic setting in which internal actions can be studied explicitly. For every object $X$ of a category of interest $\cC$, J.~M.~Casas, T.~Datuashvili, and M.~Ladra constructed in~\cite{Casas} an object $\USGA(X)$, called the \emph{universal strict general actor} of $X$, which lives in a suitable larger category~$\cC'$ and encodes actions on~$X$. When $\USGA(X)$ belongs to $\cC$, it determines a weak representation of actions on $X$.

For operadic varieties of non-associative algebras with a single bilinear multiplication, action accessibility can be characterised by the existence of the $\lambda/\mu$-rules, which give an explicit formulation of Axiom~2 for an Orzech category of interest. In~\cite{WRAAlg}, this characterisation was used to associate with every algebra $X$ in such a variety $\V$, a partial algebra $\E(X)$, called the \emph{external weak actor} of $X$, together with a natural monomorphism of functors
\[
\SplExt(-,X)\rightarrowtail\Hom_{\PAlg}(\widetilde U(-),\E(X)),
\]
where $\PAlg$ is the category of partial algebras and $\widetilde U \colon \V \to \PAlg$ denotes the forgetful functor. In the present article, we adopt the formulation of partial algebras of~\cite{SurveyActions}, which modifies the one of~\cite{WRAAlg}. Namely, a \emph{partial algebra} is a triple $(X,D,\mu)$, where $D\leq X\otimes X$ is a vector subspace and $\mu\colon D\to X$ is a linear map, and the product~$xy$ is defined if and only if $x\otimes y\in D$, in which case $xy=\mu(x\otimes y)$.

Many algebraic structures of interest in algebra, geometry and mathematical physics carry more than one bilinear multiplication, related by compatibility conditions: examples are Poisson algebras, non-commutative Poisson algebras~\cite{CasasDatuashviliNLP,CasasDatuashviliLadraNPLR}, associative dialgebras~\cite{LodayDialgebras}, matching dialgebras~\cite{ZhangBaiGuo}, compatible Lie algebras~\cite{Strohmayer}, matching Lie algebras~\cite{GaoGuoZhang} and post-Lie algebras~\cite{VallettePostLie}. Actions in some of the corresponding varieties have been studied by means of constructions specific to each of them, for instance in~\cite{CigoliManciniMetere} for Poisson algebras and in~\cite{CasasEtAlTetra} for associative dialgebras. Our motivation for a systematic study is threefold. First, for a category of interest, the universal strict general actor of~\cite{Casas} lives in a larger category and is in general more difficult to describe than the external weak actor, which, however, has only been defined for a single multiplication. Second, the behaviour of a variety with respect to actions is not determined by its multiplications taken separately. For instance, the variety of compatible Lie algebras, in which every linear combination of the two brackets is a Lie bracket, is not action accessible, and neither is the variety $\Lie^n$ of algebras with $n\geq2$ Lie brackets subject to no compatibility condition, although $\Lie$ is action representable. Conversely, a variety may be action accessible even if the variety determined by the identities involving only one of its multiplications is not; see \Cref{PropExamplesNOperations} and \Cref{ExampleAccessibleReductNotAccessible}. Third, in view of the classification recalled above, it is natural to ask whether new action representable varieties arise when several multiplications are allowed; this question was raised in~\cite[Section~6.6]{Tim}.

The first aim of the present article is to characterise action accessibility for homogeneous varieties of algebras with $n$ bilinear multiplications. We prove that such a variety is action accessible if and only if it is an Orzech category of interest. We then restrict our attention to action accessible, operadic varieties. For every algebra~$X$ of such a variety~$\V$, we construct an algebra~$\E(X)$ with $n$ bilinear partial operations, called the \emph{external weak actor} of $X$, together with a natural monomorphism of functors
\[
\SplExt(-,X)\rightarrowtail\Hom_{\PAlg^n}(\widetilde U(-),\E(X)),
\]
where $\PAlg^n$ denotes the category of algebras with $n$ bilinear partial operations and $\widetilde U\colon\V\to\PAlg^n$ is the forgetful functor. We prove that $X$ admits an actor if and only if the elements of $\E(X)$ induced by the derived actions on $X$ form a subalgebra~$S_X$ of $\E(X)$ such that the semidirect product $S_X\ltimes X$ is an object of~$\V$. In this case, $S_X$ is the actor of $X$; see~\Cref{ThmActorSX}. We deduce an obstruction to action representability: if, for some algebra $X$ of $\V$ with trivial multiplications, $\E(X)$ does not satisfy a multilinear identity of $\V$ of degree~$k$, then the algebra $X\otimes\Lambda_k$ with trivial multiplications, where $\Lambda_k=\F[\xi_1,\ldots,\xi_k]/(\xi_1^2,\ldots,\xi_k^2)$, does not admit an actor. In particular, $\V$ is not action representable. If $X$ is infinite-dimensional, then $X\otimes\Lambda_k\cong X$, so that $X$ itself does not admit an actor; see~\Cref{ThmTensorObstruction} and \Cref{CorTensorObstruction}.

We then describe explicitly the construction of the external weak actor for several varieties. For Poisson and commutative Poisson algebras, it recovers the constructions of~\cite{CigoliManciniMetere}, and for associative dialgebras the tetramultipliers of~\cite{CasasEtAlTetra}. We prove that associative and matching dialgebras are weakly action representable, and that the varieties of Poisson algebras, left-right non-commutative Poisson algebras, associative dialgebras, matching dialgebras and matching Lie algebras are not action representable.

Finally, we study the action representability of varieties of algebras with $n$ bilinear multiplications. Up to term equivalence, the multiplications of such a variety $\V$ may always be assumed to be linearly independent and each of them either commutative or anti-commutative. The numbers $p$ and $q$ of the commutative and of the anti-commutative ones do not depend on the choice of such a presentation, and the pair~$(p,q)$ is called the \emph{reduced type} of $\V$. Our main tool is the obstruction described above, from which it follows that the existence of an actor for a suitable infinite-dimensional algebra with trivial multiplications forces the identities for Axiom~2 to be encoded by a square matrix $A$ satisfying $A^2=I$ and $AEA=-EAE$, where $E$ is the diagonal matrix whose entries encode the signs of the multiplications. We deduce that, over a field $\F$ with $\bchar(\F)\neq2$, an operadic variety of reduced type $(p,0)$, $(0,q)$ or $(1,1)$ is action representable if and only if it is term equivalent to either the variety of Lie algebras or the variety of abelian algebras. This answers the question recalled above for these reduced types. For the reduced types $(p,0)$ and $(0,q)$, the arguments hold over an arbitrary field of characteristic different from~$2$. In the case of a single multiplication, they give new proofs of the corresponding results of~\cite{Tim}, without assuming the field to be infinite. The reduced type $(1,1)$ is deduced from the classification for a single multiplication recalled above, and it is the only case in which the field is required to be infinite; see \Cref{RemarkComparisonSingle} for a detailed comparison.

The article is organised as follows. In \Cref{SecPreliminaries}, we recall the categorical and algebraic background concerning internal actions, (weak) action representability, Orzech categories of interest, partial algebras, and the external weak actor for varieties with one multiplication. In \Cref{SecNOperations}, we introduce varieties of algebras with a finite number of (not necessarily associative) bilinear multiplications, we characterise action accessibility and study derived actions in this setting. In \Cref{SecExternalWeakActorN}, we construct the external weak actor of an algebra with $n$ multiplications, characterise the algebras admitting an actor, and obtain an obstruction to action representability. In \Cref{SecEx}, we describe the external weak actor for the varieties of algebras with two multiplications discussed above that are Orzech categories of interest. Finally, in \Cref{SecNLieBracketsAR2}, we introduce the reduced type $(p,q)$ of a variety of algebras with $n \geq 1$ bilinear multiplications and we classify the action representable, operadic varieties of reduced type $(p,0)$, $(0,q)$ and $(1,1)$.

\section{Preliminaries}\label{SecPreliminaries}

In this section, we recall the notions and results concerning varieties of algebras with one not necessarily associative bilinear multiplication, internal actions in semi-abelian categories, and the construction of the external weak actor that will be used throughout the text.

\subsection{Varieties of non-associative algebras}\label{SecVarieties}

Let $\F$ be a field. A \emph{non-associative algebra} over $\F$ is an $\F$-vector space $X$ endowed with a bilinear multiplication
\[
X\times X\to X\colon(x,y)\mapsto xy.
\]
No associativity or existence of a multiplicative identity is assumed. We denote by $\Alg$ the category of non-associative algebras over $\F$ and algebra homomorphisms.

A \emph{non-associative polynomial} in the variables $x_1,\ldots,x_t$ is an $\F$-linear combination of non-associative words in these variables. An \emph{identity} of an algebra $X$ is a non-associative polynomial
\[
\Phi=\Phi(x_1,\ldots,x_t)
\]
such that
\[
\Phi(a_1,\ldots,a_t)=0
\]
for all $a_1,\ldots,a_t\in X$. In this case, we say that $X$ \emph{satisfies} the identity~$\Phi$.

Given a set $I$ of identities, the \emph{variety of non-associative algebras} determined by~$I$ is the full subcategory $\V$ of $\Alg$ consisting of all algebras satisfying every identity in~$I$. 

For a monomial $\varphi=\varphi(x_1,\ldots,x_t)$, its \emph{type} is the $t$-tuple $(k_1,\ldots,k_t)\in\mathbb{N}^t$, where~$k_i$ is the number of occurrences of $x_i$ in $\varphi$. A polynomial is said to be \emph{homogeneous} if all its monomials have the same type, and it is \emph{multilinear} if all its monomials have type $(1,\ldots,1)$. A variety of non-associative algebras is called \emph{homogeneous} if it can be determined by a set of homogeneous identities, and it is called \emph{operadic}~\cite{Cosmash} if it can be determined by a set of multilinear identities. In particular, every operadic variety is homogeneous. We recall that, if $\bchar(\F)=0$, then any variety $\V$ is operadic. This follows from the multilinearisation process, see~\cite[Corollary 3.7]{Osborn}.

\begin{remark}
Throughout the paper, we exclude the trivial variety defined by the identity $x=0$, whose only object is the zero algebra.
\end{remark}

\begin{example}\label{ExamplesVarieties}
The following are standard examples of varieties of non-associative algebras.
\begin{enumerate}
\item The variety $\AbAlg$ of \emph{abelian algebras} is determined by the identity
\[
xy=0.
\]

\item The varieties $\Com$ and $\ACom$ of commutative and anti-commutative algebras are determined, respectively, by the identities
\[
xy-yx=0
\qquad\text{and}\qquad
xy+yx=0.
\]

\item The variety $\Assoc$ of associative algebras is determined by the identity
\[
(xy)z-x(yz)=0.
\]

\item The variety $\CAssoc$ of commutative associative algebras is determined by the identities
\[
xy-yx=0 \qquad\text{and}\qquad(xy)z-x(yz)=0.
\]

\item The variety $\Alt$ of alternative algebras is determined by the identities
\[
x(xy)=(xx)y
\qquad\text{and}\qquad
(yx)x=y(xx).
\]

\item The variety $\JJord$ of Jacobi--Jordan algebras~\cite{Burde} is determined by commutativity and the Jacobi identity
\[
x(yz)+y(zx)+z(xy)=0.
\]

\item The variety $\Lie$ of Lie algebras consists of the non-associative algebras whose bilinear multiplication, denoted by $[-,-]$, satisfies $[x,x]=0$ and the Jacobi identity. When $\bchar(\F)\neq2$, the identity $[x,x]=0$ is equivalent to $[x,y]+[y,x]=0$, and hence in this case the variety $\Lie$ is operadic.

\item The variety $\Leib$ of (right) Leibniz algebras~\cite{loday} consists of the non-associative algebras whose bilinear multiplication, also denoted by $[-,-]$, satisfies the (right) Leibniz identity
\[
[[x,y],z]=[[x,z],y]+[x,[y,z]].
\]

\item The variety $\PreLie$ of pre-Lie algebras~\cite{Ger63,Vin63} is determined by the identity
\[
a(x,y,z)-a(y,x,z)=0,
\]
where $a(x,y,z)=x(yz)-(xy)z$ denotes the associator.

\item Given a variety $\V$, we denote by $\Nil_k(\V)$ the subvariety of $k$-nilpotent algebras of $\V$, determined by all the identities of $\V$ together with
\[
x_1 x_2 \cdots x_{k+1}=0
\]
for all possible choices of parentheses.
\end{enumerate}
\end{example}

\subsection{Representability of actions in semi-abelian categories}\label{SecActions}

Semi-abelian categories were introduced by G.~Janelidze, L.~Márki, and W.~Tholen in~\cite{Semi-Ab} to provide a categorical setting capturing common categorical-algebraic properties of groups, rings, and algebras. We recall that a category $\cC$ is \emph{semi-abelian} if it is pointed, Barr-exact~\cite{Barr}, protomodular~\cite{Bourn}, and has finite coproducts. Standard examples include the categories of groups, Lie algebras, and, more generally, any variety of non-associative algebras over a field; see~\cite[Theorem~9.5]{VdL-NAA}.

Let $B$ and $X$ be objects of a semi-abelian category $\cC$. A \emph{split extension} of $B$ by~$X$ is a diagram
\begin{equation}\label{eq:split_extension}
\begin{tikzcd}
X \arrow[r, "k"]
&
A \arrow[r, shift left, "\alpha"]
&
B \arrow[l, shift left, "\beta"]
\end{tikzcd}
\end{equation}
in $\cC$ such that $\alpha\circ\beta=\id_B$ and $(X,k)$ is a kernel of~$\alpha$.

For every object $X$ of $\cC$, one obtains a functor
\[
\SplExt(-,X)\colon\cC^{\op}\to\Set
\]
which assigns to an object $B$ the set $\SplExt(B,X)$ of isomorphism classes of split extensions of $B$ by $X$. A morphism $f\colon B'\to B$ is sent to the change of base map
\[
f^*\colon\SplExt(B,X)\to\SplExt(B',X)
\]
defined by pulling back split extensions along~$f$.

Split extensions in a semi-abelian category are naturally equivalent to internal actions through a semidirect product construction~\cite{BJK}. More precisely, for every object $X$ of $\cC$, there is a natural isomorphism
\[
\Act(-,X)\cong\SplExt(-,X),
\]
where $\Act(B,X)$ denotes the set of internal actions of $B$ on $X$. For our purposes, we do not need an explicit description of internal actions, since split extensions are easier to work with, especially in varieties of non-associative algebras.

\begin{definition}\cite{BJK}
A semi-abelian category $\cC$ is \emph{action representable} if, for every object $X$ of $\cC$, the functor $\SplExt(-,X)$ is representable. Equivalently, for every object $X$, there exist an object $[X]$ of $\cC$, called the \emph{actor} of $X$, and a natural isomorphism of functors
\[
\SplExt(-,X)\cong\Hom_{\cC}(-,[X]).
\]
\end{definition}

For example, the category $\Grp$ of groups is action representable, with
\[
[X]=\Aut(X),
\]
the group of automorphisms of $X$, and so is the category $\Lie$ of Lie algebras, with
\[
[X]=\Der(X),
\]
the Lie algebra of derivations of $X$. Every abelian category is also action representable, with the zero object as the actor of every object.

A weaker notion was introduced by G.~Janelidze in~\cite{WRA}.

\begin{definition}\cite{WRA}
A semi-abelian category $\cC$ is \emph{weakly action representable} if, for every object $X$ of $\cC$, there exist an object $T=T(X)$ of $\cC$ and a natural monomorphism of functors
\[
\tau\colon\SplExt(-,X)\rightarrowtail\Hom_{\cC}(-,T).
\]
The pair $(T,\tau)$ is called a \emph{weak representation} of $\SplExt(-,X)$, and $T$ is called a \emph{weak representing object} of~$X$. For every object $B$ of $\cC$, the morphisms in the image of the component $\tau_B$ are called \emph{acting morphisms}.
\end{definition}

Another weakening of action representability is obtained by considering faithful split extensions. A split extension with kernel $X$ is called \emph{faithful} if it is a subterminal object in the category $\SplExtX$ of split extensions with kernel~$X$, while it is called \emph{accessible} if it has a morphism into a faithful object.

\begin{definition}\cite{act_accessible}
A semi-abelian category $\cC$ is \emph{action accessible} if, for any object~$X$ of $\cC$, every object in $\SplExtX$ is accessible.
\end{definition}

Every action representable category is weakly action representable, and every weakly action representable semi-abelian category is action accessible~\cite{WRA}. Thus,
\[
\textit{action representability}
\Longrightarrow
\textit{weak action representability}
\Longrightarrow
\textit{action accessibility}.
\]
Neither implication is reversible in general. Indeed, the variety $\Assoc$ is weakly action representable~\cite{WRA}, but not action representable~\cite{Tim}. Moreover, J.~R.~A.~Gray proved in~\cite{Gray} that the varieties $\Sol_n(\Grp)$ of $n$-solvable groups are action accessible but not weakly action representable for every $n\geq3$. More recently, this was refined in~\cite{XabiMancini}, where the same is proved for $\Sol_2(\Grp)$ and for the varieties $\Nil_k(\Grp)$ of $k$-nilpotent groups, for every $k\geq3$. The analogous results hold for the varieties $\Sol_n(\Lie)$ and $\Nil_k(\Lie)$, for every $n\geq2$ and $k\geq3$, respectively.

\subsection{Orzech categories of interest and derived actions}\label{SecCategoriesInterest}

We now recall the notion of an Orzech category of interest~\cite{Orzech}. Let $\cC$ be a category of groups with a set of operations $\Omega$ and a set of identities $\mathbb{E}$ containing the group identities. For any~$i\in \N$, let $\Omega_i$ denote the set of $i$-ary operations in~$\Omega$. The following conditions are required:
\begin{enumerate}
\item $\Omega=\Omega_0\cup\Omega_1\cup\Omega_2$.

\item The group operations, written additively as $0$, $-$, and $+$, belong respectively to $\Omega_0$, $\Omega_1$, and $\Omega_2$. Moreover, $\Omega_0=\{0\}$. We set
\[
\Omega_1'=\Omega_1\setminus\{-\},
\qquad
\Omega_2'=\Omega_2\setminus\{+\}.
\]
For every $*\in\Omega_2'$, the opposite operation
\[
x*^{\op}y=y*x
\]
also belongs to~$\Omega_2'$.

\item For every $*\in\Omega_2'$, the identity
\[
x*(y+z)=x*y+x*z
\]
belongs to~$\mathbb{E}$.

\item For every $\omega\in\Omega_1'$ and $*\in\Omega_2'$, the identities
\[
\omega(x+y)=\omega(x)+\omega(y)
\qquad\text{and}\qquad
\omega(x*y)=\omega(x)*y
\]
belong to~$\mathbb{E}$.
\end{enumerate}

The following two additional axioms are imposed:
\begin{enumerate}
\item[\textnormal{Axiom~1.}] For every $*\in\Omega_2'$, the identity
\[
x_1+(x_2*x_3)=(x_2*x_3)+x_1
\]
belongs to~$\mathbb{E}$.

\item[\textnormal{Axiom~2.}] For every ordered pair $(*,\bar{*})\in\Omega_2'\times\Omega_2'$, there exists a word $W_{*,\bar{*}}(x_1,x_2;x_3)$ in the expressions
\[
x_1(x_2x_3),\quad
x_1(x_3x_2),\quad
(x_2x_3)x_1,\quad
(x_3x_2)x_1,
\]
\[
x_2(x_1x_3),\quad
x_2(x_3x_1),\quad
(x_1x_3)x_2,\quad
(x_3x_1)x_2,
\]
where each juxtaposition represents an operation in $\Omega_2'$, such that the identity
\[
(x_1*x_2)\mathbin{\bar{*}}x_3=W_{*,\bar{*}}(x_1,x_2;x_3)
\]
belongs to~$\mathbb{E}$.
\end{enumerate}

\begin{definition}\cite{Orzech}\label{DefCategoryInterest}
A category of groups with operations satisfying all the conditions above is called an \emph{Orzech category of interest}, or simply a \emph{category of interest}.
\end{definition}

Let $\mathbb{E}'$ be the subset of $\mathbb{E}$ containing the group identities and the identities in~(3) and~(4) (and, for categories of $\F$-algebras, also the identities expressing that $+$ is commutative and that the operations $\omega_\lambda$ define a structure of $\F$-vector space), and let $\cC'$ be the corresponding category of groups with operations. Then~$\cC$ is a full subcategory of $\cC'$, and there is a forgetful functor $\overline{U}\colon\cC\to\cC'$. The category~$\cC'$ is called the \emph{general category of groups with operations} associated with~$\cC$. Here and in what follows, for every $*\in\Omega_2'$ we regard $*^{\op}$ as a derived operation, that is, we also include in $\mathbb{E}'$ the identity $x*^{\op}y=y*x$.

When the underlying group operation is commutative, Axiom~1 is automatically satisfied. In particular, for a variety of non-associative algebras $\V$ with one bilinear multiplication $*$, we have
\[
\Omega_1'=\{\omega_\lambda\mid\lambda\in\F\} \quad \text{ and } \quad \Omega_2'=\{*,*^{\op}\},
\]
where $\omega_\lambda(x)=\lambda x$, and Axiom~2 is equivalent to requiring that the products $x(yz)$ and $(yz)x$ can be expressed as linear combinations of the eight products displayed above. These expressions are precisely the $\lambda/\mu$-rules appearing in the following theorem.

\begin{theorem}\cite{XabiTim,Coco_Thesis}\label{ThmLambdaMu}
Let $\V$ be a homogeneous variety of non-associative algebras. The following conditions are equivalent:
\begin{enumerate}
\item $\V$ is an Orzech category of interest;
\item $\V$ is action accessible;
\item there exist $\lambda_1,\ldots,\lambda_8,\mu_1,\ldots,\mu_8\in\F$ such that
\begin{align}
x(yz)=\lambda_1(xy)z &+\lambda_2(yx)z+\lambda_3z(xy)+\lambda_4z(yx)\notag\\
&+\lambda_5(xz)y+\lambda_6(zx)y+\lambda_7y(xz)+\lambda_8y(zx)\label{LambdaRule}
\end{align}
and
\begin{align}
(yz)x=\mu_1(xy)z &+\mu_2(yx)z+\mu_3z(xy)+\mu_4z(yx)\notag\\
&+\mu_5(xz)y+\mu_6(zx)y+\mu_7y(xz)+\mu_8y(zx)\label{MuRule}
\end{align}
are identities of $\V$.
\end{enumerate}
\end{theorem}

The identities \eqref{LambdaRule} and \eqref{MuRule} are called the \emph{$\lambda/\mu$-rules}. Their choice need not be unique. From now on, whenever we consider an action accessible variety of non-associative algebras, we fix a choice of $\lambda/\mu$-rules.

G.~Orzech also introduced the notion of a \emph{derived action} in a category of groups with operations~\cite{Orzech}; see also~\cite{Datuashvili}. Under the correspondence between internal actions and split extensions, derived actions provide an algebraic description of internal actions; see~\cite[Theorem~3.17]{SurveyActions}. We recall this description for varieties of non-associative algebras.

Let $B$ and $X$ be objects of $\V$, and consider a split extension
\begin{equation}\label{SplitExternalActor}
\begin{tikzcd}
X \arrow[r, "k"]
&
A \arrow[r, shift left, "\alpha"]
&
B. \arrow[l, shift left, "\beta"]
\end{tikzcd}
\end{equation}
Since $(X,k)$ is a kernel of $\alpha$, we can define the bilinear maps
\[
l\colon B\times X\to X,
\qquad
r\colon X\times B\to X
\]
by
\[
b\ast x=l(b,x)=k^{-1}\bigl(\beta(b)k(x)\bigr),
\qquad
x\ast b=r(x,b)=k^{-1}\bigl(k(x)\beta(b)\bigr)
\]
for every $b\in B$ and $x\in X$. These expressions are well defined because
\[
\alpha\bigl(\beta(b)k(x)\bigr)=
\alpha\bigl(k(x)\beta(b)\bigr)=0.
\]

\begin{definition}
The pair $(l,r)$ is called the \emph{derived action} of $B$ on $X$ associated with~\eqref{SplitExternalActor}.
\end{definition}

Conversely, given two bilinear maps
\[
l\colon B\times X\to X,
\qquad
r\colon X\times B\to X,
\]
we define a bilinear multiplication on the direct sum $B\oplus X \cong B \times X$ by
\[
(b,x)(b',x')=\bigl(bb',xx'+b\ast x'+x\ast b'\bigr).
\]
The pair $(l,r)$ is a derived action of $B$ on $X$ if and only if the resulting algebra is an object of $\V$; see~\cite[Theorem~2.4]{Orzech}. In this case, the resulting algebra is called the \emph{semidirect product} of $B$ and $X$ associated with the derived action $(l,r)$ and is denoted by~$B\ltimes X$.

\begin{example}\label{ExampleDerivedAssoc}
Let $B$ and $X$ be associative algebras. A pair of bilinear maps
\[
B\times X\to X\colon(b,x)\mapsto b\ast x,
\qquad
X\times B\to X\colon(x,b)\mapsto x\ast b
\]
is a derived action of $B$ on $X$ in $\Assoc$ if and only if
\begin{align*}
(bb')\ast x&=b\ast(b'\ast x),&
x\ast(bb')&=(x\ast b)\ast b',\\
(b\ast x)\ast b'&=b\ast(x\ast b'),&
b\ast(xx')&=(b\ast x)x',\\
(xx')\ast b&=x(x'\ast b),&
(x\ast b)x'&=x(b\ast x')
\end{align*}
for every $b,b'\in B$ and $x,x'\in X$. These conditions are precisely those ensuring that the multiplication
\[
(b,x)(b',x')
=
\bigl(bb',xx'+b\ast x'+x\ast b'\bigr)
\]
on $B\oplus X$ is associative.
\end{example}

\subsection{Partial algebras and the external weak actor}\label{SecPartialActors}

We now recall the notion of a partial algebra and the construction of the external weak actor introduced in~\cite{WRAAlg}, using the modified formulation of partial algebras of~\cite[Section~4.21]{SurveyActions}.

In~\cite{WRAAlg}, the domain $\Omega\subseteq X\times X$ of a partial multiplication on a vector space $X$ was required to be a vector subspace. In the constructions under consideration, products with zero are always defined, so $(x,0),(0,y)\in\Omega$ for all $x,y\in X$. The subspace requirement then gives
\[
(x,y)=(x,0)+(0,y)\in\Omega
\]
for every $x,y\in X$. Thus $\Omega=X\times X$, and the multiplication is necessarily total. To allow non-total multiplications, the domain is instead taken to be a vector subspace of $X\otimes X$, as follows.

\begin{definition}\cite{SurveyActions}
Let $X$ be an $\F$-vector space. A \emph{bilinear partial operation} on~$X$ is a pair $(D,\mu)$ consisting of a vector subspace
\[
D\leq X\otimes X
\]
and a linear map
\[
\mu\colon D\to X.
\]
Given $x,y\in X$, we say that the product $xy$ is \emph{defined} if $x\otimes y\in D$. In this case, we set
\[
xy\coloneqq\mu(x\otimes y).
\]
\end{definition}

\begin{definition}\cite{SurveyActions}
A \emph{partial algebra} over $\F$ is a triple $(X,D,\mu)$, where $X$ is an $\F$-vector space and $(D,\mu)$ is a bilinear partial operation on $X$.

The partial algebra $(X,D,\mu)$ is said to be \emph{total} if $D=X\otimes X$. In this case, $\mu$ is the linearisation of an ordinary bilinear multiplication on $X$.
\end{definition}

\begin{definition}\cite{SurveyActions}
Let $(X,D,\mu)$ and $(X',D',\mu')$ be partial algebras over $\F$. A \emph{homomorphism of partial algebras} is a linear map $f\colon X\to X'$ such that
\[
(f\otimes f)(D)\subseteq D'
\]
and
\[
f\circ\mu=\mu'\circ\left.(f\otimes f)\right|_D.
\]
\end{definition}

We denote by $\PAlg$ the category whose objects are partial algebras over $\F$ and whose morphisms are homomorphisms of partial algebras.

Every non-associative algebra $X$ determines the total algebra
\[
\widetilde U(X)=(X,X\otimes X,\widetilde m_X),
\]
where $\widetilde m_X\colon X\otimes X\to X$ is the linearisation of the multiplication of $X$. This construction defines a canonical fully faithful functor $\widetilde U\colon\V\to\PAlg$ for every variety of non-associative algebras $\V$. We refer to~\cite[Section~4.21]{SurveyActions} for the evaluation of polynomial identities in partial algebras. The case of $n$ operations is treated in \Cref{SecPartialAlgebrasN}.

Now, let $\V$ be an action accessible, operadic variety determined by a set of multilinear identities
\[
\Phi_{t,r}(x_1,\ldots,x_t)=0,
\]
where $t=\deg \Phi_{t,r}\geq2$ (an identity of degree~$1$ only occurs in the trivial variety, which we exclude), and fix a choice of $\lambda/\mu$-rules. Let $X$ be an object of $\V$, and consider a pair
\[
f=(f\ast-,-\ast f)\in\End(X)\times\End(X).
\]
Whenever a multilinear polynomial $\Phi_{t,r}$ is evaluated on elements of $X$ together with one occurrence of $f$, we interpret a product involving $f$ by
\[
fx=f\ast x
\qquad\text{and}\qquad
xf=x\ast f.
\]
Since $f$ occurs exactly once, this determines an element of~$X$.

The \emph{external weak actor} of $X$ is the vector subspace
\[
\E(X)\subseteq\End(X)\times\End(X)
\]
consisting of all pairs $f=(f\ast-,-\ast f)$ such that
\[
\Phi_{t,r}(\alpha_1,\ldots,\alpha_t)=0
\]
for any defining identity $\Phi_{t,r}$ of $\V$ and any choice of $\alpha_j=f$ and $\alpha_k\in X$ for $k\neq j$.

The chosen $\lambda/\mu$-rules induce a bilinear map
\[
m_X\colon\E(X)\times\E(X)\to\End(X)\times\End(X).
\]
More precisely, for any $f=(f\ast-,-\ast f), \;
g=(g\ast-,-\ast g)\in\E(X)$, we define
\[
m_X(f,g)=(h\ast-,-\ast h)
\]
by
\begin{align*}
x\ast h
={}&\lambda_1(x\ast f)\ast g+\lambda_2(f\ast x)\ast g+\lambda_3g\ast(x\ast f)+\lambda_4g\ast(f\ast x)\\
&+\lambda_5(x\ast g)\ast f+\lambda_6(g\ast x)\ast f+\lambda_7f\ast(x\ast g)+\lambda_8f\ast(g\ast x)
\end{align*}
and
\begin{align*}
h\ast x
={}&\mu_1(x\ast f)\ast g+\mu_2(f\ast x)\ast g+\mu_3g\ast(x\ast f)+\mu_4g\ast(f\ast x)\\
&+\mu_5(x\ast g)\ast f+\mu_6(g\ast x)\ast f+\mu_7f\ast(x\ast g)+\mu_8f\ast(g\ast x)
\end{align*}
for every $x\in X$.

In general, $m_X(f,g)$ need not belong to $\E(X)$. Let
\[
\widetilde m_X\colon\E(X)\otimes\E(X)\to\End(X)\times\End(X)
\]
be the linearisation of $m_X$. Then
\[
D_X\coloneqq\widetilde m_X^{-1}\bigl(\E(X)\bigr)
\]
is a vector subspace of $\E(X)\otimes\E(X)$, and
\[
\mu_X\coloneqq\left.\widetilde m_X\right|_{D_X}\colon D_X\to\E(X)
\]
is a linear map. Thus, $(\E(X),D_X,\mu_X)$ is a partial algebra, which we still denote by $\E(X)$. Its partial multiplication is denoted by $\langle-,-\rangle$: for $f,g\in\E(X)$, the product $\langle f,g\rangle$ is defined if and only if $f\otimes g\in D_X$, or, equivalently, if and only if $m_X(f,g)\in\E(X)$, and in this case
\[
\langle f,g\rangle=\mu_X(f\otimes g)=m_X(f,g).
\]

We further observe that the split extension~\eqref{SplitExternalActor} determines a homomorphism of partial algebras
\[
\widetilde U(B)\to\E(X)\colon b\mapsto(b\ast-,-\ast b),
\]
where $(b\ast-,-\ast b)$ is induced by the associated derived action. The following result shows that these maps define a natural monomorphism of functors.

\begin{theorem}\cite{WRAAlg,SurveyActions}\label{ThmExternalActorOneOperation}
Let $\V$ be an action accessible, operadic variety of non-associative algebras, and fix a choice of $\lambda/\mu$-rules. For every object $X$ of $\V$, there exists a natural monomorphism of functors
\[
\tau\colon\SplExt(-,X)\rightarrowtail\Hom_{\PAlg}(\widetilde U(-),\E(X)).
\]
For every object $B$ of $\V$, the component~$\tau_B$ sends a split extension of $B$ by $X$ to the homomorphism of partial algebras
\[
\widetilde U(B)\to\E(X)\colon b\mapsto(b\ast-,-\ast b)
\]
induced by the corresponding derived action of $B$ on $X$. Moreover, if $\E(X)$ is total and the resulting algebra, with multiplication $\langle-,-\rangle$, is an object of $\V$, then $(\E(X),\tau)$ is a weak representation of $\SplExt(-,X)$.
\end{theorem}

We refer to~\cite[Theorem~5.14]{SurveyActions} for a proof of \Cref{ThmExternalActorOneOperation} in the present formulation of partial algebras.

The pair $(\E(X),\tau)$ is called an \emph{external weak representation} of $\SplExt(-,X)$. If~$\tau$ is a natural isomorphism, then $\E(X)$ is called \emph{external actor} of~$X$.

\begin{remark}\label{RemarkExternalComACom}
Suppose that $\V$ is a subvariety of $\Com$ or $\ACom$. Then every element $(f\ast-,-\ast f)$ of $\E(X)$ satisfies
\[
x\ast f=\varepsilon(f\ast x)
\]
for every $x\in X$, where $\varepsilon=1$ in the commutative case and $\varepsilon=-1$ in the anti-commutative case. Hence, $\E(X)$ is isomorphic, as a vector space, to
\[
\left\{
f\in\End(X)
\mid
\Phi_{t,r}(x_1,\ldots,x_{j-1},f,x_{j+1},\ldots,x_t)=0
\right\},
\]
where the condition is required for every defining identity
$\Phi_{t,r}(x_1,\ldots,x_t)=0$ of~$\V$, every $1\leq j\leq t$ and every choice of the other variables in $X$, and products involving~$f$ are evaluated using $fx=f(x)$ and $xf=\varepsilon f(x)$. The isomorphism is given by
\[
f\longmapsto(f,\varepsilon f).
\]
In this case, the $\lambda/\mu$-rules reduce to an identity of the form
\[
x(yz)=\alpha(xy)z+\beta(xz)y
\]
for some $\alpha,\beta\in\F$, and one has
\[
(yz)x=\varepsilon (\alpha z(yx)+ \beta y(zx)).
\]
Thus, under this identification, the product $\langle f,g\rangle$ of $f,g\in\E(X)$ is defined if and only if $\varepsilon\bigl(\alpha\,(g\circ f)+\beta\,(f\circ g)\bigr)$ belongs to $\E(X)$, and in this case
\[
\langle f,g\rangle=\varepsilon\bigl(\alpha\,(g\circ f)+\beta\,(f\circ g)\bigr).
\]
\end{remark}

We conclude this section with some examples of external weak actors.

\begin{example}\label{ExampleExternalAssoc}
Let $\V=\Assoc$. The standard choice of the $\lambda/\mu$-rules is given by
\[
\lambda_1=\mu_8=1,
\]
with all the other coefficients equal to zero. In this case,
\[
\E(X)\cong\Bim(X),
\]
where
\[
\Bim(X)=
\left\{
(f\ast-,-\ast f)\in\End(X)^2
\ \middle|\
\begin{aligned}
f\ast(xy)&=(f\ast x)y,\\
(xy)\ast f&=x(y\ast f),\\
x(f\ast y)&=(x\ast f)y
\end{aligned} \quad \forall x,y\in X
\right\}
\]
is the associative algebra of \emph{bimultipliers} of~$X$; see~\cite{MacLane58}. Its multiplication is given~by
\[
\langle f,g\rangle\ast x=f\ast(g\ast x),
\qquad
x\ast\langle f,g\rangle=(x\ast f)\ast g.
\]
Thus, $\E(X)$ is an associative algebra, so $\Assoc$ is weakly action representable.
\end{example}

\begin{example}\label{ExampleExternalCAssoc}
Let $\V=\CAssoc$, and let $X$ be a commutative associative algebra. In the notation of \Cref{RemarkExternalComACom}, choose the $\lambda/\mu$-rules with $\alpha=0$ and $\beta=1$. Then $\langle f,g\rangle=f\circ g$, and $\E(X)$ is isomorphic to the associative algebra
\[
\operatorname{M}(X)
=
\left\{
f\in\End(X)
\mid
f(xy)=f(x)y, \quad \forall x,y\in X
\right\}
\]
of \emph{multipliers} of~$X$~\cite{Casas}, whose multiplication is given by the usual composition of functions. Thus, $\operatorname{M}(X)$ is an associative algebra. Moreover, there is a natural isomorphism
\[
\SplExt_{\CAssoc}(-,X)
\cong
\Hom_{\Assoc}\bigl(\overline{U}(-),\operatorname{M}(X)\bigr),
\]
where $\overline{U}\colon\CAssoc\to\Assoc$ is the forgetful functor; see~\cite[Proof of Theorem~2.6]{BJK2}. Hence, $\E(X)\cong\operatorname{M}(X)$ is an external actor of~$X$. Although $\M(X)$ is the natural candidate for a weakly representing object, it need not be commutative and hence need not be an object of $\CAssoc$.
\end{example}

\begin{example}\label{ExampleExternalLie}
Let $\bchar(\F) \neq 2$ and consider the variety $\Lie$ of Lie algebras. Then, $\E(X)$ is isomorphic to the Lie algebra
\[
\Der(X)=
\left\{
d\in\End(X)
\mid
d([x,y])=[d(x),y]+[x,d(y)], \; \forall x,y\in X
\right\}
\]
of derivations of~$X$, endowed with the Lie bracket
\[
[d,d']=d \circ d' - d' \circ d.
\]
In this case, the external weak representation is a natural isomorphism, and $\Der(X)$ is the actor of~$X$.
\end{example}

\begin{example}\label{ExampleExternalLeib}
Let $\V=\Leib$. Consider the presentation of the right Leibniz identity given by
\[
[x,[y,z]]=[[x,y],z]-[[x,z],y]
\]
and
\[
[[y,z],x]=[[y,x],z]-[y,[x,z]].
\]
The corresponding external weak actor is isomorphic to the Leibniz algebra
\[
\Bider(X)=
\left\{
(d,D)\in\End(X)^2
\ \middle|\
\begin{aligned}
d([x,y])&=[d(x),y]+[x,d(y)],\\
D([x,y])&=[D(x),y]-[D(y),x],\\
[x,d(y)]&=[x,D(y)]
\end{aligned}
\quad \forall x,y\in X
\right\}
\]
of \emph{biderivations} of~$X$~\cite{loday}, endowed with the bracket
\[
[(d,D),(d',D')]
=
\bigl(d\circ d'-d'\circ d,\,
D\circ d'-d'\circ D\bigr).
\]
The isomorphism is given by
\[
\Bider(X)\to\E(X)\colon(d,D)\mapsto(D,-d);
\]
see~\cite{Casas,CigoliManciniMetere}. In particular, $\E(X)$ is a Leibniz algebra, so $\Leib$ is weakly action representable.
\end{example}

\section{Varieties of algebras with \texorpdfstring{$n$}{n} bilinear multiplications}\label{SecNOperations}

An \emph{algebra with $n$ bilinear multiplications} over $\F$ is an $\F$-vector space $X$ equipped with $n$ bilinear maps
\[
\ast_i \colon X\times X\to X \colon (x,y) \mapsto x \ast_i y, \qquad i=1,\ldots,n.
\]
No associativity or compatibility condition between the $n$ multiplications is assumed. We denote such an algebra by $(X,\ast_1,\ldots,\ast_n)$.

We recall the following definitions.

\begin{definition}
A homomorphism of algebras with $n$ bilinear multiplications
\[
f\colon(X,\ast_1,\ldots,\ast_n)\to(Y,\ast_1',\ldots,\ast_n')
\]
is a linear map satisfying
\[
f(x\ast_i y)=f(x)\ast_i' f(y)
\]
for every $i=1,\ldots,n$ and $x,y\in X$.
\end{definition}

We denote by $\Alg^n$ the category of $\F$-algebras with $n$ bilinear multiplications and their homomorphisms.

\begin{definition}
Let $(X,\ast_1,\ldots,\ast_n)$ be an object of $\Alg^n$. A vector subspace $Y\subseteq X$ is a \emph{subalgebra} of $X$ if
\[
y_1\ast_i y_2\in Y
\]
for any $i=1,\ldots,n$ and $y_1,y_2\in Y$. A subalgebra $I\subseteq X$ is an \emph{ideal} of $X$ if
\[
x\ast_i j,\ j\ast_i x\in I
\]
for any $i=1,\ldots,n$, $x\in X$ and $j\in I$.
\end{definition}

\begin{remark}
From a categorical point of view, subalgebras are precisely subobjects, while ideals are precisely normal subobjects, i.e., kernels of morphisms. Indeed, given a morphism $f\colon X\to Y$ in $\Alg^n$, its kernel $\Ker f$ is an ideal of $X$. Conversely, every ideal $I$ of $X$ is the kernel of the canonical projection $\pi\colon X\to X/I$.
\end{remark}

Let $S$ be a set of variables. We denote by $\mathsf{W}(S)$ the smallest set satisfying the following conditions:
\begin{enumerate}
\item every element of $S$ belongs to $\mathsf{W}(S)$;
\item if $u,v\in\mathsf{W}(S)$, then $u\ast_i v\in\mathsf{W}(S)$, for any $i=1,\ldots,n$.
\end{enumerate}
Thus, the elements of $\mathsf{W}(S)$ are the non-associative words in the elements of $S$ formed using the $n$ formal multiplications $\ast_1,\ldots,\ast_n$.

We denote by $\F\{S\}_{\ast_1,\ldots,\ast_n}$
the $\F$-vector space with basis $\mathsf{W}(S)$, endowed with the $n$ bilinear multiplications obtained by extending the formal multiplications $\ast_1,\ldots,\ast_n$ linearly. The map
\[
S\longmapsto\F\{S\}_{\ast_1,\ldots,\ast_n}
\]
defines a functor
\[
\F\{-\}_{\ast_1,\ldots,\ast_n}\colon\Set\to\Alg^n
\]
which is left adjoint to the forgetful functor
\[
U_n\colon\Alg^n\to\Set.
\]
Equivalently, for every algebra $(X,\ast_1,\ldots,\ast_n)$ and every map $f\colon S\to U_n(X)$, there exists a unique algebra homomorphism
\[
\overline{f}\colon\F\{S\}_{\ast_1,\ldots,\ast_n}\to X
\]
whose restriction to $S$ is~$f$. Thus, $\F\{S\}_{\ast_1,\ldots,\ast_n}$ is the free algebra with $n$ bilinear multiplications generated by~$S$.

An element of $\F\{S\}_{\ast_1,\ldots,\ast_n}$ is called a \emph{non-associative polynomial with $n$ multiplications}. An \emph{identity} of an algebra $(X,\ast_1,\ldots,\ast_n)$ is a polynomial
\[
\Phi=\Phi(x_1,\ldots,x_t)
\in
\F\{x_1,\ldots,x_t\}_{\ast_1,\ldots,\ast_n}
\]
such that
\[
\Phi(a_1,\ldots,a_t)=0
\]
for every $a_1,\ldots,a_t\in X$.

Given a set $I$ of identities, the \emph{variety of algebras with $n$ bilinear multiplications} determined by $I$ is the full subcategory $\V$ of $\Alg^n$ consisting of all algebras satisfying every identity in~$I$. Any such variety is a semi-abelian category.

For a monomial $\varphi=\varphi(x_1,\ldots,x_t)$, its \emph{type} is the $t$-tuple
\[
(k_1,\ldots,k_t)\in\mathbb{N}^t,
\]
where $k_i$ is the number of occurrences of $x_i$ in $\varphi$. Its \emph{degree} is
\[
k_1+\cdots+k_t.
\]
A polynomial is \emph{homogeneous} if all its monomials have the same type, and it is \emph{multilinear} if all its monomials have type $(1,\ldots,1)$.

\begin{definition}
A variety of algebras with $n$ bilinear multiplications is \emph{homogeneous} if it can be determined by a set of homogeneous identities, and it is \emph{operadic} if it can be determined by a set of multilinear identities.
\end{definition}

\begin{remark}
For every $i=1,\ldots,n$, let $\V_{\ast_i}$ denote the variety of non-associative algebras determined by all the identities of $\V$ involving only the multiplication $\ast_i$. There is a canonical forgetful functor
\[
U_{\ast_i}\colon\V\to\V_{\ast_i},
\qquad
(X,\ast_1,\ldots,\ast_n)\longmapsto(X,\ast_i).
\]
\end{remark}

We now recall some examples of varieties of algebras with $n$ bilinear multiplications. Whenever one of the multiplications is a Lie bracket, we assume that $\bchar(\F)\neq2$, so that the corresponding variety is operadic.

\begin{example}
The variety $\AbAlg^n$ of \emph{abelian algebras with $n$ multiplications} $\ast_1,\ldots,\ast_n$ is determined by the identities
\[
x \ast_1 y=\ldots=x \ast_n y=0.
\]
Clearly, $\AbAlg^n$ is isomorphic to the category $\Vect$ of vector spaces over $\F$ and is therefore an abelian category.
\end{example}

\begin{example}
The variety $\Lie^n$ of \emph{algebras with $n$ Lie brackets}
\[
[-,-]_1,\ldots,[-,-]_n
\]
is determined by the identities $[x,x]_i=0$ and
\[
[x,[y,z]_i]_i+[y,[z,x]_i]_i+[z,[x,y]_i]_i=0
\]
for any $i=1,\ldots,n$. No compatibility condition between the brackets is required.
\end{example}

We now focus on the case $n=2$.

\begin{example}\label{ExamplePoissonVariety}
A \emph{Poisson algebra} is a vector space $X$ equipped with an associative multiplication $\cdot$
and a Lie bracket $[-,-]$ satisfying the Poisson identity
\[
[x,y\cdot z]=[x,y]\cdot z+y\cdot[x,z]
\]
for every $x,y,z\in X$. We denote by $\Pois$ the variety of Poisson algebras. The subvariety of $\Pois$ determined by the additional identity
\[
x\cdot y=y\cdot x
\]
is denoted by $\CPois$ and is called the variety of \emph{commutative Poisson algebras}.

We observe that every associative algebra becomes a Poisson algebra with the standard Lie bracket $[x,y]=x \cdot y - y \cdot x$.
\end{example}

\begin{example}\cite{CasasDatuashviliLadraNPLR}\label{ExampleNoncommutativePoisson}
A \emph{left non-commutative Poisson algebra} is a vector space $X$ equipped with an associative multiplication $\cdot$ and a right Leibniz bracket $[-,-]$ satisfying
\[
[x\cdot y,z]=x\cdot[y,z]+[x,z]\cdot y
\]
for every $x,y,z\in X$. These algebras were originally introduced under the name \emph{non-commutative Leibniz--Poisson algebras} in \cite{CasasDatuashviliNLP}. A \emph{right non-commutative Poisson algebra} is defined by requiring instead
\[
[x,y\cdot z]=y\cdot[x,z]+[x,y]\cdot z.
\]
A \emph{left-right non-commutative Poisson algebra} is an associative algebra endowed with a right Leibniz bracket satisfying both identities. We denote the corresponding varieties by $\NP^{l}$, $\NP^{r}$, and $\NP^{lr}$, respectively.

Clearly, every Poisson algebra is a left-right non-commutative Poisson algebra. Moreover, every Leibniz algebra becomes a left-right non-commutative Poisson algebra by taking the associative multiplication to be trivial.
\end{example}

\begin{example}\cite{LodayDialgebras}\label{ExampleAssociativeDialgebras}
An \emph{associative dialgebra} is a vector space $X$ equipped with two associative multiplications
\[
\dashv,\vdash\colon X\times X\to X
\]
satisfying
\begin{align*}
x\dashv (y\dashv z)&=x\dashv(y\vdash z),\\
(x\vdash y)\dashv z&=x\vdash(y\dashv z),\\
(x\dashv y)\vdash z&=(x\vdash y)\vdash z
\end{align*}
for every $x,y,z\in X$. We denote the corresponding variety by $\Dias$.
\end{example}

\begin{example}\cite{ZhangBaiGuo}\label{ExampleMatchingDialgebras}
A \emph{matching dialgebra} is a vector space $X$ equipped with two associative multiplications $\ast$ and $\bullet$ satisfying the matching associativity
\[
(x\ast y)\bullet z=x\ast(y\bullet z),
\]
\[
(x\bullet y)\ast z=x\bullet(y\ast z)
\]
for every $x,y,z\in X$. Equivalently,
\[
(x\diamond y)\circ z=x\diamond(y\circ z)
\]
for every $\diamond,\circ\in\{\ast,\bullet\}$. We denote the corresponding variety by $\Match$.
\end{example}

\begin{remark}
Every associative algebra determines both an associative dialgebra and a matching dialgebra by taking the two multiplications to coincide.

We further observe that every associative dialgebra defines a left non-commutative Poisson algebra with respect to the multiplications
\[
x\cdot y=x\vdash y
\qquad\text{and}\qquad
[x,y]=x\dashv y-y\vdash x,
\]
see~\cite[Example~2.1.2]{CasasDatuashviliNLP}. In general, it does not define an object of $\NP^{lr}$.
\end{remark}

\begin{example}\cite{Strohmayer}\label{ExampleCompatibleLie}
A \emph{compatible Lie algebra} is a vector space $X$ equipped with two Lie brackets $[-,-]$ and $\{-,-\}$ satisfying the mixed Jacobi identity
\begin{align*}
\{[x,y],z\}+\{[y,z],x\}+\{[z,x],y\}
&+[\{x,y\},z]\\
&+[\{y,z\},x]+[\{z,x\},y]=0
\end{align*}
for every $x,y,z\in X$. Equivalently, every linear combination
\[
\lambda[-,-]+\mu\{-,-\},
\qquad
\lambda,\mu\in\F,
\]
is a Lie bracket on~$X$. We denote the corresponding variety by $\CLie$.
\end{example}

\begin{example}\cite{GaoGuoZhang}\label{ExampleMatchingLie}
A \emph{matching Lie algebra} is a vector space $X$ endowed with two Lie brackets $[-,-]$ and $\{-,-\}$ satisfying the matching Jacobi identities
\begin{align*}
[x,\{y,z\}]+\{y,[z,x]\}+\{z,[x,y]\}&=0,\\
\{x,[y,z]\}+[y,\{z,x\}]+[z,\{x,y\}]&=0
\end{align*}
for all $x,y,z\in X$. We denote the corresponding variety by $\MatchLie$.

Every matching Lie algebra is a compatible Lie algebra. Indeed, the sum of the two matching Jacobi identities gives the mixed Jacobi identity defining $\CLie$. Thus, $\MatchLie$ is a subvariety of $\CLie$.
\end{example}

\begin{example}\cite{VallettePostLie}
A \emph{post-Lie algebra} over $\F$ is a vector space $X$ endowed with a Lie bracket $[-,-]$ and a bilinear multiplication $\cdot$ satisfying
\[
x\cdot[y,z]=[x\cdot y,z]+[y,x\cdot z]
\]
and
\[
[x,y]\cdot z=a(x,y,z)-a(y,x,z),
\]
for all $x,y,z\in X$, where
\[
a(x,y,z)=x\cdot(y\cdot z)-(x\cdot y)\cdot z
\]
denotes the associator on $(X,\cdot)$. We denote by $\PostLie$ the variety of post-Lie algebras over $\F$. 

Every pre-Lie algebra $(X,\cdot)$ can be seen as a post-Lie algebra with trivial Lie bracket. Under this identification, $\PreLie$ is the subvariety of $\PostLie$ determined by the identity $[x,y]=0$. Since $\PreLie$ is not action accessible, it follows that $\PostLie$ is not action accessible either; see~\cite[Proposition~3.5]{AlbanoStefanelli}.
\end{example}

\subsection{Action accessibility and categories of interest}\label{SecNCategoriesInterest}

Let $n\geq1$, and let $\V$ be a homogeneous variety of $\F$-algebras endowed with bilinear multiplications $\ast_1,\ldots,\ast_n$. Since the underlying additive group of every object of $\V$ is abelian, Axiom~1 of \Cref{DefCategoryInterest} is satisfied. Moreover, the sets of non-group unary and binary operations of $\V$ are
\[
\Omega_1'=\{\omega_\lambda\mid\lambda\in\F\} \quad \text{ and } \quad
\Omega_2'=\{\ast_1,\ast_1^{\op},\ldots,\ast_n,\ast_n^{\op}\},
\]
where $\omega_\lambda(x)=\lambda x$ and
\[
x\ast_i^{\op}y=y\ast_i x.
\]
Thus, $\V$ is an Orzech category of interest if and only if Axiom~2 is satisfied for every $(\circ,\diamond)\in\Omega_2'\times\Omega_2'$. For such a pair, consider the formal expressions
\[
x_1(x_2x_3),\quad x_1(x_3x_2),\quad (x_2x_3)x_1,\quad (x_3x_2)x_1,
\]
\[
x_2(x_1x_3),\quad x_2(x_3x_1),\quad (x_1x_3)x_2,\quad (x_3x_1)x_2,
\]
where each juxtaposition denotes an operation of $\Omega_2'$. Axiom~2 requires the existence of a linear combination $W_{\circ,\diamond}(x_1,x_2;x_3)$ of the monomials above such that
\[
(x_1\circ x_2)\diamond x_3
=
W_{\circ,\diamond}(x_1,x_2;x_3)
\]
is an identity of $\V$.

\begin{definition}
A family of identities
\[
(x_1\circ x_2)\diamond x_3
=
W_{\circ,\diamond}(x_1,x_2;x_3),
\qquad
(\circ,\diamond)\in\Omega_2'\times\Omega_2',
\]
satisfying Axiom~2 is called a \emph{choice of identities for Axiom~2} for $\V$.
\end{definition}

A choice of identities for Axiom~2 need not be unique. Throughout the article, whenever $\V$ is an Orzech category of interest, we fix one such choice.

\begin{remark}
When $n=1$, the identities above reduce to the usual $\lambda/\mu$-rules for varieties of algebras with one bilinear multiplication.
\end{remark}

The following result is a generalisation of \Cref{ThmLambdaMu} to varieties of algebras with $n$ bilinear multiplications.

\begin{theorem}\label{ThmActionAccessibilityNOperations}
Let $\V$ be a homogeneous variety of algebras with $n$ bilinear multiplications $\ast_1,\ldots,\ast_n$. The following conditions are equivalent:
\begin{tfae}
\item $\V$ is an Orzech category of interest.
\item $\V$ is action accessible.
\item For every $(\circ,\diamond)\in\Omega_2'\times\Omega_2'$, there exists a linear combination
\[
W_{\circ,\diamond}(x_1,x_2;x_3)
\]
of the degree-three monomials above such that
\[
(x_1\mathbin{\circ}x_2)\mathbin{\diamond}x_3
=
W_{\circ,\diamond}(x_1,x_2;x_3)
\]
is an identity of $\V$.
\end{tfae}
\end{theorem}

\begin{proof}
We only need to prove $(2)\Rightarrow(3)$. Indeed, the implication $(1)\Rightarrow(2)$ follows from~\cite[Theorem~4.7]{Montoli}, while $(3)\Rightarrow(1)$ follows from the discussion above.

Let $\V$ be action accessible. We follow the argument used in the proof of~\cite[Theorem~3.4.1]{Coco_Thesis}. By~\cite[Lemma~2.9]{CigoliMantovaniMetere}, if $I$ and $J$ are subalgebras of an algebra $X$ of~$\V$, with $I$ an ideal of $X$, then
\begin{equation}\label{EqCommutatorNormalClosure}
[I,J]=0
\quad\Longrightarrow\quad
[I,\operatorname{cl}_X(J)]=0,
\end{equation}
where $[-,-]$ denotes the Higgins commutator and $\operatorname{cl}_X(J)$ is the normal closure of~$J$ in $X$, that is, the smallest ideal of $X$ containing $J$.

Let $W$ be the free algebra in $\V$ generated by $x,y,z$. Consider the ideal $I=(x)$ of $W$ generated by $x$ and the subalgebra $J$ of $W$ generated by $y$. Let $N$ be the ideal of $W$ generated by
\[
\{i\mathbin{\circ}j
\mid i\in I,\ j\in J,\ \circ\in \Omega_2'\}
\]
and let
\[
q\colon W\to Q\coloneqq W/N
\]
be the canonical projection.

In $Q$, all products between elements of $q(I)$ and elements of $q(J)$ vanish for every multiplication in $\Omega_2'$ and in both orders. Hence,
\[
[q(I),q(J)]=0.
\]
The image $q(I)$ of $I$ is an ideal of $Q$, while the image $q(J)$ of $J$ is a subalgebra. Thus, \cite[Lemma~2.9]{CigoliMantovaniMetere} implies that
\[
[q(I),\operatorname{cl}_Q(q(J))]=0.
\]
Since quotients preserve normal closures and Higgins commutators, we have
\[
0=[q(I),\operatorname{cl}_Q(q(J))]=\bigl[q(I),q(\operatorname{cl}_W(J))\bigr]=q\bigl([I,\operatorname{cl}_W(J)]\bigr).
\]
For every $\circ,\diamond\in\{\ast_1,\ldots,\ast_n\}$, the elements $y\mathbin{\circ}z,z\mathbin{\circ}y$ belong to $\operatorname{cl}_W(J)$. Since $x\in I$, it follows that
\[
x\mathbin{\diamond}(y\mathbin{\circ}z), \,(y\mathbin{\circ}z)\mathbin{\diamond}x
\in[I,\operatorname{cl}_W(J)].
\]
Their images under $q$ are therefore zero, and hence
\[
x\mathbin{\diamond}(y\mathbin{\circ}z),\,(y\mathbin{\circ}z)\mathbin{\diamond}x
\in N.
\]
Since $\V$ is homogeneous, the free algebra $W$ is graded by type. The ideal $I$ and the subalgebra $J$ are graded, since they are generated by the type-homogeneous elements $x$ and $y$, respectively. Hence, by bilinearity, $N$ is generated as an ideal by the products $i\mathbin{\circ}j$ with $i\in I$ and $j\in J$ type-homogeneous, so that $N$ is a graded ideal. Thus, the two elements above, which have type $(1,1,1)$, belong to the component of~$N$ of type $(1,1,1)$.

Each of these generators $i\mathbin{\circ}j$ involves $x$, through $i$, and $y$, through $j$. Hence, the terms contributing to the component of type $(1,1,1)$ arise in two ways. First, a generator $i\mathbin{\circ}j$ may itself have type $(1,1,1)$. In this case, $j$ is a scalar multiple of $y$, while $i$ has type $(1,0,1)$ and is a linear combination of the monomials
\[
x\mathbin{\diamond}z,
\qquad
z\mathbin{\diamond}x,
\qquad
\diamond\in\Omega_2'.
\]
Second, a term may be obtained by multiplying a generator of type $(1,1,0)$ once, on either side, by $z$. Such a generator is a linear combination of the monomials
\[
x\mathbin{\diamond}y,
\qquad
y\mathbin{\diamond}x,
\qquad
\diamond\in\Omega_2'.
\]
Therefore, the component of $N$ of type $(1,1,1)$ is spanned by monomials of the forms
\[
(xy)z,\quad
(yx)z,\quad
z(xy),\quad
z(yx),
\]
\[
(xz)y,\quad
(zx)y,\quad
y(xz),\quad
y(zx),
\]
where each juxtaposition may independently denote an
operation in $\Omega_2'$.

It follows that, for any $(\circ,\diamond)\in
\Omega_2' \times \Omega_2'$, there exists a linear combination
\[
W_{\circ,\diamond}(x_1,x_2;x_3)
\]
of the monomials above such that, after taking
\[
y=x_1,\qquad z=x_2,\qquad x=x_3,
\]
we obtain the identity
\[
(x_1\mathbin{\circ}x_2)\mathbin{\diamond}x_3
=
W_{\circ,\diamond}(x_1,x_2;x_3).
\]
This proves~$(3)$.
\end{proof}

In view of \Cref{ThmActionAccessibilityNOperations}, we have the following.

\begin{proposition}\label{PropExamplesNOperations}
The varieties $\AbAlg^n$, $\Pois$, $\CPois$, $\NP^{lr}$, $\Dias$, $\Match$, and $\MatchLie$ are action accessible. On the other hand, the varieties $\NP^{l}$, $\NP^{r}$, $\CLie$, and $\PostLie$, as well as the varieties $\Lie^n$ for every $n\geq2$, are not action accessible.
\end{proposition}

\begin{proof}
For $\AbAlg^n$, $\Pois$, $\CPois$, $\Dias$, $\Match$, and $\MatchLie$, their defining identities provide a choice of identities for Axiom~2; explicit choices are given in \Cref{SecPoisson,SecExternalDias,SecExternalMatching,SecExternalMatchingLie}. The case of $\NP^{lr}$ is treated in~\cite[Example~2.5]{CasasDatuashviliLadraNPLR}, where it is also shown that $\NP^l$ and $\NP^r$ do not satisfy Axiom~2. Hence, by \Cref{ThmActionAccessibilityNOperations}, they are not action accessible. The variety $\PostLie$ is not action accessible by~\cite[Proposition~3.5]{AlbanoStefanelli}. 

It remains to consider $\CLie$ and $\Lie^n$. Let $X$ be the vector space with basis
\[
\{a,b,c,u,v,w\}.
\]
We define two bilinear operations $[-,-]$ and $\{-,-\}$ by
\[
[a,b]=-[b,a]=u,
\qquad
[v,c]=-[c,v]=-w,
\]
and
\[
\{a,b\}=-\{b,a\}=v,
\qquad
\{u,c\}=-\{c,u\}=w,
\]
with all the remaining brackets between elements of the basis being zeros. It is immediate to check that both $(X,[-,-])$ and $(X,\{-,-\})$ are $2$-nilpotent algebras, and hence they are Lie algebras. Moreover, the only non-zero mixed multiplications occurring in the mixed Jacobi identity are
\[
\{[a,b],c\}=w
\qquad\text{and}\qquad
[\{a,b\},c]=-w.
\]
Hence, $(X,[-,-],\{-,-\})$ is a compatible Lie algebra.

Now consider Axiom~2 for the monomial $\{[a,b],c\}$. Every monomial allowed on its right-hand side vanishes, since
\[
[a,c]=[b,c]=\{a,c\}=\{b,c\}=0.
\]
On the other hand,
\[
\{[a,b],c\}=w\neq0.
\]
Therefore, $\CLie$ does not satisfy Axiom~2 and, by \Cref{ThmActionAccessibilityNOperations}, is not action accessible. Finally, let $n\geq2$ and see $X$ as an object of $\Lie^n$ by setting
\[
[-,-]_1=[-,-],
\qquad
[-,-]_2=\{-,-\},
\qquad
[-,-]_i=0
\quad\text{for }3\leq i\leq n .
\]
Since $[a,c]_i=[b,c]_i=0$ for every $i=1,\ldots,n$, every monomial allowed on the right-hand side of Axiom~2 for $[[a,b]_1,c]_2$ vanishes, whereas $[[a,b]_1,c]_2=\{[a,b],c\}=w\neq0$. Hence, for every $n\geq2$, the variety $\Lie^n$ does not satisfy Axiom~2 and, by \Cref{ThmActionAccessibilityNOperations}, it is not action accessible.
\end{proof}

The following example shows that action accessibility of a variety with $n$ bilinear multiplications $\ast_1,\ldots,\ast_n$ does not imply action accessibility of the varieties of non-associative algebras $\V_{\ast_i}$.

\begin{example}\label{ExampleAccessibleReductNotAccessible}
Assume that $\F$ is infinite and $\bchar(\F)\neq2$. Let $\V$ be the variety of algebras with two bilinear multiplications $\ast$ and $\bullet$ satisfying
\[
x\ast y=y\ast x,
\qquad
x\bullet y=-y\bullet x
\]
and
\begin{align}\label{id_pol}
\begin{split}
(x\ast y)\ast z
+(x\bullet y)\ast z
&+(x\ast y)\bullet z
+(x\bullet y)\bullet z
\\&=x\ast(y\ast z)
+x\ast(y\bullet z)
+x\bullet(y\ast z)
+x\bullet(y\bullet z).
\end{split}
\end{align}
For every object $X$ of $\V$, we define the bilinear multiplication
\[
xy=x\ast y+x\bullet y.
\]
The identity \eqref{id_pol} is precisely the associativity of this multiplication. Thus, any algebra $X$ of $\V$ defines an object of $\Assoc$. Conversely, every associative algebra determines an object of $\V$ by setting
\[
x\ast y=\frac12(xy+yx),
\qquad
x\bullet y=\frac12(xy-yx).
\]
Thus, $\V$ is term-equivalent to $\Assoc$ and is therefore action accessible.

On the other hand, the variety $\V_\ast$ determined by the identities involving only $\ast$ is not action accessible. Since $\F$ is infinite, the variety $\V_\ast$ is homogeneous. Hence, if it were action accessible, \Cref{ThmLambdaMu} and the commutativity of $\ast$ would give an identity
\[
x\ast(y\ast z)
=
\alpha(x\ast y)\ast z
+
\beta(x\ast z)\ast y
\]
for some $\alpha,\beta\in\F$. Now, consider the free associative algebra generated by $x,y,z$. 
We have
\[
x\ast(y\ast z)
=
\frac14\bigl(xyz+xzy+yzx+zyx\bigr),
\]
\[
(x\ast y)\ast z
=
\frac14\bigl(xyz+yxz+zxy+zyx\bigr),
\]
and
\[
(x\ast z)\ast y
=
\frac14\bigl(xzy+zxy+yxz+yzx\bigr).
\]
Since the associative words of degree~$3$ are linearly independent, the identity
\[
x\ast(y\ast z)
=
\alpha(x\ast y)\ast z
+
\beta(x\ast z)\ast y
\]
would imply $\alpha=1$ by comparing the coefficients of $xyz$, $\beta=1$ by comparing those of $xzy$, and $\alpha+\beta=0$ by comparing those of $yxz$. This is a contradiction since $\bchar(\F)\neq2$. Therefore, $\V_\ast$ does not satisfy the $\lambda/\mu$-rules and hence is not action accessible.
\end{example}

\subsection{Derived actions of algebras with $n$ multiplications}\label{SecDerivedActionsNOperations}

Let $\V$ be a subvariety of $\Alg^n$, and consider a split extension in $\V$
\begin{equation}\label{SplitExtensionNOperations}
\begin{tikzcd}
X \arrow[r, "k"]
&
A \arrow[r, shift left, "\alpha"]
&
B. \arrow[l, shift left, "\beta"]
\end{tikzcd}
\end{equation}
The associated derived action of $B$ on $X$ consists of the $2n$ bilinear maps
\[
l_{\ast_i}\colon B\times X\to X,
\qquad
r_{\ast_i}\colon X\times B\to X,
\qquad
i=1,\ldots,n,
\]
defined by
\[
l_{\ast_i}(b,x)=b\ast_i x
=
k^{-1}\bigl(\beta(b)\ast_i k(x)\bigr),
\qquad
r_{\ast_i}(x,b)=x\ast_i b
=
k^{-1}\bigl(k(x)\ast_i\beta(b)\bigr).
\]
These expressions are well defined because the elements on the right-hand sides belong to $\ker(\alpha)\cong X$.

Conversely, let $B$ and $X$ be objects of $\V$. The $2n$ bilinear maps
\[
l_{\ast_i} \colon B\times X\to X,\qquad r_{\ast_i} \colon X\times B\to X, \qquad i=1,\ldots,n
\]
denoted by
\[
l_{\ast_i} \colon (b,x)\longmapsto b\ast_i x,\qquad
r_{\ast_i} \colon (x,b)\longmapsto x\ast_i b,
\]
determine $n$ bilinear multiplications on the direct sum $B\oplus X \cong B \times X$ by
\begin{align}\label{SemidirectProductNOperations}
(b,x)\ast_i(b',x')
&=
\bigl(b\ast_i b',\,x\ast_i x'+b\ast_ix'+x\ast_i b'\bigr).
\end{align}
This construction gives rise to the
split extension in $\Alg^n$
\begin{equation}
\begin{tikzcd}
X \arrow[r, "\iota_2"]
&
(B\oplus X, \ast_1, \ldots, \ast_n) \arrow[r, shift left, "\pi_1"]
&
B \arrow[l, shift left, "\iota_1"],
\end{tikzcd}
\end{equation}
where $\iota_2(x)=(0,x)$, $\iota_1(b)=(b,0)$, and $\pi_1(b,x)=b$. This is a split extension in~$\V$ if and only if
$(B\oplus X,\ast_1,\ldots,\ast_n)$ is an object of $\V$, that is, if and only if it satisfies the identities defining $\V$. We therefore obtain the
following result, analogous to~\cite[Theorem~2.4]{Orzech}.

\begin{proposition}
The $2n$ bilinear maps above define a \emph{derived action} of $B$ on~$X$ in $\V$ if and only if the algebra $B\oplus X$, endowed with the $n$ multiplications~\eqref{SemidirectProductNOperations}, is an object of $\V$. In this case, it is called the \emph{semidirect product} of~$B$ and $X$ with respect to the derived action $(l_{\ast_1},r_{\ast_1},\ldots,l_{\ast_n},r_{\ast_n})$ and is denoted by $B\ltimes X$. \noproof
\end{proposition}

\begin{remark}
The construction above gives a natural bijection between derived actions of $B$ on $X$ and isomorphism classes of split extensions of $B$ by $X$; see~\cite[Theorem~3.17]{SurveyActions}.
\end{remark}

\begin{remark}\label{rem_b*}
For every split extension~\eqref{SplitExtensionNOperations} and its associated derived
action, there is an isomorphism of split extensions
\begin{equation*}
\begin{tikzcd}
X
    \arrow[r, "\iota_2"]
    \arrow[d, equal]
&
B\ltimes X
    \arrow[r, shift left, "\pi_1"]
    \arrow[d, "\theta"']
&
B
    \arrow[l, shift left, "\iota_1"]
    \arrow[d, equal]
\\
X \arrow[r, "k"]
&
A \arrow[r, shift left, "\alpha"]
&
B \arrow[l, shift left, "\beta"],
\end{tikzcd}
\end{equation*}
where
$\theta
(b,x)=\beta(b)+k(x)$.
Thus, for each $\circ\in\{\ast_1,\ldots,\ast_n\}$, when writing $b\mathbin{\circ}x$ (respectively,~$x\mathbin{\circ}b$), one may identify it with the product $(b,0)\mathbin{\circ}(0,x)$ (respectively, $(0,x)\mathbin{\circ}(b,0)$) in $B\ltimes X$.
\end{remark}

\begin{example}\label{DerivedActionMatching}
Let $B$ and $X$ be matching dialgebras. Four bilinear maps
\[
l_\ast, l_\bullet \colon B\times X\to X,\qquad r_\ast, r_\bullet \colon X\times B\to X
\]
denoted as above, define a derived action of $B$ on $X$ if and only if
\begin{align*}
(b\circ b')\diamond x&=b\circ(b'\diamond x),&
(x\circ b)\diamond b'&=x\circ(b\diamond b'),\\
(b\circ x)\diamond b'&=b\circ(x\diamond b'),&
(b\circ x)\diamond y&=b\circ(x\diamond y),\\
(x\circ b)\diamond y&=x\circ(b\diamond y),&
(x\circ y)\diamond b&=x\circ(y\diamond b),
\end{align*}
for every $b,b'\in B$, $x,y\in X$, and $\circ,\diamond\in\{\ast,\bullet\}$. These identities are obtained from
\[
(u\circ v)\diamond w=u\circ(v\diamond w)
\]
by allowing one or two of the variables to belong to $B$ and the remaining variables to belong to $X$.

We observe that when $\circ=\diamond$, we get the usual conditions defining an action in the variety $\Assoc$.
\end{example}

\begin{example}\label{DerivedActionMatchingLie}
Let $B$ and $X$ be matching Lie algebras, and write
\[
[-,-]_1=[-,-],
\qquad
[-,-]_2=\{-,-\}.
\]
A derived action of $B$ on $X$ is determined by four bilinear maps
\[
B\times X\to X\colon (b,x)\mapsto[b,x]_i,
\qquad
X\times B\to X\colon (x,b)\mapsto[x,b]_i,
\qquad i=1,2,
\]
satisfying
\[
[x,b]_i=-[b,x]_i
\]
and
\begin{align*}
[b,[x,y]_j]_i
&=[[b,x]_i,y]_j+[x,[b,y]_i]_j,\\
[[b,b']_i,x]_j
&=[b,[b',x]_j]_i-[b',[b,x]_i]_j,\\
[[b,x]_i,y]_j
&=[[b,x]_j,y]_i,\\
[[b,x]_i,b']_j
&=[[b,x]_j,b']_i,
\end{align*}
for all $b,b'\in B$, $x,y\in X$, and $i,j\in\{1,2\}$.

Let us explain the origin of the last two identities. If we apply the matching Jacobi identity to the triple $(x,b,y)$, we get
\[
[x,[b,y]_j]_i-[b,[x,y]_i]_j+[[b,x]_i,y]_j=0.
\]
On the other hand, the first identity, with $i$ and $j$ interchanged, gives
\[
[b,[x,y]_i]_j
=
[[b,x]_j,y]_i+[x,[b,y]_j]_i.
\]
Substituting the latter equality into the former, we obtain
\[
[[b,x]_i,y]_j
=
[[b,x]_j,y]_i.
\]
The last identity is obtained in a similar way, by applying the matching Jacobi identity to the triple $(x,b,b')$. 

We observe that when $i=j$, we get the usual conditions defining an action in the variety $\Lie$.
\end{example}

\begin{remark}
Given a derived action of $B$ on $X$ in $\V$, each $b\in B$ determines $2n$ linear endomorphisms of $X$
\[
\bigl(
b\ast_1-,-\ast_1 b,\ldots,
b\ast_n-,-\ast_n b
\bigr).
\]
These endomorphisms will be used in the construction of the external weak actor.
\end{remark}

\subsection{Algebras with $n$ bilinear partial operations}\label{SecPartialAlgebrasN}

We now extend the notions recalled in \Cref{SecPartialActors} to the case of $n$ operations.

\begin{definition}
An \emph{algebra with $n$ bilinear partial operations} over $\F$ is a tuple
\[
\bigl(X,(D_i,\mu_i)_{i=1}^n\bigr)=(X,D_1,\mu_1,\ldots,D_n,\mu_n),
\]
where $X$ is an $\F$-vector space and $(D_i,\mu_i)$ is a bilinear partial operation on $X$ for every $i=1,\ldots,n$, that is, $D_i\leq X\otimes X$ is a vector subspace and $\mu_i\colon D_i\to X$ is a linear map. Given $x,y\in X$ and $i\in\{1,\ldots,n\}$, we say that the product $x\ast_i y$ is \emph{defined} if $x\otimes y\in D_i$. In this case, we set
\[
x\ast_i y\coloneqq\mu_i(x\otimes y).
\]

The algebra $\bigl(X,(D_i,\mu_i)_{i=1}^n\bigr)$ is said to be \emph{total} if $D_i=X\otimes X$ for every $i=1,\ldots,n$. In this case, each $\mu_i$ is the linearisation of an ordinary bilinear multiplication on $X$.
\end{definition}

\begin{definition}
Let $\bigl(X,(D_i,\mu_i)_{i=1}^n\bigr)$ and $\bigl(X',(D_i',\mu_i')_{i=1}^n\bigr)$ be algebras with $n$ bilinear partial operations over $\F$. A \emph{homomorphism of algebras with $n$ bilinear partial operations} is a linear map $f\colon X\to X'$ such that
\[
(f\otimes f)(D_i)\subseteq D_i'
\]
and
\[
f\circ\mu_i=\mu_i'\circ\left.(f\otimes f)\right|_{D_i}
\]
for every $i=1,\ldots,n$.
\end{definition}

We denote by $\PAlg^n$ the category whose objects are algebras with $n$ bilinear partial operations over $\F$ and whose morphisms are homomorphisms of algebras with~$n$ bilinear partial operations. For $n=1$, we recover the category $\PAlg$ of \Cref{SecPartialActors}.

\begin{definition}\label{DefPartialEvaluationN}
Let $\bigl(X,(D_i,\mu_i)_{i=1}^n\bigr)$ be an algebra with $n$ bilinear partial operations, and let
\[
t\in\mathsf{W}(\{x_1,\ldots,x_r\})
\]
be a non-associative word in the $n$ formal multiplications $\ast_1,\ldots,\ast_n$. For every $\mathbf a=(a_1,\ldots,a_r)\in X^r$, we define the \emph{partial evaluation} of $t$ at $\mathbf a$ recursively as follows.

If $t=x_j$, then $t$ is defined at $\mathbf a$, and $t(\mathbf a)=a_j$.

If $t=t_1\ast_i t_2$, assume that $t_1$ and $t_2$ are defined at $\mathbf a$, with respective values
\[
t_1(\mathbf a)=b
\qquad\text{and}\qquad
t_2(\mathbf a)=c.
\]
Then $t$ is defined at $\mathbf a$ if and only if $b\otimes c\in D_i$, and in this case
\[
t(\mathbf a)=\mu_i(b\otimes c).
\]
If one of these conditions is not satisfied, then $t$ is not defined at $\mathbf a$.

Let $p\in\F\{x_1,\ldots,x_r\}_{\ast_1,\ldots,\ast_n}$ have the unique expansion
\[
p=\sum_{j=1}^{s}\lambda_jt_j,
\]
where the $t_j$ are pairwise distinct elements of $\mathsf{W}(\{x_1,\ldots,x_r\})$ and $\lambda_j\in\F\setminus\{0\}$. We say that $p$ is \emph{defined} at $\mathbf a$ if every $t_j$ is defined at $\mathbf a$. In this case, we set
\[
p(\mathbf a)=\sum_{j=1}^{s}\lambda_jt_j(\mathbf a).
\]
We allow $s=0$; thus the zero polynomial is defined at every $\mathbf a\in X^r$ and has value zero.

We say that $\bigl(X,(D_i,\mu_i)_{i=1}^n\bigr)$ \emph{satisfies the polynomial identity}
\[
p(x_1,\ldots,x_r)=0
\]
if $p(\mathbf a)=0$ for every $\mathbf a\in X^r$ at which $p$ is defined.
\end{definition}

\begin{example}
Let $x,y,z\in X$ and $i,j\in\{1,\ldots,n\}$. The monomial $x\ast_i(y\ast_j z)$ is defined at $(x,y,z)$ if and only if
\[
y\otimes z\in D_j \quad \text{ and } \quad x\otimes\mu_j(y\otimes z)\in D_i.
\]
When these conditions hold, its value is
\[
x\ast_i(y\ast_j z)
=
\mu_i\bigl(
x\otimes\mu_j(y\otimes z)
\bigr).
\]

Similarly, the monomial $(x\ast_i y)\ast_j z$ is defined if and only if
\[
x\otimes y\in D_i \quad \text{ and } \quad \mu_i(x\otimes y)\otimes z\in D_j,
\]
and in this case
\[
(x\ast_i y)\ast_j z
=
\mu_j\bigl(
\mu_i(x\otimes y)\otimes z
\bigr).
\]
\end{example}

\begin{example}
An algebra with $n$ bilinear partial operations $\bigl(X,(D_i,\mu_i)_{i=1}^n\bigr)$ satisfies the identity
\[
(x\ast_i y)\ast_j z=x\ast_i(y\ast_j z)
\]
if this equality holds for every $x,y,z\in X$ for which both sides are defined. For $n=2$, the algebras with two bilinear partial operations satisfying this identity for all $i,j\in\{1,2\}$ are the partial analogues of the matching dialgebras of \Cref{ExampleMatchingDialgebras}.
\end{example}

\begin{remark}\label{RemUTildeN}
Let $X$ be an algebra with $n$ bilinear multiplications and, for every $i=1,\ldots,n$, let
\[
\widetilde{\ast}_i\colon X\otimes X\to X
\]
be the linearisation of the multiplication $\ast_i$. Then
\[
\widetilde U(X)
=
\bigl(X,(X\otimes X,\widetilde{\ast}_i)_{i=1}^n\bigr)
\]
is a total algebra with $n$ bilinear partial operations. This construction defines a canonical fully faithful functor $\widetilde U\colon\Alg^n\to\PAlg^n$ and, by restriction, a fully faithful functor $\widetilde U\colon\V\to\PAlg^n$ for every variety of algebras $\V$ with $n$ bilinear multiplications.
\end{remark}

\begin{lemma}\label{LemMorphismFromTotalN}
Let $B$ be an algebra with $n$ bilinear multiplications, and let
\[
\bigl(X,(D_i,\mu_i)_{i=1}^n\bigr)
\]
be an algebra with $n$ bilinear partial operations. A linear map $f\colon B\to X$ defines a homomorphism of algebras with $n$ bilinear partial operations
\[
f\colon\widetilde U(B)\to\bigl(X,(D_i,\mu_i)_{i=1}^n\bigr)
\]
if and only if
\[
f(b)\otimes f(b')\in D_i
\]
and
\[
f(b\ast_i b')=\mu_i\bigl(f(b)\otimes f(b')\bigr),
\]
for every $i=1,\ldots,n$ and $b,b'\in B$.
\end{lemma}

\begin{proof}
The pure tensors $b\otimes b'$ span $B\otimes B$. Therefore,
\[
(f\otimes f)(B\otimes B)\subseteq D_i
\]
if and only if $f(b)\otimes f(b')\in D_i$ for all $b,b'\in B$.
Likewise, the equality
\[
f\circ\widetilde{\ast}_i
=
\mu_i\circ(f\otimes f)
\]
between linear maps $B\otimes B\to X$ may be checked on pure tensors.
\end{proof}

\begin{remark}\label{RemarkNonPureTensors}
The domain $D_i$ of a bilinear partial operation may contain tensors which are not linear combinations of pure tensors belonging to~$D_i$. Such tensors do not affect the defined products, but they do affect the homomorphisms of algebras with bilinear partial operations. For example, let $n=1$, let $X=\F^2$, with canonical basis $\{e_1,e_2\}$, and let $D$ be the subspace of $X\otimes X$ spanned by $e_1\otimes e_2+e_2\otimes e_1$. Since this tensor has rank~$2$, $D$ contains no non-zero pure tensors. Hence the partial algebras $(X,D,0)$ and $(X,0,0)$ have the same defined products, namely those involving~$0$, but they are not isomorphic in~$\PAlg$, since an isomorphism $f$ would satisfy $(f\otimes f)(D)=0$, whereas $f\otimes f$ is injective. On the other hand, by \Cref{LemMorphismFromTotalN}, the homomorphisms whose domain is a total algebra only depend on the defined products of the codomain.
\end{remark}

\begin{lemma}\label{LemPartialFromBilinearN}
Let $Y$ be a vector subspace of a vector space $Z$, and let
\[
m_i\colon Y\times Y\to Z,
\qquad i=1,\ldots,n,
\]
be bilinear maps. Denote by $\widetilde m_i\colon Y\otimes Y\to Z$ their linearisations. Then, for every $i=1,\ldots,n$,
\[
D_i\coloneqq\widetilde m_i^{-1}(Y)
\]
is a vector subspace of $Y\otimes Y$, and
\[
\mu_i\coloneqq\left.\widetilde m_i\right|_{D_i}
\colon D_i\to Y
\]
is a linear map, so that $\bigl(Y,(D_i,\mu_i)_{i=1}^n\bigr)$ is an algebra with $n$ bilinear partial operations. For $u,v\in Y$, the product $u\ast_i v$ is defined precisely when $m_i(u,v)\in Y$, and in this case
\[
u\ast_i v=m_i(u,v).
\]
\end{lemma}

\begin{proof}
Since $Y\leq Z$ and $\widetilde m_i$ is linear, the inverse image
$D_i=\widetilde m_i^{-1}(Y)$ is a vector subspace of $Y\otimes Y$.
The remaining assertions follow immediately from the definitions.
\end{proof}

\section{The external weak actor of an algebra with \texorpdfstring{$n$}{n} multiplications}\label{SecExternalWeakActorN}

Let $\V$ be an action accessible, operadic variety of algebras with $n$ bilinear multiplications $\ast_1,\ldots,\ast_n$, determined by a set of multilinear identities
\[
\Phi_{t,r}(x_1,\ldots,x_t)=0,
\]
where $t=\deg \Phi_{t,r}\geq2$, and let $X$ be an object of $\V$. In particular, every algebra with $n$ trivial multiplications is an object of $\V$.

Consider a $2n$-tuple
\[
f=
(f\ast_1-,-\ast_1 f,\ldots,f\ast_n-,-\ast_n f)
\in\End(X)^{2n}.
\]
Whenever a multilinear polynomial $\Phi_{t,r}$ is evaluated on elements of $X$ together with one occurrence of $f$, products involving $f$ are interpreted by
\[
f\ast_i x=(f\ast_i-)(x),
\qquad
x\ast_i f=(-\ast_i f)(x),
\]
for any $i=1,\ldots,n$. Since $f$ occurs exactly once, this determines an element of $X$.

\begin{definition}\label{ExternalWeakActorUnderlyingSpace}
For an object $X$ of $\V$, let $\E(X)$ be the subspace of $\End(X)^{2n}$ consisting of all $2n$-tuples
\[
f=(f\ast_1-,-\ast_1 f,\ldots,f\ast_n-,-\ast_n f)
\]
such that
\[
\Phi_{t,r}(\alpha_1,\ldots,\alpha_t)=0
\]
for every defining identity $\Phi_{t,r}=0$, every $j=1,\ldots,t$, and every choice of $x_k\in X$, with $k\neq j$, where
\[
\alpha_j=f
\qquad\text{and}\qquad
\alpha_k=x_k\qquad\text{for }k\neq j.
\]
\end{definition}

The multilinearity of the defining identities implies that $\E(X)$ is an $\F$-vector subspace of $\End(X)^{2n}$.

\begin{remark}\label{RemarkExternalActorAbelian}
Suppose that the $n$ multiplications on the algebra $X$ are trivial. Then the identities of degree at least~$3$ impose no condition on the elements of $\E(X)$. Indeed, every evaluation of such an identity involving one element of $\End(X)^{2n}$ and all the remaining arguments in $X$ contains a product of two elements of $X$ and hence vanishes. Hence, $\E(X)$ is the subspace of $\End(X)^{2n}$ determined by the identities of degree~$2$ satisfied by $\V$. In particular, if $\V$ satisfies no non-trivial identity of degree~$2$, then
\[
\E(X)=\End(X)^{2n}.
\]
\end{remark}

Now, let
\[
(x_1\circ x_2)\diamond x_3
=
W_{\circ,\diamond}(x_1,x_2;x_3),
\qquad
\circ,\diamond\in\Omega_2'
\]
be a choice of identities for Axiom~2, and let $f,g\in\End(X)^{2n}$. For every $\circ,\diamond\in\Omega_2'$ and $x\in X$, the expression
\[
W_{\circ,\diamond}(f,g;x)
\]
is obtained by substituting $f$, $g$, and $x$ for $x_1$, $x_2$, and $x_3$, respectively, in the chosen set of identities for Axiom~2. Every monomial occurring in this expression has one of the forms
\[
f(gx),\quad f(xg),\quad (gx)f,\quad (xg)f,
\]
\[
g(fx),\quad g(xf),\quad (fx)g,\quad (xf)g,
\]
where each juxtaposition denotes an operation in $\Omega_2'$. Hence, each monomial is evaluated using the endomorphisms determined by $f$ and $g$, and therefore belongs to $X$. Thus,
\[
W_{\circ,\diamond}(f,g;-)\colon
X\to X \colon
x\mapsto W_{\circ,\diamond}(f,g;x),
\]
is a linear endomorphism of $X$.

For each $i=1,\ldots,n$, we define the bilinear map
\[
m_X^{(i)}\colon
\End(X)^{2n}\times\End(X)^{2n}\to\End(X)^{2n}
\]
by
\begin{align*}
m_X^{(i)}(f,g)
=
\bigl(
&W_{\ast_i,\ast_1}(f,g;-),
W_{\ast_i,\ast_1^{\op}}(f,g;-),\\
&\ldots,
W_{\ast_i,\ast_n}(f,g;-),
W_{\ast_i,\ast_n^{\op}}(f,g;-)
\bigr).
\end{align*}
The $2n$ components describe the endomorphisms
\[
m_X^{(i)}(f,g)\ast_j-
\quad \text{ and } \quad
-\ast_j m_X^{(i)}(f,g)
\]
for any $j=1,\ldots,n$.

In general, for $f,g\in\E(X)$, the $2n$-tuple $m_X^{(i)}(f,g)$ need not belong to $\E(X)$. For every $i=1,\ldots,n$, let
\[
\widetilde m_X^{(i)}\colon\E(X)\otimes\E(X)\to\End(X)^{2n}
\]
be the linearisation of the restriction of $m_X^{(i)}$ to $\E(X)\times\E(X)$, and define
\[
D_X^{(i)}\coloneqq\bigl(\widetilde m_X^{(i)}\bigr)^{-1}\bigl(\E(X)\bigr)\leq\E(X)\otimes\E(X)
\]
and
\[
\mu_X^{(i)}\coloneqq\left.\widetilde m_X^{(i)}\right|_{D_X^{(i)}}\colon D_X^{(i)}\to\E(X).
\]
By \Cref{LemPartialFromBilinearN},
\[
\bigl(\E(X),(D_X^{(i)},\mu_X^{(i)})_{i=1}^n\bigr)
\]
is an algebra with $n$ bilinear partial operations. By abuse of notation, we shall continue to denote it simply by $\E(X)$.

The $i$-th partial operation of $\E(X)$ will be denoted by $\langle-,-\rangle_i$. More precisely, for $f,g\in\E(X)$, the product $\langle f,g\rangle_i$ is defined if and only if
\[
f\otimes g\in D_X^{(i)},
\]
or, equivalently, if and only if
\[
m_X^{(i)}(f,g)\in\E(X).
\]
Whenever it is defined, we have
\[
\langle f,g\rangle_i=\mu_X^{(i)}(f\otimes g)=m_X^{(i)}(f,g).
\]

Notice that each $D_X^{(i)}$ contains the kernel of $\widetilde m_X^{(i)}$, which in general contains tensors that are not pure, so that the phenomenon described in \Cref{RemarkNonPureTensors} may occur. This does not affect the homomorphisms $\widetilde U(B)\to\E(X)$ used below to represent actions, since, by \Cref{LemMorphismFromTotalN}, they only depend on the defined products of~$\E(X)$.

\begin{remark}
Let
\[
\begin{tikzcd}
X \arrow[r, "k"]
&
A \arrow[r, shift left, "\alpha"]
&
B \arrow[l, shift left, "\beta"]
\end{tikzcd}
\]
be a split extension in $\V$, and consider the corresponding derived action of $B$ on $X$. For each $b\in B$, we define
\[
\varphi(b)
=
(b\ast_1-,-\ast_1 b,\ldots,b \ast_n -, - \ast_n b).
\]
Since $B\ltimes X$ is an object of $\V$, evaluating each defining identity of $\V$ in $B\ltimes X$, with $b$ in one argument and elements of $X$ in the remaining ones, shows that $\varphi(b)$ satisfies the conditions of \Cref{ExternalWeakActorUnderlyingSpace}, that is, $\varphi(b)\in\E(X)$. We therefore obtain a linear map
\[
\varphi\colon B\to\E(X).
\]
Moreover, evaluating the chosen identities for Axiom~2 in the semidirect product $B\ltimes X$, with $x_1=b$ and $x_2=b'$ in $B$ and $x_3=x$ in $X$, we obtain
\[
m_X^{(i)}\bigl(\varphi(b),\varphi(b')\bigr)=\varphi(b\ast_i b')\in\E(X)
\]
for every $i=1,\ldots,n$ and $b,b'\in B$. Hence, $\varphi(b)\otimes\varphi(b')\in D_X^{(i)}$ and
\[
\langle\varphi(b),\varphi(b')\rangle_i=\varphi(b\ast_i b').
\]
By \Cref{LemMorphismFromTotalN}, $\varphi\colon\widetilde U(B)\to\E(X)$ is a morphism in $\PAlg^n$.
\end{remark}

We can now state the main result of this section.

\begin{theorem}\label{ExternalWeakRepresentationN}
Let $\V$ be an action accessible, operadic variety of algebras with~$n$ bilinear multiplications, and fix a choice of identities for Axiom~2. For every object~$X$ of $\V$, there exists a natural monomorphism of functors
\[
\tau\colon
\SplExt(-,X)
\rightarrowtail
\Hom_{\PAlg^n}(\widetilde U(-),\E(X)),
\]
where $\widetilde U\colon\V\to\PAlg^n$ is the functor of \Cref{RemUTildeN}. For every object $B$ of $\V$, the component~$\tau_B$ is the injection which sends an element of $\SplExt(B,X)$ to the homomorphism
\[
\widetilde U(B)\to\E(X) \colon
b \mapsto
(b\ast_1-,-\ast_1 b,\ldots,b \ast_n -, - \ast_n b).
\]
\end{theorem}

\begin{proof}
It follows from the discussion above that, for any object $B$ of $\V$, the map
\[
\tau_B \colon \SplExt(B,X) \to \Hom_{\PAlg^n}(\widetilde U(B),\E(X))
\]
is well defined. Moreover, it is an injection, since every element of $\SplExt(B,X)$ is uniquely determined by the corresponding derived action of $B$ on $X$.

Finally, the collection $\{\tau_B\}_B$ gives rise to a natural transformation since, for every morphism $h\colon B'\to B$ in $\V$, the diagram in $\Set$
\[
\begin{tikzcd}
\SplExt(B,X)
    \arrow[r, "\tau_B"]
    \arrow[d, "{h^*}"']
&
\Hom_{\PAlg^n}(\widetilde U(B),\E(X))
    \arrow[d, "{- \circ h}"]
\\
\SplExt(B',X)
    \arrow[r, "\tau_{B'}"']
&
\Hom_{\PAlg^n}(\widetilde U(B'),\E(X))
\end{tikzcd}
\]
commutes. Thus, $\tau$ is a natural monomorphism of functors.
\end{proof}

\begin{proposition}\label{ImageExternalRepresentationTwo}
Let $B$ and $X$ be objects of $\V$, and let
\[
\varphi\colon\widetilde U(B)\to\E(X)
\]
be a morphism in $\PAlg^n$. Then $\varphi \in \Imm(\tau_B)$ if and only if
\[
\Phi_{t,r}(\alpha_1,\ldots,\alpha_t)=0
\]
for every defining identity $\Phi_{t,r}=0$ of $\V$ and every choice of $\alpha_1,\ldots,\alpha_t\in B\cup X$ such that at least one $\alpha_j$ belongs to $B$ and at least one $\alpha_k$ belongs to $X$. Here, the mixed products between elements of $B$ and $X$ are determined using the notation of \Cref{rem_b*}.
\end{proposition}

\begin{proof}
If $\varphi \in \Imm(\tau_B)$, then it is induced by a derived action of $B$ on $X$. Hence, the corresponding semidirect product $B\ltimes X$ is an object of $\V$ and satisfies all the identities above.

Conversely, suppose that all the mixed evaluations of the defining identities vanish. Since the identities $\Phi_{t,r}=0$ are multilinear, their evaluation on arbitrary elements of $B\oplus X$ is a sum of evaluations obtained by choosing each variable in either $B \times \{0\}$ or $\{0\} \times X$. The evaluations involving only elements of $B$ or only elements of $X$ vanish because $B$ and $X$ are objects of $\V$, while the remaining evaluations vanish by hypothesis. Therefore, $B\ltimes X$ satisfies all the defining identities of $\V$ and is an object of $\V$. Thus, $\varphi$ defines a derived action of $B$ on $X$, i.e., $\varphi \in \Imm(\tau_B)$.
\end{proof}

\begin{definition}
The algebra $\E(X)$ is called the \emph{external weak actor} of $X$, and the pair $(\E(X),\tau)$ is called an \emph{external weak representation} of the functor $\SplExt(-,X)$. If~$\tau$ is a natural isomorphism, then $\E(X)$ is called the \emph{external actor} of $X$. Finally, a morphism $\varphi\colon\widetilde U(B)\to\E(X)$ is called an \emph{acting morphism} if $\varphi\in\Imm(\tau_B)$.
\end{definition}

As in the case of varieties of algebras with one multiplication, a connection between the construction of the external weak actor and weak action representability is given by the following.

\begin{proposition}\label{ExternalActorTwoWRA}
If $\E(X)$ is total and the resulting algebra with $n$ bilinear multiplications is an object of $\V$, then $(\E(X),\tau)$ is a weak representation of $\SplExt(-,X)$.
\end{proposition}

\begin{proof}
If $D_X^{(i)}=\E(X)\otimes\E(X)$ for every $i=1,\ldots,n$, then
\[
\bigl(\E(X),(D_X^{(i)},\mu_X^{(i)})_{i=1}^n\bigr)=\widetilde U(\E(X)),
\]
where $\E(X)$ is seen as an algebra with the $n$ bilinear multiplications $\langle-,-\rangle_i$, for $i=1,\ldots,n$. If this algebra is an object of $\V$, then, since $\widetilde U$ is fully faithful, we have
\[
\Hom_{\PAlg^n}(\widetilde U(-),\E(X))\cong\Hom_\V(-,\E(X)),
\]
and we obtain a natural monomorphism of functors
\[
\tau\colon
\SplExt(-,X)
\rightarrowtail
\Hom_\V(-,\E(X)).
\]
Thus, $(\E(X),\tau)$ is a weak representation of $\SplExt(-,X)$.
\end{proof}

Every algebra $X$ of $\V$ acts on itself via left and right multiplications, giving rise to an inner homomorphism into $\E(X)$, which is described as follows.

\begin{proposition}\label{InnerMultipliersN}
For every object $X$ of $\V$, the linear map
\[
\operatorname{Inn}\colon \widetilde U(X)\to\E(X)\colon
x\mapsto
\bigl(
x\ast_1-,-\ast_1x,\ldots,
x\ast_n-,-\ast_nx
\bigr)
\]
is a morphism in $\PAlg^n$. Its image
\[
\Inn(X)=
\left\{
\bigl(
x\ast_1-,-\ast_1x,\ldots,
x\ast_n-,-\ast_nx
\bigr)
\mid 
x\in X
\right\}
\]
is an object of $\V$.
\end{proposition}

\begin{proof}
Since the defining identities of $\V$ hold in $X$, the $2n$-tuple
\[
(x\ast_1-,-\ast_1 x,\ldots,x\ast_n-,-\ast_n x)
\]
belongs to $\E(X)$ for every $x\in X$. Moreover, the chosen identities for Axiom~2 give
\[
m_X^{(i)}\bigl(
\operatorname{Inn}(x),
\operatorname{Inn}(y)
\bigr)
=
\operatorname{Inn}(x\ast_i y)\in\E(X)
\]
for every $i=1,\ldots,n$ and $x,y\in X$. Hence, $\operatorname{Inn}(x)\otimes\operatorname{Inn}(y)\in D_X^{(i)}$ and
\[
\langle\operatorname{Inn}(x),\operatorname{Inn}(y)\rangle_i=\operatorname{Inn}(x\ast_iy),
\]
so that $\operatorname{Inn}$ is a morphism in $\PAlg^n$ by \Cref{LemMorphismFromTotalN}. In particular, for all $f,g\in\Inn(X)$ and $i=1,\ldots,n$, the product $\langle f,g\rangle_i$ is defined and belongs to $\Inn(X)$, so that the~$n$ partial operations of $\E(X)$ restrict to bilinear multiplications on $\Inn(X)$. Finally, since the corestriction $\operatorname{Inn} \colon X \to \Inn(X)$ is a surjective homomorphism, $\Inn(X)$ satisfies all the defining identities of $\V$, and is therefore an object of $\V$.
\end{proof}

\begin{remark}
The kernel of $\operatorname{Inn}$ is the annihilator
\[
\Ann(X)
=
\left\{
x\in X
\mid
x\ast_i y=y\ast_i x=0, \; \forall i=1,\ldots,n, \, \forall y\in X
\right\}.
\]
Thus, one has
\[
\Inn(X)\cong X/\Ann(X).
\]
\end{remark}

\begin{remark}\label{RemarkCommutativeOperationTwo}
Suppose that one multiplication $\ast_i$ is commutative or anti-commutative, that is
\[
x\ast_i y=\varepsilon_i(y\ast_i x),
\]
with $\varepsilon_i\in\{1,-1\}$. Then
\[
x\ast_i f=\varepsilon_i(f\ast_i x)
\]
for every $f\in\E(X)$ and $x\in X$. In particular, if every multiplication is commutative or anti-commutative, then $\E(X)$ is isomorphic to a subspace of $\End(X)^n$. Explicitly, the isomorphism is given by
\[
(f_1,\ldots,f_n)
\longmapsto
(f_1,\varepsilon_1f_1,\ldots,f_n,\varepsilon_nf_n).
\]
In this situation, unless otherwise specified, we choose the identities for Axiom~2 corresponding to the pairs $(\circ,\ast_i^{\op})$, with $\circ\in\Omega_2'$, to be
\[
(x_1\circ x_2)\ast_i^{\op}x_3=\varepsilon_i\,W_{\circ,\ast_i}(x_1,x_2;x_3),
\]
which is possible since $(x_1\circ x_2)\ast_i^{\op}x_3=x_3\ast_i(x_1\circ x_2)=\varepsilon_i\,(x_1\circ x_2)\ast_ix_3$. With this choice, for all $f,g\in\End(X)^{2n}$ and $j=1,\ldots,n$, the right $\ast_i$-component of $m_X^{(j)}(f,g)$ is $\varepsilon_i$ times its left $\ast_i$-component.
\end{remark}

\begin{remark}\label{RemarkOneOperationCase}
Suppose that all the underlying varieties of algebras with one multiplication $\V_{\ast_1}, \ldots, \V_{\ast_n}$ are action accessible and operadic. If $f\in\E(X)$, then the pair
\[
f_i=(f\ast_i-,-\ast_i f)\in\E_{\V_{\ast_i}}(X,\ast_i)
\]
for any $i=1,\ldots,n$. It follows that $\E(X)$ consists of all $n$-tuples
\[
(f_1,\ldots,f_n)
\in
\E_{\V_{\ast_1}}(X,\ast_1)\times \cdots \times \E_{\V_{\ast_n}}(X,\ast_n)
\]
which satisfy all the mixed identities of $\V$. Thus, once the external weak actors associated with the individual multiplications are known, the construction of $\E(X)$ reduces to determining the additional compatibility conditions arising from the identities involving more than one multiplication.

If, for every $i=1,\ldots,n$, the identities for Axiom~2 corresponding to the pairs $(\ast_i,\ast_i)$ and $(\ast_i,\ast_i^{\op})$ are chosen to be $\lambda/\mu$-rules of $\V_{\ast_i}$, then the two components of $m_X^{(i)}(f,g)$ corresponding to the multiplication $\ast_i$ are determined by these $\lambda/\mu$-rules.
\end{remark}

\subsection{An obstruction to action representability}\label{SecObstruction}

Since $\V$ is an Orzech category of interest, every object $X$ of $\V$ admits a universal strict general actor $\USGA(X)$ in $\Alg^n$, together with a natural monomorphism of functors
\[
\mu\colon
\SplExt(-,X)
\rightarrowtail
\Hom_{\Alg^n}(\overline U(-),\USGA(X)),
\]
where $\overline U\colon\V\to\Alg^n$ is the forgetful functor; see~\cite{Casas}. The constructions of $\USGA(X)$ and $\E(X)$ use the same choice of identities for Axiom~2. Thus, for every object $X$ of $\V$, the actions on $X$ are encoded both by $\mu$ and by the natural monomorphism of functors
\[
\tau\colon
\SplExt(-,X)
\rightarrowtail
\Hom_{\PAlg^n}(\widetilde U(-),\E(X))
\]
of \Cref{ExternalWeakRepresentationN}. It is often more convenient to work with the external weak actor~$\E(X)$, since its construction is simpler than that of the universal strict general actor $\USGA(X)$; see~\cite[Section~4]{Casas}.
In general, $\USGA(X)$ and $\E(X)$ are different, since $\USGA(X)$ is generated by the elements which arise from derived actions on~$X$. For instance, for $n=1$, $\V=\Nil_2(\Alg)$ and $X=\F$ with trivial multiplication, every derived action on $X$ is trivial, so that $\USGA(X)=0$, whereas $\E(X)=\End(X)^2\cong\F^2$.

Moreover, by~\cite[Theorem~3.9]{Casas}, an object $X$ of $\V$ admits an actor if and only if the semidirect product $\USGA(X)\ltimes X$ is an object of $\V$, and in this case $[X]\cong\USGA(X)$. We now give an analogous characterisation in terms of the external weak actor, and we deduce from it an obstruction to action representability.

Let $X$ be an object of $\V$. We denote by
\[
S_X\subseteq\E(X)
\]
the set of all $2n$-tuples $\varphi(b)$, where $B$ is an object of $\V$, $b\in B$, and $\varphi\colon\widetilde U(B)\to\E(X)$ is an acting morphism. In other words, $S_X$ is the union of the images of all acting morphisms, that is, the set of all $2n$-tuples
\[
(b\ast_1-,-\ast_1b,\ldots,b\ast_n-,-\ast_nb)
\]
arising from the derived actions on $X$ of the objects of $\V$.

The main result of this subsection characterises the objects of $\V$ admitting an actor in terms of $S_X$. One of its conditions is that $S_X$ be a vector subspace of $\E(X)$ such that, for all $f,g\in S_X$ and $i=1,\ldots,n$, the product $\langle f,g\rangle_i$ is defined and belongs to $S_X$. In this case, $S_X$, endowed with the restrictions of the operations $\langle-,-\rangle_1,\ldots,\langle-,-\rangle_n$, is an algebra with $n$ bilinear multiplications, and we denote by $S_X\ltimes X$ the algebra with $n$ bilinear multiplications obtained from~\eqref{SemidirectProductNOperations} by means of the \emph{canonical action} of $S_X$ on $X$, defined by
\[
f\ast_ix=(f\ast_i-)(x),
\qquad
x\ast_if=(-\ast_if)(x),
\]
for every $f\in S_X$, $x\in X$ and $i=1,\ldots,n$. The other condition is that $S_X\ltimes X$ be an object of $\V$. The following result is an analogue of~\cite[Theorem~3.9]{Casas} for the external weak actor.

\begin{theorem}\label{ThmActorSX}
Let $X$ be an object of $\V$. Then $X$ admits an actor if and only if the following conditions hold:
\begin{enumerate}
\item $S_X$ is a vector subspace of $\E(X)$ and, for all $f,g\in S_X$ and $i=1,\ldots,n$, the product $\langle f,g\rangle_i$ is defined and belongs to $S_X$;
\item the semidirect product $S_X\ltimes X$ with respect to the canonical action is an object of $\V$.
\end{enumerate}
In this case, $S_X$ is an object of $\V$, and $[X]\cong S_X$.
\end{theorem}

\begin{proof}
Suppose that $X$ admits an actor $[X]$, and let
\[
\theta\colon\SplExt(-,X) \cong \Hom_{\V}(-,[X])
\]
be a natural isomorphism. Let $\sigma_0\in\SplExt([X],X)$ be such that $\theta_{[X]}(\sigma_0)=\id_{[X]}$, and set
\[
\varphi_0=\tau_{[X]}(\sigma_0)\colon\widetilde U([X])\to\E(X).
\]
Let $B$ be an object of $\V$, let $\sigma\in\SplExt(B,X)$, and set $\chi=\theta_B(\sigma)$. By the naturality of $\theta$, we have $\sigma=\chi^*(\sigma_0)$, and the naturality of $\tau$ gives
\[
\tau_B(\sigma)=\tau_B\bigl(\chi^*(\sigma_0)\bigr)=\varphi_0\circ\chi.
\]
Thus, every acting morphism factors through $\varphi_0$, so that $S_X=\varphi_0([X])$. In particular, $S_X$ is a vector subspace of $\E(X)$.

We now show that $\varphi_0$ is injective. Let $K$ be its kernel. By \Cref{LemMorphismFromTotalN},
\[
\varphi_0(k\ast_ik')=\bigl\langle\varphi_0(k),\varphi_0(k')\bigr\rangle_i=\langle0,0\rangle_i=0
\]
for all $k,k'\in K$ and $i=1,\ldots,n$. Hence, $K$ is a subalgebra of $[X]$, and therefore an object of $\V$. Let $\iota\colon K\hookrightarrow[X]$ be the inclusion and $0\colon K\to[X]$ the zero morphism. Then
\[
\tau_K\bigl(\iota^*(\sigma_0)\bigr)=\varphi_0\circ\iota=0=\varphi_0\circ0=\tau_K\bigl(0^*(\sigma_0)\bigr).
\]
Since $\tau_K$ is injective, we get $\iota^*(\sigma_0)=0^*(\sigma_0)$, and the naturality of $\theta$ gives
\[
\iota=\theta_K\bigl(\iota^*(\sigma_0)\bigr)=\theta_K\bigl(0^*(\sigma_0)\bigr)=0.
\]
Hence, $K=0$.

It follows that $\varphi_0$ is a linear bijection from $[X]$ onto $S_X$. By \Cref{LemMorphismFromTotalN}, for all $a,a'\in[X]$ and $i=1,\ldots,n$, the product $\langle\varphi_0(a),\varphi_0(a')\rangle_i$ is defined and
\[
\bigl\langle\varphi_0(a),\varphi_0(a')\bigr\rangle_i=\varphi_0(a\ast_ia')\in S_X.
\]
This proves~(1), and shows that $\varphi_0\colon[X]\to S_X$ is an isomorphism of algebras with~$n$ bilinear multiplications. In particular, $S_X$ is an object of $\V$, and $S_X\cong[X]$. Moreover, by the definition of $\tau$, the derived action of $[X]$ on $X$ associated with $\sigma_0$ is given by
\[
a\ast_ix=\bigl(\varphi_0(a)\ast_i-\bigr)(x),
\qquad
x\ast_ia=\bigl(-\ast_i\varphi_0(a)\bigr)(x),
\]
for all $a\in[X]$, $x\in X$ and $i=1,\ldots,n$. Hence, the map
\[
[X]\ltimes X\to S_X\ltimes X\colon(a,x)\mapsto(\varphi_0(a),x)
\]
is an isomorphism of algebras with $n$ bilinear multiplications. Since $[X]\ltimes X$ is an object of $\V$ by \Cref{rem_b*}, so is $S_X\ltimes X$, which proves~(2).

Conversely, suppose that~(1) and~(2) hold. Then $S_X$ is an object of $\V$, being isomorphic to the subalgebra $S_X \times \{0\}$ of $S_X\ltimes X$, and
\[
\begin{tikzcd}
X \arrow[r, "\iota_2"]
&
S_X\ltimes X \arrow[r, shift left, "\pi_1"]
&
S_X \arrow[l, shift left, "\iota_1"]
\end{tikzcd}
\]
is a split extension in $\V$, whose class in $\SplExt(S_X,X)$ we denote by $\sigma_S$. Let $j\colon S_X\hookrightarrow\E(X)$ be the inclusion. For every object $B$ of $\V$ and every morphism $\psi\colon B\to S_X$ in $\V$, the pullback $\psi^*(\sigma_S)$ is the split extension of $B$ by $X$ with derived action
\[
b\ast_ix=\bigl(\psi(b)\ast_i-\bigr)(x),
\qquad
x\ast_ib=\bigl(-\ast_i\psi(b)\bigr)(x),
\]
so that $\tau_B\bigl(\psi^*(\sigma_S)\bigr)=j\circ\psi$. This defines a natural transformation
\[
\Hom_{\V}(-,S_X)\to\SplExt(-,X)\colon\psi\mapsto\psi^*(\sigma_S).
\]
Its components are injective: if $\psi^*(\sigma_S)=\psi'^*(\sigma_S)$, then $j\circ\psi=j\circ\psi'$, and hence $\psi=\psi'$. They are also surjective. Indeed, let $\sigma\in\SplExt(B,X)$ and $\varphi=\tau_B(\sigma)$. By the definition of $S_X$, the map $\varphi$ takes values in $S_X$, and, by \Cref{LemMorphismFromTotalN},
\[
\varphi(b\ast_ib')=\bigl\langle\varphi(b),\varphi(b')\bigr\rangle_i
\]
for all $b,b'\in B$ and $i=1,\ldots,n$. Hence, the corestriction $\psi\colon B\to S_X$ of $\varphi$ is a morphism in $\V$, and
\[
\tau_B\bigl(\psi^*(\sigma_S)\bigr)=j\circ\psi=\varphi=\tau_B(\sigma).
\]
Since $\tau_B$ is injective, $\psi^*(\sigma_S)=\sigma$. Therefore, $\SplExt(-,X)\cong\Hom_\V(-,S_X)$, and~$S_X$ is the actor of $X$.
\end{proof}

When $X$ has trivial multiplications, some elements of $S_X$ can be described explicitly.

\begin{lemma}\label{LemSquareZero}
Let $X$ be an object of $\V$ with trivial multiplications, and let
\[
f=(f\ast_1-,-\ast_1f,\ldots,f\ast_n-,-\ast_nf)\in\E(X)
\]
be such that $u\circ v=0$ for all, not necessarily distinct, components $u$ and $v$ of $f$. Let $B=\F b$ be the one-dimensional algebra with trivial multiplications. Then the bilinear maps determined by
\[
b\ast_ix=(f\ast_i-)(x),
\qquad
x\ast_ib=(-\ast_if)(x),
\qquad
i=1,\ldots,n,
\]
define a derived action of $B$ on $X$, whose acting morphism sends $b$ to $f$. In particular, $f\in S_X$.
\end{lemma}

\begin{proof}
Since the defining identities of $\V$ have degree at least~$2$, the algebra $B$ is an object of $\V$. As in the proof of \Cref{ImageExternalRepresentationTwo}, by multilinearity it suffices to show that every defining identity $\Phi_{t,r}=0$ of $\V$ vanishes in $B\oplus X$, endowed with the multiplications~\eqref{SemidirectProductNOperations}, at every choice of $\alpha_1,\ldots,\alpha_t\in\{b\}\cup X$, where $b$ and $x\in X$ are identified with $(b,0)$ and $(0,x)$, respectively. If none of the $\alpha_j$ is equal to $b$, this holds because $X$ is an object of $\V$; if all of them are equal to $b$, it holds because $B$ is an object of $\V$; and if exactly one of them is equal to $b$, it holds because $f\in\E(X)$.

Suppose now that at least two of the $\alpha_j$ are equal to $b$ and at least one of them belongs to $X$. In $B\oplus X$, all products of two elements of $B$ and all products of two elements of $X$ vanish, while a product of $b$ and an element $y\in X$, in either order, is of the form $u(y)$ for some component $u$ of $f$. By induction on the degree, a monomial evaluated at elements of $\{b\}\cup X$ vanishes if it contains at least two elements of $X$, while, if it contains exactly one element $x$ of $X$ and $s$ occurrences of $b$, its value is either zero or of the form $u_1\circ\cdots\circ u_s(x)$, where $u_1,\ldots,u_s$ are components of $f$. Since $s\geq2$ and $u_1\circ u_2=0$, every monomial of $\Phi_{t,r}$ vanishes at $(\alpha_1,\ldots,\alpha_t)$. Therefore, $B\ltimes X$ is an object of $\V$, and the corresponding acting morphism sends $b$ to~$f$.
\end{proof}

In order to state the obstruction, we fix some further notation. For an object~$X$ of $\V$, the expressions $W_{\circ,\diamond}(f,g;x)$, and hence the bilinear maps $m_X^{(1)},\ldots,m_X^{(n)}$, are defined for all $f,g\in\End(X)^{2n}$. We denote by $\widehat{\E}(X)$ the algebra with $n$ bilinear multiplications whose underlying vector space is $\End(X)^{2n}$, with the multiplications $m_X^{(1)},\ldots,m_X^{(n)}$. If, moreover, $X$ has trivial multiplications, we denote by $\widehat{\E}(X)\ltimes X$ the algebra with $n$ bilinear multiplications whose underlying vector space is $\End(X)^{2n}\oplus X$, with multiplications
\[
(f,x)\ast_i(g,y)=\bigl(m_X^{(i)}(f,g),\,(f\ast_i-)(y)+(-\ast_ig)(x)\bigr),
\qquad i=1,\ldots,n.
\]
We regard $\E(X)$ and $X$ as subsets of $\widehat{\E}(X)\ltimes X$ via $f\mapsto(f,0)$ and $x\mapsto(0,x)$. If a polynomial $\Phi(x_1,\ldots,x_k)$ is defined, in the sense of \Cref{DefPartialEvaluationN}, at some $(f_1,\ldots,f_k)\in\E(X)^k$, then its value coincides with $\Phi(f_1,\ldots,f_k)$ computed in $\widehat{\E}(X)$.

Finally, for every $d\geq1$, we consider the commutative associative algebra
\[
\Lambda_d=\F[\xi_1,\ldots,\xi_d]/(\xi_1^2,\ldots,\xi_d^2),
\]
which has a basis consisting of the monomials $\xi_T=\prod_{s\in T}\xi_s$, for $T\subseteq\{1,\ldots,d\}$. If~$X$ has trivial multiplications, we regard the vector space $X\otimes\Lambda_d$ as an algebra with $n$ trivial multiplications, which is an object of $\V$.

\begin{theorem}\label{ThmTensorObstruction}
Let $X$ be an object of $\V$ with trivial multiplications, let $d\geq1$, and suppose that $X\otimes\Lambda_d$ admits an actor. Let $\Phi(x_1,\ldots,x_k)=0$ be a multilinear identity of $\V$, and let $J\subseteq\{1,\ldots,k\}$ with $|J|\leq d$. Then
\[
\Phi(\alpha_1,\ldots,\alpha_k)=0
\]
in $\widehat{\E}(X)\ltimes X$, for all $\alpha_j\in\E(X)$ with $j\in J$ and all $\alpha_j\in X$ with $j\notin J$.
\end{theorem}

\begin{proof}
Let $X'=X\otimes\Lambda_d$ and, for $s=1,\ldots,d$, let $\ell_s\in\End(\Lambda_d)$ be the multiplication by $\xi_s$. The endomorphisms $\ell_1,\ldots,\ell_d$ commute pairwise and satisfy $\ell_s^2=0$. For $T\subseteq\{1,\ldots,d\}$, we set $\ell_T=\prod_{s\in T}\ell_s$, so that $\ell_T(1)=\xi_T\neq0$. For $u=(u_1,\ldots,u_{2n})\in\End(X)^{2n}$ and $e\in\End(\Lambda_d)$, we write
\[
u\otimes e=(u_1\otimes e,\ldots,u_{2n}\otimes e)\in\End(X')^{2n}.
\]

By the definition of $m_X^{(i)}$, every component of $m_X^{(i)}(u,v)$, for $u,v\in\End(X)^{2n}$, is a linear combination of composites $p\circ q$ and $q\circ p$, where $p$ is a component of $u$ and~$q$ is a component of $v$, with coefficients that only depend on the chosen identities for Axiom~2. Since $(p\otimes e)\circ(q\otimes e')=(p\circ q)\otimes(e\circ e')$, we obtain
\begin{equation}\label{EqTensorProducts}
m_{X'}^{(i)}(u\otimes e,v\otimes e')=m_X^{(i)}(u,v)\otimes(e\circ e')
\end{equation}
for all $u,v\in\End(X)^{2n}$, all commuting $e,e'\in\End(\Lambda_d)$ and all $i=1,\ldots,n$. Moreover,
\[
(u\otimes e)\ast_i(x\otimes w)=(u\ast_ix)\otimes e(w)
\qquad\text{and}\qquad
(x\otimes w)\ast_i(u\otimes e)=(x\ast_iu)\otimes e(w)
\]
for all $x\in X$ and $w\in\Lambda_d$, where $u\ast_ix$ and $x\ast_iu$ denote the images of $x$ under the corresponding components of $u$.

Next, we show that $f\otimes e\in\E(X')$ for all $f\in\E(X)$ and $e\in\End(\Lambda_d)$. Since~$X$ and $X'$ have trivial multiplications, by \Cref{RemarkExternalActorAbelian} the subspaces $\E(X)$ and $\E(X')$ are determined by the defining identities of degree~$2$ of $\V$. When one of the two variables of such an identity is replaced by a $2n$-tuple of endomorphisms and the other one by an element of $X$ (respectively,~$X'$), one obtains a linear relation between the components of the $2n$-tuple, which is the same for $X$ and for $X'$. Since these relations are preserved by $u\mapsto u\otimes e$, the claim follows.

Now, let $\rho\colon J\to\{1,\ldots,d\}$ be an injective map, and set
\[
\alpha'_j=\alpha_j\otimes\ell_{\rho(j)}\quad\text{for } j\in J,
\qquad
\alpha'_j=\alpha_j\otimes1\quad\text{for } j\notin J.
\]
For $j\in J$, we have $\alpha'_j\in\E(X')$, and all the composites of two components of $\alpha'_j$ vanish, since $\ell_{\rho(j)}^2=0$. Hence, $\alpha'_j\in S_{X'}$ by \Cref{LemSquareZero}. Since $X'$ admits an actor, \Cref{ThmActorSX} shows that $S_{X'}$ is closed under the operations $\langle-,-\rangle_i$, which coincide with $m_{X'}^{(i)}$ on $S_{X'}$, and that $S_{X'}\ltimes X'$ is an object of $\V$. Since $S_{X'}\ltimes X'$ is a subalgebra of $\widehat{\E}(X')\ltimes X'$ and $\Phi=0$ is an identity of $\V$, we get
\[
\Phi(\alpha'_1,\ldots,\alpha'_k)=0
\]
in $\widehat{\E}(X')\ltimes X'$.

On the other hand, let $w$ be a monomial in the variables $x_j$, $j\in L$, for some $L\subseteq\{1,\ldots,k\}$, each of them occurring exactly once, and denote by $w(\alpha)$ and $w(\alpha')$ its evaluations at $(\alpha_j)_{j\in L}$ in $\widehat{\E}(X)\ltimes X$ and at $(\alpha'_j)_{j\in L}$ in $\widehat{\E}(X')\ltimes X'$, respectively. Using~\eqref{EqTensorProducts}, the formulas above, and the fact that all products of two elements of~$X$, or of $X'$, vanish, an induction on the degree of $w$ shows that
\[
w(\alpha')=
\begin{cases}
w(\alpha)\otimes\ell_{\rho(L)} & \text{if } L\subseteq J,\\[0.2em]
w(\alpha)\otimes\xi_{\rho(L\cap J)} & \text{if } L\not\subseteq J.
\end{cases}
\]
Since $\Phi$ is multilinear, it follows that $\Phi(\alpha')$ is equal to $\Phi(\alpha)\otimes\ell_{\rho(J)}$ if $J=\{1,\ldots,k\}$, and to $\Phi(\alpha)\otimes\xi_{\rho(J)}$ otherwise. Since $\ell_{\rho(J)}(1)=\xi_{\rho(J)}\neq0$, for $u\in\End(X)^{2n}$ we have $u\otimes\ell_{\rho(J)}=0$ only if $u=0$ (evaluate its components at the elements $x\otimes1$), and for $y\in X$ we have $y\otimes\xi_{\rho(J)}=0$ only if $y=0$. Hence, the equality $\Phi(\alpha')=0$ implies $\Phi(\alpha)=0$.
\end{proof}

\begin{corollary}\label{CorTensorObstruction}
Let $X$ be an object of $\V$ with trivial multiplications.
\begin{enumerate}
\item Suppose that $\V$ is action representable, or that $X$ is infinite-dimensional and admits an actor. Then every multilinear identity of $\V$ holds in $\widehat{\E}(X)\ltimes X$ whenever its variables are evaluated at elements of $\E(X)\cup X$. Moreover, $\E(X)$ is total and, endowed with its $n$ bilinear operations, it is an object of~$\V$, and so is the semidirect product $\E(X)\ltimes X$ with respect to the canonical action
\[
f\ast_ix=(f\ast_i-)(x),
\qquad
x\ast_if=(-\ast_if)(x).
\]
\item Let $\Phi(x_1,\ldots,x_k)=0$ be a multilinear identity of $\V$. If $\Phi(\alpha_1,\ldots,\alpha_k)\neq0$ in $\widehat{\E}(X)\ltimes X$ for some $\alpha_1,\ldots,\alpha_k\in\E(X)\cup X$, then $X\otimes\Lambda_k$ does not admit an actor. In particular, $\V$ is not action representable. This is the case, for instance, if $\E(X)$ does not satisfy the identity $\Phi=0$.
\end{enumerate}
\end{corollary}

\begin{proof}{\ }
\begin{enumerate}
\item[(1)] If $\V$ is action representable, then $X\otimes\Lambda_d$ admits an actor for every $d\geq1$. The same holds if $X$ is infinite-dimensional and admits an actor, because in this case $X\otimes\Lambda_d$ and $X$ have the same dimension, and hence they are isomorphic as algebras with trivial multiplications. The first assertion then follows from \Cref{ThmTensorObstruction}. We now show that $\E(X)$ is total. Let $f,g\in\E(X)$ and $i\in\{1,\ldots,n\}$, and let $X'=X\otimes\Lambda_2$, which admits an actor. With the notation of the proof of \Cref{ThmTensorObstruction}, the elements $f\otimes\ell_1$ and $g\otimes\ell_2$ belong to $S_{X'}$, so that, by \Cref{ThmActorSX} and the identities~\eqref{EqTensorProducts},
\[
m_X^{(i)}(f,g)\otimes(\ell_1\circ\ell_2)=m_{X'}^{(i)}(f\otimes\ell_1,g\otimes\ell_2)\in S_{X'}\subseteq\E(X').
\]
As observed in the proof of \Cref{ThmTensorObstruction}, $\E(X)$ and $\E(X')$ are determined by the same linear relations between the components of $2n$-tuples of endomorphisms. Since $(\ell_1\circ\ell_2)(1)=\xi_1\xi_2\neq0$, a $2n$-tuple $u\otimes(\ell_1\circ\ell_2)$ satisfies such a relation only if $u$ does. Hence, $m_X^{(i)}(f,g)\in\E(X)$, so that $\E(X)$ is total, and $\E(X)\oplus X$ is a subalgebra of $\widehat{\E}(X)\ltimes X$, namely the semidirect product $\E(X)\ltimes X$. By multilinearity, the first assertion shows that the defining identities of $\V$ are satisfied by $\E(X)\ltimes X$, which is thus an object of $\V$. Hence, so is its subalgebra $\E(X)\cong\E(X)\times\{0\}$.

\item[(2)] Since $\Phi(\alpha_1,\ldots,\alpha_k)\neq0$, every $\alpha_j$ is non-zero, and hence it belongs to exactly one of $\E(X)$ and $X$. Thus, the first assertion follows from \Cref{ThmTensorObstruction}, applied with $d=k$ and $J=\{j\mid\alpha_j\in\E(X)\}$. If $\E(X)$ does not satisfy $\Phi=0$, then there exist $f_1,\ldots,f_k\in\E(X)$ such that $\Phi$ is defined at $(f_1,\ldots,f_k)$ and $\Phi(f_1,\ldots,f_k)\neq0$. Since this value coincides with the evaluation of $\Phi$ at $(f_1,\ldots,f_k)$ in $\widehat{\E}(X)\ltimes X$, the first assertion applies.
\end{enumerate}
\end{proof}

\begin{remark}\label{RemarkTensorNecessary}
In \Cref{CorTensorObstruction}(2), the algebra $X\otimes\Lambda_k$ cannot, in general, be replaced by $X$. For instance, let $n=1$ and $\V=\Nil_c(\Assoc)$, with $c\geq2$, and choose the identities for Axiom~2 given by associativity, as in \Cref{ExampleExternalAssoc}. Let $X=\F$ be the one-dimensional algebra with trivial multiplication. If an object $B$ of $\V$ acts on~$X$, with $b\ast x=\beta(b)x$ and $x\ast b=\beta'(b)x$ for some linear maps $\beta,\beta'\colon B\to\F$, then the $c$-nilpotency of $B\ltimes X$ gives
\[
\beta(b)^cx=b\ast(b\ast(\cdots(b\ast x)\cdots))=0
\qquad\text{and}\qquad
\beta'(b)^cx=(\cdots((x\ast b)\ast b)\cdots)\ast b=0
\]
for all $b\in B$ and $x\in X$, so that $\beta=\beta'=0$. Thus, every derived action on~$X$ is trivial, and the zero algebra is the actor of~$X$. On the other hand, $\E(X)=\End(X)^2$ is total, with multiplication
\[
\langle f,g\rangle=\bigl(f\ast(g\ast-),(-\ast f)\ast g\bigr),
\]
and it does not satisfy the identity $(\cdots((x_1x_2)x_3)\cdots)x_{c+1}=0$ of $\V$, since $f=(\id_X,\id_X)$ satisfies $\langle f,f\rangle=f\neq0$. By \Cref{CorTensorObstruction}(2), the algebra $X\otimes\Lambda_{c+1}$ does not admit an actor. In particular, $\Nil_c(\Assoc)$ is not action representable.
\end{remark}

\section{Examples}\label{SecEx}

We now specialise the preceding construction to the varieties of algebras introduced in \Cref{SecNOperations} that are Orzech categories of interest.

\subsection{Abelian algebras with $n$ multiplications}\label{SecAb}

Let $X$ be a vector space endowed with $n$ trivial multiplications $x\ast_1 y=\ldots=x\ast_n y=0$.
Every derived action on~$X$ is trivial. Hence, $\E(X)=0$ and there is a natural isomorphism
\[
\SplExt(-,X)\cong\Hom_{\AbAlg^n}(-,0),
\]
so the zero algebra is the actor of $X$. Thus, $\AbAlg^n$ is action representable. This also follows from the fact that $\AbAlg^n$ is an abelian category.

\subsection{Poisson algebras}\label{SecPoisson}

Let $\bchar(\F) \neq 2$. The representability of actions of Poisson algebras was studied in~\cite[Section~5]{CigoliManciniMetere}. We briefly recall the construction given there and interpret it in terms of the external weak actor constructed in the previous section.

Let $X=(X,\cdot,[-,-])$ be a Poisson algebra. Using the same choice of identities for Axiom~2 as in~\cite{CigoliManciniMetere}, that is,
\begin{align*}
(x\cdot y)\cdot z&=x\cdot(y\cdot z), &
z\cdot(x\cdot y)&=(z\cdot x)\cdot y,\\
[[x,y],z]&=[[x,z],y]+[x,[y,z]], &
[z,[x,y]]&=[[z,x],y]-[[z,y],x],\\
[x\cdot y,z]&=[x,z]\cdot y+x\cdot[y,z], &
[z,x\cdot y]&=[z,x]\cdot y+x\cdot[z,y],\\
[x,y]\cdot z&=x\cdot[y,z]-[y,x\cdot z], &
z\cdot[x,y]&=[z\cdot x,y]-[z,y]\cdot x,
\end{align*}
the external weak actor $\E(X)$ is isomorphic, as a vector space, to
\[
\left\{
(f \ast -, - \ast f,[f,-])\in\Bim(X)\times\Der(X)
\ \middle|\
\begin{aligned}
f \ast [x,y]&=[f \ast x,y]-[f,y]\cdot x,\\
[x,y] \ast f&=[x \ast f,y]-x\cdot [f,y],\\
[f,x \cdot y]&=[f,x] \cdot y+x \cdot [f,y]
\end{aligned}
\right\},
\]
where the three identities are required to hold for all $x,y\in X$, $\Bim(X)$ is the associative algebra of bimultipliers of $(X,\cdot)$ and $\Der(X)$ is the Lie algebra of derivations of $(X,[-,-])$.

The two bilinear partial operations on $\E(X)$ are given by
\[
\langle f,g\rangle_\ast=(f \ast (g \ast -), (- \ast f) \ast g, f \ast [g,-] + [f,-] \ast g)
\]
and
\[
\langle f,g\rangle_{[-,-]}=(f \ast [g,-] - [g,f \ast -], [g,-] \ast f - [g,- \ast f], [f,[g,-]]-[g,[f,-]]),
\]
for any $f=(f \ast -, - \ast f,[f,-]), g=(g \ast -, - \ast g,[g,-]) \in \E(X)$. It turns out that $\E(X)$ is a total algebra, which coincides with the algebra with two bilinear operations constructed in~\cite[Theorem~5.6]{CigoliManciniMetere}, where it is denoted by~$[X]$. Thus, the natural monomorphism of \Cref{ExternalWeakRepresentationN} recovers the one of~\cite[Theorem~5.6]{CigoliManciniMetere}.

Now, let $X$ be a commutative Poisson algebra. Since the multiplication $\cdot$ is commutative, \Cref{RemarkCommutativeOperationTwo} gives $x\ast f=f\ast x$ for every $f\in\E(X)$ and every $x\in X$. Hence, the external weak actor in the category $\CPois$ is isomorphic to
\[
\left\{
(f \ast -, [f,-])\in\operatorname{M}(X)\times\Der(X)
\ \middle|\
\begin{aligned}
f \ast [x,y]&=[f \ast x,y]-[f,y]\cdot x,\\
[f,x \cdot y]&=[f,x] \cdot y+x \cdot [f,y],
\end{aligned}
\; \; \forall x,y\in X
\right\},
\]
where $\M(X)$ is the associative algebra of \emph{multipliers} of $(X,\cdot)$, with
\[
\langle f,g\rangle_\ast=(f \ast (g \ast -), f \ast [g,-] + [f,-] \ast g)
\]
and
\[
\langle f,g\rangle_{[-,-]}=(f \ast [g,-] - [g,f \ast -], [f,[g,-]]-[g,[f,-]]).
\]
In this case, the external weak representation is a natural isomorphism, so $\E(X)$ is an external actor of $X$; see~\cite[Section~5]{CigoliManciniMetere}, where it is denoted by~$[X]_c$. Nevertheless, $\E(X)$ need not be a commutative Poisson algebra.

We now show that neither $\Pois$ nor $\CPois$ is action representable.

\begin{theorem}\label{TheoremPoisNotAR}
Let $\bchar(\F)\neq2$. The varieties $\Pois$ and $\CPois$ are not action representable.
\end{theorem}

\begin{proof}
Let $X$ be the two-dimensional vector space with basis $\{e_1,e_2\}$, endowed with trivial operations, and let $u,v\in\End(X)$ be defined by
\[
u(e_1)=0,\qquad u(e_2)=e_1,\qquad v(e_1)=e_2,\qquad v(e_2)=0.
\]
By \Cref{RemarkExternalActorAbelian}, $\E(X)$ is the subspace of $\End(X)^4$ determined by the identities of degree~$2$ of the variety.

Consider first $\V=\Pois$. Its identities of degree~$2$ only express the anti-commutativity of the bracket, and hence impose no restriction on the components $f\ast-$ and $-\ast f$ of the elements $f\in\E(X)$. Hence, we may choose $f,g\in\E(X)$ with all components zero except $f\ast-=u$ and $-\ast g=v$, respectively. In $\widehat{\E}(X)\ltimes X$ we then have
\[
(f\ast e_1)\ast g=v\bigl(u(e_1)\bigr)=0\neq e_1=u\bigl(v(e_1)\bigr)=f\ast(e_1\ast g).
\]
Thus, the associativity of the multiplication of $\Pois$, which is a multilinear identity of degree~$3$, fails at $(f,e_1,g)$. By \Cref{CorTensorObstruction}(2), the sixteen-dimensional algebra $X\otimes\Lambda_3$ with trivial operations does not admit an actor in $\Pois$.

Consider now $\V=\CPois$. In this case, $\E(X)\cong\End(X)\times\End(X)$, its elements being the pairs $f=(f\ast-,[f,-])$, and the description of $\E(X)$ given above yields
\[
\langle f,g\rangle_\ast-\langle g,f\rangle_\ast=\bigl((f\ast-)\circ(g\ast-)-(g\ast-)\circ(f\ast-),\,0\bigr)
\]
for all $f,g\in\E(X)$. Taking $f=(u,0)$ and $g=(v,0)$, we get $u\circ v-v\circ u\neq0$, since $(u\circ v-v\circ u)(e_1)=e_1$. Hence, $\E(X)$ does not satisfy the multilinear identity $x\cdot y-y\cdot x=0$ of $\CPois$, which has degree~$2$, and by \Cref{CorTensorObstruction}(2) the eight-dimensional algebra $X\otimes\Lambda_2$ with trivial operations does not admit an actor in~$\CPois$.
\end{proof}

\subsection{Left-right non-commutative Poisson algebras}\label{SecExternalNPlr}

Let $X$ be a left-right non-commutative Poisson algebra. Thus, $X$ is endowed with an associative multiplication~$\cdot$ and a right Leibniz bracket $[-,-]$ satisfying
\[
[x\cdot y,z]=x\cdot[y,z]+[x,z]\cdot y
\]
and
\[
[x,y\cdot z]=y\cdot[x,z]+[x,y]\cdot z
\]
for all $x,y,z\in X$.

The external weak actor $\E(X)$ consists of all quadruples
\[
f=(f\ast-,-\ast f,[f,-],[-,f])\in\End(X)^4
\]
satisfying the following identities:
\begin{align}
f\ast(x\cdot y)&=(f\ast x)\cdot y, \label{EqNPlrE1}\\
(x\cdot y)\ast f&=x\cdot(y\ast f), \label{EqNPlrE2}\\
x\cdot(f\ast y)&=(x\ast f)\cdot y, \label{EqNPlrE3}\\
[f,[x,y]]&=[[f,x],y]-[[f,y],x], \label{EqNPlrE4}\\
[[x,y],f]&=[[x,f],y]+[x,[y,f]], \label{EqNPlrE5}\\
[x,[f,y]]&=-[x,[y,f]], \label{EqNPlrE6}\\
[f\ast x,y]&=f\ast[x,y]+[f,y]\cdot x, \label{EqNPlrE7}\\
[x\ast f,y]&=x\cdot[f,y]+[x,y]\ast f, \label{EqNPlrE8}\\
[x\cdot y,f]&=x\cdot[y,f]+[x,f]\cdot y, \label{EqNPlrE9}\\
[f,x\cdot y]&=x\cdot[f,y]+[f,x]\cdot y, \label{EqNPlrE10}\\
[x,f\ast y]&=f\ast[x,y]+[x,f]\cdot y, \label{EqNPlrE11}\\
[x,y\ast f]&=y\cdot[x,f]+[x,y]\ast f, \label{EqNPlrE12}
\end{align}
for all $x,y\in X$.

The first three identities say that
\[
(f\ast-,-\ast f)\in\Bim(X,\cdot),
\]
while \eqref{EqNPlrE4}--\eqref{EqNPlrE6} state that $([f,-],[-,f])$ belongs to the external weak actor of the Leibniz algebra $(X,[-,-])$, that is, by \Cref{ExampleExternalLeib},
\[
(-[-,f],[f,-]) \in \Bider(X,[-,-]).
\]
The remaining identities express the compatibility between the associative multiplication and the Leibniz bracket.

With the following choice of identities for Axiom~2:
\begin{align}
(x\cdot y)\cdot z&=x\cdot(y\cdot z), \label{EqNPlrAxiom1}\\
z\cdot(x\cdot y)&=(z\cdot x)\cdot y, \label{EqNPlrAxiom2}\\
[[x,y],z]&=[[x,z],y]+[x,[y,z]], \label{EqNPlrAxiom7}\\
[z,[x,y]]&=[[z,x],y]-[[z,y],x], \label{EqNPlrAxiom8}\\
[x\cdot y,z]&=x\cdot[y,z]+[x,z]\cdot y, \label{EqNPlrAxiom3}\\
[z,x\cdot y]&=x\cdot[z,y]+[z,x]\cdot y, \label{EqNPlrAxiom4}\\
[x,y]\cdot z&=[x,y\cdot z]-y\cdot[x,z], \label{EqNPlrAxiom5}\\
z\cdot[x,y]&=[z\cdot x,y]-[z,y]\cdot x, \label{EqNPlrAxiom6}
\end{align}
we get that
\[
\langle f,g\rangle_\ast=
\bigl(
f\ast(g\ast-),
(-\ast f)\ast g,
f\ast[g,-]+[f,-]\ast g,
f\ast[-,g]+[-,f]\ast g
\bigr)
\]
and
\begin{align*}
\langle f,g\rangle_{[-,-]}=
\bigl(
[f,g\ast-]-g\ast[f,-]&,[-\ast f,g]-[-,g]\ast f,\\
&[[f,-],g]+[f,[g,-]],
[[-,f],g]-[[-,g],f]
\bigr).
\end{align*}

The following example shows that the external weak actor constructed above need not satisfy the defining identities of the variety $\NP^{lr}$.

\begin{example}\label{ExampleExternalNPlrNotLeibniz}
Let $X=\F$ be the one-dimensional algebra with trivial multiplication and trivial bracket. Thus, $\E(X)=\End(X)^4$ is a total algebra. Consider the element $f=(0,0,\id_X,0)\in\E(X)$. We have
\[
\langle f,f\rangle_{[-,-]}=\bigl(
0,
0,
\id_X\circ\id_X,
0
\bigr)=f.
\]
Denote this operation simply by $[-,-]_{\E(X)}$. If it were a right Leibniz bracket, then
\[
[[f,f]_{\E(X)},f]_{\E(X)}
=
[[f,f]_{\E(X)},f]_{\E(X)}
+
[f,[f,f]_{\E(X)}]_{\E(X)}.
\]
Since $[f,f]_{\E(X)}=f$, we obtain $f=2f$, which is a contradiction because $f\neq0$. Therefore, the second operation on $\E(X)$ does not satisfy the right Leibniz identity, $\E(X)$ is not a left-right non-commutative Poisson algebra, and we have the following.
\end{example}

\begin{theorem}\label{TheoremNPlrNotAR}
The variety $\NP^{lr}$ of left-right non-commutative Poisson algebras is not action representable.
\end{theorem}

\begin{proof}
By \Cref{ExampleExternalNPlrNotLeibniz}, the external weak actor of the one-dimensional algebra $X=\F$ with trivial multiplication and trivial bracket is a total algebra which does not satisfy the right Leibniz identity, a multilinear identity of $\NP^{lr}$ of degree~$3$. By \Cref{CorTensorObstruction}(2), the eight-dimensional algebra $X\otimes\Lambda_3$, with trivial multiplication and trivial bracket, does not admit an actor. Hence, $\NP^{lr}$ is not action representable.
\end{proof}

\begin{remark}\label{RemarkNPlrTotalityOpen}
The representability of actions in $\NP^{lr}$ was studied in~\cite[Section~4]{CasasDatuashviliLadraNPLR}, where the authors construct the universal strict general actor $\USGA(X)$ and study conditions under which an actor exists. However, they did not provide an explicit description of $\USGA(X)$. As in the case of algebras with a single bilinear multiplication, the external weak actor admits an explicit description, and \Cref{ExampleExternalNPlrNotLeibniz} shows that it suffices to compute it for the one-dimensional algebra with trivial operations in order to prove that $\NP^{lr}$ is not action representable.
\end{remark}

\subsection{Associative dialgebras}\label{SecExternalDias}

Let $X=(X,\dashv,\vdash)$ be an associative dialgebra. Recall that a \emph{tetramultiplier} of $X$ is a quadruple
\[
f=(f\dashv-,-\dashv f,f\vdash-,-\vdash f)\in\End(X)^4
\]
satisfying the fifteen identities obtained from the five defining identities of associative dialgebras by replacing, in turn, one of the three variables by $f$; see~\cite[Definition~2.4]{CasasEtAlTetra}. More explicitly,
\[
(f\dashv-,-\dashv f)\in\Bim(X,\dashv),
\qquad
(f\vdash-,-\vdash f)\in\Bim(X,\vdash)
\]
and
\begin{align*}
f\dashv(x\vdash y)&=(f\dashv x)\dashv y,\\
(x\dashv f)\dashv y&=x\dashv(f\vdash y),\\
(x\dashv y)\dashv f&=x\dashv(y\vdash f),\\
f\vdash(x\dashv y)&=(f\vdash x)\dashv y,\\
(x\vdash f)\dashv y&=x\vdash(f\dashv y),\\
(x\vdash y)\dashv f&=x\vdash(y\dashv f),\\
f\vdash(x\vdash y)&=(f\dashv x)\vdash y,\\
(x\dashv f)\vdash y&=x\vdash(f\vdash y),\\
(x\dashv y)\vdash f&=x\vdash(y\vdash f),
\end{align*}
for any $x,y \in X$. We denote by $\Tetra(X)$ the vector space of tetramultipliers of~$X$. It follows that
\[
\E(X)=\Tetra(X)
\]
as vector spaces. Moreover, with the following choice of identities for Axiom~2, each of which is one of the defining identities of $\Dias$,
\begin{align*}
(x\dashv y)\dashv z&=x\dashv(y\vdash z), &
z\dashv(x\dashv y)&=(z\dashv x)\dashv y,\\
(x\dashv y)\vdash z&=x\vdash(y\vdash z), &
z\vdash(x\dashv y)&=(z\vdash x)\dashv y,\\
(x\vdash y)\dashv z&=x\vdash(y\dashv z), &
z\dashv(x\vdash y)&=(z\dashv x)\dashv y,\\
(x\vdash y)\vdash z&=x\vdash(y\vdash z), &
z\vdash(x\vdash y)&=(z\dashv x)\vdash y,
\end{align*}
the two partial operations on $\E(X)$ are given by
\[
\langle f,g\rangle_{\dashv}
=
\bigl(
f\dashv(g\vdash-),
(-\dashv f)\dashv g,
f\vdash(g\vdash-),
(-\vdash f)\dashv g
\bigr)
\]
and
\[
\langle f,g\rangle_{\vdash}
=
\bigl(
f\vdash(g\dashv-),
(-\dashv f)\dashv g,
f\vdash(g\vdash-),
(-\dashv f)\vdash g
\bigr),
\]
and they coincide with the dialgebra operations on $\Tetra(X)$ described in~\cite{CasasEtAlTetra}. We thus get the following.

\begin{proposition}\cite{CasasEtAlTetra}\label{PropExternalDias}
The external weak actor $\E(X) \cong \Tetra(X)$ is an associative dialgebra. \noproof
\end{proposition}

As an immediate consequence, we obtain a new example of a weakly action representable category.

\begin{theorem}\label{TheoremDiasWRA}
The variety $\Dias$ of associative dialgebras is weakly action representable.
\end{theorem}

\begin{proof}
By \Cref{PropExternalDias}, $\E(X)$ is an object of $\Dias$ for every associative dialgebra~$X$. The claim then follows from \Cref{ExternalActorTwoWRA}.
\end{proof}

\subsection{Matching dialgebras}\label{SecExternalMatching}

Let $X=(X,\ast,\bullet)$ be a matching dialgebra. Thus, both multiplications are associative and satisfy
\[
(x\ast y)\bullet z=x\ast(y\bullet z)
\qquad\text{and}\qquad
(x\bullet y)\ast z=x\bullet(y\ast z)
\]
for all $x,y,z\in X$.

The external weak actor $\E(X)$ consists of all quadruples
\[
f=(f\ast-,-\ast f,f\bullet-,-\bullet f)\in\End(X)^4
\]
such that
\[
(f\ast-,-\ast f)\in\Bim(X,\ast),
\qquad
(f\bullet-,-\bullet f)\in\Bim(X,\bullet),
\]
and
\begin{align*}
(f\ast x)\bullet y&=f\ast(x\bullet y), & (x\ast f)\bullet y&=x\ast(f\bullet y), & (x\ast y)\bullet f&=x\ast(y\bullet f),\\
(f\bullet x)\ast y&=f\bullet(x\ast y), & (x\bullet f)\ast y&=x\bullet(f\ast y), & (x\bullet y)\ast f&=x\bullet(y\ast f),
\end{align*}
for all $x,y\in X$.

With the following choice of identities for Axiom~2,
\[
(x\circ y)\diamond z=x\circ(y\diamond z),
\qquad
\circ,\diamond\in\{\ast,\bullet\},
\]
the two bilinear partial operations on $\E(X)$ are given by
\[
\langle f,g\rangle_\ast=
\bigl(
f\ast(g\ast-),
(-\ast f)\ast g,
f\ast(g\bullet-),
(-\bullet f)\ast g
\bigr)
\]
and
\[
\langle f,g\rangle_\bullet=
\bigl(
f\bullet(g\ast-),
(-\ast f)\bullet g,
f\bullet(g\bullet-),
(-\bullet f)\bullet g
\bigr),
\]
whenever the resulting quadruples belong to $\E(X)$.

\begin{proposition}\label{PropExternalMatching}
The external weak actor $\E(X)$ is a matching dialgebra.
\end{proposition}

\begin{proof}
Let
\[
\circ_1=\ast
\qquad\text{and}\qquad
\circ_2=\bullet.
\]
The identities defining $\E(X)$ can be written as
\begin{align}
f\circ_i(x\circ_j y)&=(f\circ_i x)\circ_j y, \label{EqMatchingUniform1}\\
(x\circ_j y)\circ_i f&=x\circ_j(y\circ_i f), \label{EqMatchingUniform2}\\
x\circ_j(f\circ_i y)&=(x\circ_j f)\circ_i y, \label{EqMatchingUniform3}
\end{align}
for all $x,y\in X$ and $i,j\in\{1,2\}$.

Let $f,g\in\E(X)$, let $k\in\{1,2\}$, and set $u=\langle f,g\rangle_{\circ_k}$. By definition,
\[
u\circ_i x=f\circ_k(g\circ_i x)
\qquad\text{and}\qquad
x\circ_i u=(x\circ_i f)\circ_k g
\]
for every $x\in X$ and $i\in\{1,2\}$. Using \eqref{EqMatchingUniform1}, we obtain
\[
u\circ_i(x\circ_j y)
=f\circ_k\bigl(g\circ_i(x\circ_j y)\bigr)
=f\circ_k\bigl((g\circ_i x)\circ_j y\bigr)
=(u\circ_i x)\circ_j y.
\]
Similarly, \eqref{EqMatchingUniform2} gives
\[
(x\circ_j y)\circ_i u
=\bigl((x\circ_j y)\circ_i f\bigr)\circ_k g
=\bigl(x\circ_j(y\circ_i f)\bigr)\circ_k g
=x\circ_j(y\circ_i u).
\]
Finally, by \eqref{EqMatchingUniform3}, we get
\[
x\circ_j(u \circ_i y)
=x\circ_j\bigl(f\circ_k(g\circ_i y)\bigr)
=(x\circ_j f)\circ_k(g\circ_i y)
=\bigl((x\circ_j f)\circ_k g\bigr)\circ_i y
=(x\circ_j u)\circ_i y.
\]
Therefore, $u \in\E(X)$, and $\E(X)$ is a total algebra. It remains to verify that it satisfies the defining identities of $\Match$. We must prove that
\[
\langle\langle f,g\rangle_{\circ_k},h\rangle_{\circ_\ell}
=
\langle f,\langle g,h\rangle_{\circ_\ell}\rangle_{\circ_k}
\]
for all $f,g,h\in\E(X)$ and $k,\ell\in\{1,2\}$. Let
\[
p=\langle\langle f,g\rangle_{\circ_k},h\rangle_{\circ_\ell}
\qquad\text{and}\qquad
q=\langle f,\langle g,h\rangle_{\circ_\ell}\rangle_{\circ_k}.
\]
To prove that $p=q$, it is enough to show that their four endomorphism components coincide. Let $i\in\{1,2\}$ and $x\in X$. By the definition of the operations on $\E(X)$, we have
\begin{align*}
p\circ_i x
&=\langle f,g\rangle_{\circ_k}\circ_\ell(h\circ_i x)\\
&=f\circ_k\bigl(g\circ_\ell(h\circ_i x)\bigr),
\end{align*}
whereas
\begin{align*}
q\circ_i x
&=f\circ_k\bigl(\langle g,h\rangle_{\circ_\ell}\circ_i x\bigr)\\
&=f\circ_k\bigl(g\circ_\ell(h\circ_i x)\bigr).
\end{align*}
Thus, the left components of $p$ and $q$ corresponding to $\circ_i$ coincide. Similarly,
\begin{align*}
x\circ_i p
&=\bigl(x\circ_i\langle f,g\rangle_{\circ_k}\bigr)\circ_\ell h\\
&=\bigl((x\circ_i f)\circ_k g\bigr)\circ_\ell h,
\end{align*}
while
\begin{align*}
x\circ_i q
&=(x\circ_i f)\circ_k\langle g,h\rangle_{\circ_\ell}\\
&=\bigl((x\circ_i f)\circ_k g\bigr)\circ_\ell h.
\end{align*}
Hence, the right components corresponding to $\circ_i$ also coincide. Since this holds for $i=1,2$, we conclude that $p=q$. Therefore,
\[
\langle\langle f,g\rangle_{\circ_k},h\rangle_{\circ_\ell}
=
\langle f,\langle g,h\rangle_{\circ_\ell}\rangle_{\circ_k},
\]
so $\E(X)$ is a matching dialgebra.
\end{proof}

As an immediate consequence, we obtain another example of a weakly action representable variety.

\begin{theorem}
The variety $\Match$ of matching dialgebras is weakly action representable.
\end{theorem}

\begin{proof}
By \Cref{PropExternalMatching}, $\E(X)$ is an object of $\Match$ for every matching dialgebra~$X$. The claim then follows from \Cref{ExternalActorTwoWRA}.
\end{proof}

Although $\Dias$ and $\Match$ are weakly action representable, they are not action representable.

\begin{theorem}\label{TheoremDiasMatchNotAR}
The varieties $\Dias$ and $\Match$ are not action representable.
\end{theorem}

\begin{proof}
Let $\V$ be one of the varieties $\Dias$ and $\Match$, and let $\ast$ denote, respectively, the multiplication $\dashv$ of $\Dias$ or the multiplication~$\ast$ of $\Match$, which is associative. Let $X$ be the two-dimensional vector space with trivial operations considered in the proof of \Cref{TheoremPoisNotAR}. Since $\V$ satisfies no non-trivial identity of degree~$2$, \Cref{RemarkExternalActorAbelian} gives $\E(X)=\End(X)^4$. The argument used for $\Pois$ in the proof of \Cref{TheoremPoisNotAR} then applies verbatim, and shows that the sixteen-dimensional algebra $X\otimes\Lambda_3$ with trivial operations does not admit an actor in $\V$.
\end{proof}

\subsection{Matching Lie algebras}\label{SecExternalMatchingLie}

We conclude this section with matching Lie algebras. Let $\bchar(\F) \neq 2$, let $X=(X,[-,-],\{-,-\})$ be an object of $\MatchLie$, and write
\[
[-,-]_1=[-,-],
\qquad
[-,-]_2=\{-,-\}.
\]
For Axiom~2, we choose the identities
\[
[[x,y]_i,z]_j
=
[[x,z]_i,y]_j+[x,[y,z]_j]_i,
\]
for all $i,j\in\{1,2\}$. When $i=j$, these are the usual Jacobi identities, while the cases in which $i\neq j$ are equivalent to the matching Jacobi identities.

Since both brackets are anti-commutative, the external weak actor $\E(X)$ consists of all pairs
\[
f=([f,-],\{f,-\})\in\End(X)^2
\]
satisfying
\begin{align*}
[f,[x,y]]
&=[[f,x],y]+[x,[f,y]],\\
\{f,\{x,y\}\}
&=\{\{f,x\},y\}+\{x,\{f,y\}\},\\
[f,\{x,y\}]
&=\{[f,x],y\}+\{x,[f,y]\},\\
\{f,[x,y]\}
&=[\{f,x\},y]+[x,\{f,y\}],\\
\{[f,x],y\}
&=[\{f,x\},y],
\end{align*}
for all $x,y\in X$. Thus, both $[f,-]$ and $\{f,-\}$ are derivations with respect to both Lie brackets, and the last identity expresses the additional compatibility arising from the matching Jacobi identities.

The two bilinear partial operations on $\E(X)$ are given by
\[
\langle f,g\rangle_{[-,-]}
=
\bigl(
[f,[g,-]]-[g,[f,-]],
[f,\{g,-\}]-\{g,[f,-]\}
\bigr)
\]
and
\[
\langle f,g\rangle_{\{-,-\}}
=
\bigl(
\{f,[g,-]\}-[g,\{f,-\}],
\{f,\{g,-\}\}-\{g,\{f,-\}\}
\bigr).
\]
The following example shows that $\E(X)$ is not in general a matching Lie algebra.

\begin{example}\label{ExampleExternalMatchingLie}
Let $X$ be the two-dimensional vector space with basis $\{e_1,e_2\}$, endowed with two trivial Lie brackets. Clearly, $X$ is a matching Lie algebra, and
\[
\E(X)=\End(X)^2,
\]
so it is a total algebra.

Consider the linear endomorphisms $D,F\in\End(X)$ defined by
\[
D(e_1)=e_1,
\qquad
D(e_2)=0,
\]
and
\[
F(e_1)=0,
\qquad
F(e_2)=e_1.
\]
Then
\[
D\circ F=F,
\qquad
F\circ D=0.
\]
Let $f=(D,F)\in\E(X)$. It follows that
\[
\langle f,f\rangle_{[-,-]}
=
\bigl(
D\circ D-D\circ D,
D\circ F-F\circ D
\bigr)
=
(0,F)\neq0.
\]
Therefore, since $\bchar(\F)\neq2$, the bracket $\langle-,-\rangle_{[-,-]}$ is not skew-symmetric, that is, $\E(X)$ does not satisfy the multilinear identity $[x,y]+[y,x]=0$ of $\MatchLie$. In particular, $\E(X)$ is not a matching Lie algebra.
\end{example}

We thus have the following.

\begin{theorem}\label{TheoremMatchingLieNotAR}
Let $\F$ be a field with $\bchar(\F)\neq2$. Then the variety $\MatchLie$ of matching Lie algebras over $\F$ is not action representable.
\end{theorem}

\begin{proof}
By \Cref{ExampleExternalMatchingLie}, the external weak actor of the two-dimensional algebra~$X$ with trivial brackets does not satisfy the multilinear identity $[x,y]+[y,x]=0$ of $\MatchLie$, which has degree~$2$. By \Cref{CorTensorObstruction}(2), the eight-dimensional algebra $X\otimes\Lambda_2$ with trivial brackets does not admit an actor. Hence, $\MatchLie$ is not action representable.
\end{proof}

\section{Action representability of varieties with \texorpdfstring{$n$}{n} multiplications}\label{SecNLieBracketsAR2}

The non-action-representability of matching Lie algebras proved in \Cref{TheoremMatchingLieNotAR} is a special case of a much more general result. In this section we associate with every variety of algebras with $n$ bilinear multiplications a pair of non-negative integers $(p,q)$, called its \emph{reduced type}, and we classify the action representable operadic varieties whose reduced type satisfies $p=0$, $q=0$ or $(p,q)=(1,1)$. Throughout,~$\F$ is a field with $\bchar(\F)\neq2$; we shall need $\F$ to be infinite only for the reduced type~$(1,1)$. The results of this section are compared with the case of a single multiplication in \Cref{RemarkComparisonSingle}.

Our first observation is that, up to term equivalence, one may always assume that the multiplications of $\V$ are linearly independent and that each of them is either commutative or anti-commutative.

\begin{lemma}\label{LemReduction}
Let $\V$ be a variety of $\F$-algebras with $n$ bilinear multiplications $\ast_1,\ldots,\ast_n$. Then $\V$ is term equivalent to a variety $\V'$ of $\F$-algebras with $m\leq2n$ bilinear multiplications $\circ_1,\ldots,\circ_m$ such that
\begin{enumerate}
\item each $\circ_a$ is either commutative or anti-commutative, say
\[
x\circ_ay=\varepsilon_a\,y\circ_ax,
\qquad
\varepsilon_a\in\{1,-1\},
\qquad
a=1,\ldots,m;
\]
\item the $m$ multiplications are linearly independent, that is,
\[
\sum_{a=1}^{m}\lambda_a\,x\circ_ay=0
\]
is an identity of $\V'$ only when $\lambda_1=\cdots=\lambda_m=0$.
\end{enumerate}
Moreover, $\V'$ is operadic as soon as $\V$ is, and the number $p$ of commutative and the number $q=m-p$ of anti-commutative multiplications among $\circ_1,\ldots,\circ_m$ depend only on $\V$, and not on the choice of $\V'$.
\end{lemma}

\begin{proof}
Let $F$ be the free algebra of $\V$ on two generators $x$ and $y$, and let $W\subseteq F$ be the subspace spanned by the $2n$ elements
\[
x\ast_1y,\ \ldots,\ x\ast_ny,
\qquad
y\ast_1x,\ \ldots,\ y\ast_nx .
\]
An element of $W$ is thus of the form
\[
w=w(x,y)=\sum_{i=1}^{n}\bigl(\lambda_i\,x\ast_iy+\mu_i\,y\ast_ix\bigr),
\qquad
\lambda_i,\mu_i\in\F.
\]
For $u,v$ in an object $X$ of $\V$, let $w(u,v)$ denote the image of $w$ under the unique homomorphism $F\to X$ sending $x$ to $u$ and $y$ to $v$, which exists because $F$ is free. Then
\[
w(u,v)=\sum_{i=1}^{n}\bigl(\lambda_i\,u\ast_iv+\mu_i\,v\ast_iu\bigr),
\]
so that $w$ determines a bilinear multiplication $(u,v)\mapsto w(u,v)$ on every object of~$\V$.

Let $\sigma\colon W\to W$ be the restriction to $W$ of the endomorphism of $F$ which interchanges the two generators $x$ and $y$. It is well defined because $\sigma(x\ast_iy)=y\ast_ix$ and $\sigma(y\ast_ix)=x\ast_iy$ for every $i=1,\ldots,n$. In the coordinates above, $\sigma$ is given by $(\lambda,\mu)\mapsto(\mu,\lambda)$, so that it is a linear involution of $W$. Since $\bchar(\F)\neq2$, we have
\[
W=W^+\oplus W^-,
\]
where $W^{\pm}=\{w\in W\mid\sigma(w)=\pm w\}$. Explicitly, $W^+$ is spanned by the elements $x\ast_iy+y\ast_ix$ and $W^-$ by the elements $x\ast_iy-y\ast_ix$. If $w\in W$ satisfies $\sigma(w)=\varepsilon\,w$ with $\varepsilon\in\{1,-1\}$, and if $\circ$ denotes the corresponding multiplication, then
\[
v\circ u=w(v,u)=\sigma(w)(u,v)=\varepsilon\,w(u,v)=\varepsilon\,(u\circ v)
\]
for all $u,v\in X$, so that $\circ$ is commutative when $\varepsilon=1$ and anti-commutative when $\varepsilon=-1$.

Let now $\circ_1,\ldots,\circ_p$ be a basis of $W^+$ and $\circ_{p+1},\ldots,\circ_m$ a basis of $W^-$, so that $m=\dim_\F W\leq2n$ and each $\circ_a$ is a bilinear multiplication on the objects of $\V$, commutative if $a\leq p$ and anti-commutative otherwise.

Every $\ast_i$ lies in $W$, hence is a linear combination of $\circ_1,\ldots,\circ_m$, and conversely every $\circ_a$ is a linear combination of the $x\ast_iy$ and the $y\ast_ix$. Therefore the two maps
\[
(X,\ast_1,\ldots,\ast_n)\longmapsto(X,\circ_1,\ldots,\circ_m),
\qquad
(X,\circ_1,\ldots,\circ_m)\longmapsto(X,\ast_1,\ldots,\ast_n)
\]
are given by terms and are mutually inverse, so that they define a term equivalence between $\V$ and the class $\V'$ of the algebras $(X,\circ_1,\ldots,\circ_m)$, where $X$ is an object of~$\V$. Since a term equivalence preserves products, subalgebras and quotients, $\V'$ is a variety. Moreover, it is operadic when $\V$ is, because its identities are obtained from those of~$\V$ through the substitutions above, which preserve multilinearity, and~(2) holds because the bilinear operations $\circ_a$ form a basis of $W$. Finally, $W$, $W^+$ and $W^-$ are constructed from $\V$ alone, so that $m=\dim_\F W$, $p=\dim_\F W^+$ and $q=\dim_\F W^-$ do not depend on the choice of the bases of $W^+$ and $W^-$, that is, on the choice of $\V'$.
\end{proof}

\begin{definition}\label{DefReducedType}
The pair $(p,q)$ of \Cref{LemReduction} is called the \emph{reduced type} of $\V$.
\end{definition}

The terminology is meant in the sense of universal algebra: $(p,q)$ is the signature of the reduced variety $\V'$, which has $p$ commutative and $q$ anti-commutative bilinear multiplications. In particular, $\V$ has reduced type $(p,0)$ precisely when all the multiplications of $\V'$ are commutative, and reduced type $(0,q)$ precisely when they are all anti-commutative. Thus a variety with $n$ commutative multiplications has reduced type $(p,0)$ with $p\leq n$, and a variety with $n$ anti-commutative multiplications has reduced type $(0,q)$ with $q\leq n$. Both inequalities may be strict, because the multiplications of $\V'$ are linearly independent: a multiplication of $\V$ which is identically zero, or which is a linear combination of the others, does not contribute to $p$ or to $q$. For instance, $\AbAlg^n$ has reduced type $(0,0)$ for every $n$, whereas $\Lie^n$ has reduced type $(0,n)$ because its $n$ brackets are linearly independent.

For the term equivalence constructed in the proof of \Cref{LemReduction}, the varieties~$\V$ and $\V'$ have the same reduced type. Indeed, a term equivalence preserves free algebras, so that the free algebra of $\V'$ on $x$ and $y$ is again $F$, and the involution associated with $\V'$ is again $\sigma$. Moreover, the subspace associated with $\V'$ is spanned by the elements $x\circ_ay$ and $y\circ_ax=\varepsilon_a\,x\circ_ay$, hence equals $W$, because the $\circ_a$ form a basis of $W$.

\begin{example}\label{ExReducedType}
A variety with one bilinear multiplication has reduced type $(0,0)$ if its multiplication is trivial, $(1,0)$ if it is commutative and non-trivial, $(0,1)$ if it is anti-commutative and non-trivial, and $(1,1)$ otherwise. In the last case the two multiplications of \Cref{LemReduction} are
\[
x\odot y=\tfrac12(xy+yx),
\qquad
[x,y]=\tfrac12(xy-yx).
\]
Thus $\AbAlg$ has reduced type $(0,0)$, $\Lie$ has reduced type $(0,1)$, $\Com$ and $\CAssoc$ have reduced type $(1,0)$, and $\Assoc$ has reduced type $(1,1)$. Moreover, $\Lie^n$ has reduced type $(0,n)$, $\CPois$ has reduced type $(1,1)$, $\Pois$ has reduced type $(1,2)$, and $\Dias$ and $\Match$ have reduced type $(2,2)$.
\end{example}

By \Cref{LemReduction} there is no loss of generality in the following assumptions, which we adopt for the rest of the section: $\V$ is an operadic variety of algebras with $m$ linearly independent bilinear multiplications $\circ_1,\ldots,\circ_m$ and
\[
x\circ_ay=\varepsilon_a\,y\circ_ax,
\qquad
\varepsilon_a\in\{1,-1\}, \qquad
a=1,\ldots,m.
\]
We write $E$ for the diagonal $m^2 \times m^2$ matrix with rows and columns indexed by the ordered pairs $(a,b)$, whose diagonal entry at $(a,b)$ is $\varepsilon_a$. Thus $E^2=I_{m^2}$, and multiplying a matrix by $E$ on the left, respectively on the right, multiplies its row, respectively its column, of index $(a,b)$ by $\varepsilon_a$.

We outline the strategy of the rest of this section. If an infinite-dimensional algebra~$X$ of $\V$ with trivial multiplications admits an actor, then, by \Cref{CorTensorObstruction}(1), $\E(X)$ is total and the semidirect product $\E(X)\ltimes X$ is an object of $\V$. In particular, it satisfies the chosen identities for Axiom~2. Taking for $X$ the underlying vector space of a free associative algebra, and comparing the coefficients of linearly independent composites of left multiplications by its generators, we obtain two matrix equations, $A^2=I_{m^2}$ and $AEA=-EAE$, for the matrix $A$ of the coefficients of the identities for Axiom~2. When $p=0$ or $q=0$, these equations determine $A$ completely and force $\V$ to satisfy cyclic mixed Jacobi identities, from which we derive a contradiction unless $\V$ has reduced type $(0,0)$ or $(0,1)$. In the latter case, $\V$ turns out to be term equivalent to $\Lie$.

We now choose a set of identities for Axiom~2 for $\V$. For $\eta\in\{1,-1\}$ and $c,d\in\{1,\ldots,m\}$ we consider the multilinear polynomial
\[
M^{\eta}_{c,d}(x,y,z)=(x\circ_cz)\circ_dy+\eta\,(y\circ_cz)\circ_dx.
\]
It is easy to check that $M^{\eta}_{c,d}(y,x,z)=\eta\,M^{\eta}_{c,d}(x,y,z)$.

\begin{proposition}\label{TheoremIndependentLieBracketsNotAR2}
Let $\F$ be a field with $\bchar(\F)\neq2$ and let $\V$ be as above. If $\V$ is an Orzech category of interest, then a choice of identities for Axiom~2 can be made in the form
\begin{equation}\label{EqMixedAxiomLieN2}
(x\circ_ay)\circ_bz
=
\sum_{c,d=1}^{m}\alpha^{ab}_{cd}\,M^{\varepsilon_a}_{c,d}(x,y,z),
\qquad
1\leq a,b\leq m,
\end{equation}
for suitable $\alpha^{ab}_{cd}\in\F$. Moreover, let $\mathcal A$ be the free unital associative algebra on the $3m$ generators $t_1,\ldots,t_m$, $s_1,\ldots,s_m$, $r_1,\ldots,r_m$, let $X$ be its underlying vector space endowed with $m$ trivial multiplications, and let $A$ be the $m^2 \times m^2$ matrix with entries $\alpha^{ab}_{cd}$, whose rows and columns are indexed by the ordered pairs $(a,b)$ and $(c,d)$. If~$X$ admits an actor in $\V$, then
\begin{equation}\label{EqForcedCoefficients}
A^2=I_{m^2}
\qquad\text{and}\qquad
AEA=-EAE .
\end{equation}
\end{proposition}

\begin{proof}
Let $a,b\in\{1,\ldots,m\}$ and let
\[
(x\circ_ay)\circ_bz=T(x,y,z)
\]
be an identity of $\V$ provided by Axiom~2, so that $T$ is a linear combination of the monomials allowed on the right-hand side of Axiom~2. Interchanging $x$ and $y$, and using the identity $y\circ_ax=\varepsilon_a\,x\circ_ay$, we obtain $\varepsilon_a\,(x\circ_ay)\circ_bz=T(y,x,z)$. Thus,
\[
2\,(x\circ_ay)\circ_bz=T(x,y,z)+\varepsilon_a\,T(y,x,z).
\]
Since $\bchar(\F)\neq2$, we may therefore replace $T$ by one half of its $\varepsilon_a$-symmetrisation in $x$ and $y$. This is again an admissible right-hand side, because the set of monomials allowed by Axiom~2 is stable under the interchange of $x$ and $y$. Now, since each~$\circ_c$ is commutative or anti-commutative, every such monomial is, up to sign, of the form $(x\circ_cz)\circ_dy$ or $(y\circ_cz)\circ_dx$. Hence its $\varepsilon_a$-symmetrisation in $x$ and $y$ is, up to sign, one of the $M^{\varepsilon_a}_{c,d}$. This proves~\eqref{EqMixedAxiomLieN2}.

Now, let $X$ be the underlying vector space of the free unital associative algebra~$\mathcal{A}$, endowed with $m$ trivial multiplications. Since $\V$ is operadic, $X$ is an object of $\V$. Moreover, since the multiplications are linearly independent, the variety $\V$ satisfies no identity of degree two other than those expressing the commutativity or the anti-commutativity of the $\circ_a$. It follows from \Cref{RemarkExternalActorAbelian,RemarkCommutativeOperationTwo} that
\[
\E(X)=\End(X)^m .
\]
For any $f=(F_1,\ldots,F_m)\in\E(X)$ and $x\in X$ we denote
\[
f\circ_ax=F_a(x),
\qquad
x\circ_af=\varepsilon_a\,F_a(x),
\qquad
a=1,\ldots,m .
\]
For any $f=(F_1,\ldots,F_m)$ and $g=(G_1,\ldots,G_m)$ in $\E(X)$, the $m$ components of $\langle f,g\rangle_a$ are obtained by evaluating the right-hand side of~\eqref{EqMixedAxiomLieN2} at $(f,g,x)$. Since
\[
(f\circ_cx)\circ_dg=\varepsilon_d\,G_dF_c(x),
\qquad
(g\circ_cx)\circ_df=\varepsilon_d\,F_dG_c(x),
\]
where juxtaposition denotes composition of endomorphisms, we obtain
\begin{equation}\label{EqPartialOperationSigns}
\bigl(\langle f,g\rangle_a\bigr)_b
=
\theta^{ab}_{f,g}
:=
\sum_{c,d=1}^{m}\varepsilon_d\,\alpha^{ab}_{cd}\bigl(\varepsilon_a\,F_dG_c+G_dF_c\bigr).
\end{equation}
In particular, the $m$ bilinear partial operations on $\E(X)$ are total.

Now, suppose that $X$ admits an actor in $\V$. Since $X$ is infinite-dimensional and~$\E(X)$ is total, \Cref{CorTensorObstruction}(1) shows that the semidirect product $\E(X)\ltimes X$, with respect to the canonical action of $\E(X)$ on $X$, is an object of $\V$. In particular, $\E(X)$ is an object of $\V$, and the identities~\eqref{EqMixedAxiomLieN2} hold in $\E(X)\ltimes X$.

We evaluate them at $(f,x,g)$. Since
\[
(f\circ_ax)\circ_bg=\varepsilon_b\,G_bF_a(x),
\qquad
(x\circ_cg)\circ_df=\varepsilon_c\varepsilon_d\,F_dG_c(x),
\qquad
(f\circ_cg)\circ_dx=\theta^{cd}_{f,g}(x),
\]
and since $x$ is arbitrary, \eqref{EqMixedAxiomLieN2} gives the following identities in $\End(X)$:
\begin{equation}\label{EqSemidirectSigns}
\varepsilon_b\,G_bF_a
=
\sum_{c,d=1}^{m}\alpha^{ab}_{cd}
\bigl(
\varepsilon_c\varepsilon_d\,\varepsilon_a\,F_dG_c
+
\theta^{cd}_{f,g}
\bigr),
\qquad
1\leq a,b\leq m .
\end{equation}
Let $F_a=L_{t_a}$ and $G_a=L_{s_a}$, for $a=1,\ldots,m$, be the left multiplication endomorphisms associated with the generators $t_a$ and $s_a$ of $\mathcal A$. Then the $2m^2$ endomorphisms
\[
F_dG_c,
\qquad
G_dF_c,
\qquad
1\leq c,d\leq m,
\]
are linearly independent. Indeed, if a linear combination of them were zero, evaluating it at the unit of $\mathcal A$ would give a linear relation between the $2m^2$ distinct associative monomials $t_ds_c$ and $s_dt_c$, which is impossible.

Substituting~\eqref{EqPartialOperationSigns} into~\eqref{EqSemidirectSigns} and collecting the terms, we obtain
\begin{align*}
\varepsilon_b\,G_bF_a
=
\sum_{c',d'=1}^{m}\varepsilon_{d'}
\Bigl(
\varepsilon_a\varepsilon_{c'}\alpha^{ab}_{c'd'}
+&
\sum_{c,d=1}^{m}\varepsilon_c\,\alpha^{ab}_{cd}\,\alpha^{cd}_{c'd'}
\Bigr)F_{d'}G_{c'}\\
+&
\sum_{c',d'=1}^{m}\varepsilon_{d'}
\Bigl(
\sum_{c,d=1}^{m}\alpha^{ab}_{cd}\,\alpha^{cd}_{c'd'}
\Bigr)G_{d'}F_{c'} .
\end{align*}
Comparing the coefficients of $G_{d'}F_{c'}$ and of $F_{d'}G_{c'}$, and dividing by $\varepsilon_{d'}$, we get
\[
\sum_{c,d=1}^{m}\alpha^{ab}_{cd}\,\alpha^{cd}_{c'd'}=\delta_{(a,b),(c',d')}
\qquad\text{and}\qquad
\varepsilon_a\varepsilon_{c'}\,\alpha^{ab}_{c'd'}+\sum_{c,d=1}^{m}\varepsilon_c\,\alpha^{ab}_{cd}\,\alpha^{cd}_{c'd'}=0
\]
for all $a,b,c',d'$, where $\delta$ denotes the Kronecker delta. The first family is $A^2=I_{m^2}$ and the second one is $EAE+AEA=0$.
\end{proof}

\begin{corollary}\label{CorCriterion}
Let $\V$ be as above, and suppose that $\V$ is an Orzech category of interest. If, for some choice of identities for Axiom~2 of the form~\eqref{EqMixedAxiomLieN2}, the matrix~$A$ does not satisfy~\eqref{EqForcedCoefficients}, then the algebra $X$ of \Cref{TheoremIndependentLieBracketsNotAR2} does not admit an actor. In particular, $\V$ is not action representable.
\end{corollary}

\begin{proof}
This follows from \Cref{TheoremIndependentLieBracketsNotAR2}, whose proof applies to every choice of identities for Axiom~2 of the form~\eqref{EqMixedAxiomLieN2}.
\end{proof}

\begin{remark}\label{RemarkPureTypes}
Suppose that $\varepsilon_1=\cdots=\varepsilon_m=\eta$, that is, that the reduced type of $\V$ is $(m,0)$ or $(0,m)$, and that $X$ admits an actor in $\V$, so that~\eqref{EqForcedCoefficients} holds. Then $E=\eta\,I_{m^2}$, so that $AEA=\eta A^2=\eta\,I_{m^2}$ and $-EAE=-A$. Thus,~\eqref{EqForcedCoefficients} gives
\[
A=-\eta\,I_{m^2}.
\]
Substituting $\alpha^{ab}_{cd}=-\eta\,\delta_{(a,b),(c,d)}$ into~\eqref{EqMixedAxiomLieN2} and using $x\circ_az=\eta\,z\circ_ax$, we find that~$\V$ satisfies the cyclic identities
\begin{equation}\label{EqForcedMixedIdentity2}
(x\circ_ay)\circ_bz+(y\circ_az)\circ_bx+(z\circ_ax)\circ_by=0,
\qquad
1\leq a,b\leq m,
\end{equation}
and that
\begin{equation}\label{EqBracketOnE}
\bigl(\langle f,g\rangle_a\bigr)_b=-\varepsilon_b\bigl(F_bG_a+\eta\,G_bF_a\bigr),
\qquad
1\leq a,b\leq m .
\end{equation}
For $\eta=-1$ the identity~\eqref{EqForcedMixedIdentity2} becomes $(x\circ_ay)\circ_bz=(x\circ_az)\circ_by-(y\circ_az)\circ_bx$, and for $a=b$ it is the Jacobi identity: thus, for a variety of reduced type $(0,m)$, the existence of an actor for $X$ already forces the $m$ multiplications to be Lie brackets. Symmetrically, for $\eta=1$ and $a=b$ the identity~\eqref{EqForcedMixedIdentity2} is the Jacobi identity for the commutative multiplication $\circ_a$. Thus, for a variety of reduced type $(m,0)$, the existence of an actor for $X$ forces the $m$ multiplications to be Jacobi--Jordan products, in the sense of \Cref{ExamplesVarieties}. 
\end{remark}

We now study action representability of varieties of reduced type $(0,q)$ and $(p,0)$, for which, by \Cref{RemarkPureTypes}, the conditions~\eqref{EqForcedCoefficients} determine $A$ completely.

\begin{theorem}\label{TheoremAntiCommOne}
Let $\F$ be a field with $\bchar(\F)\neq2$ and let $\V$ be an operadic variety of algebras of reduced type $(0,1)$. If $\V$ is action representable, then $\V$ is term equivalent to $\Lie$.
\end{theorem}

\begin{proof}
By \Cref{LemReduction}, and since action representability is preserved under term equivalences, we may assume that $\V$ has a single bilinear multiplication, which is anti-commutative, and we prove that $\V=\Lie$. Being action representable, $\V$ is action accessible, hence an Orzech category of interest by \Cref{ThmActionAccessibilityNOperations}, and every object of $\V$ admits an actor. Applying \Cref{RemarkPureTypes} with $\eta=-1$ and $m=1$, the identity~\eqref{EqForcedMixedIdentity2} is the Jacobi identity, so that $\V$ is a subvariety of $\Lie$. Moreover,~\eqref{EqBracketOnE} gives
\[
\langle f,g\rangle=fg-gf .
\]
Let now $M$ be a Lie algebra, let $U(M)$ be its universal enveloping algebra, and let~$Y$ be the underlying vector space of $U(M)$ endowed with the trivial bracket. Then~$Y$ is an object of $\V$, because $\V$ is operadic, and, as in the proof of \Cref{TheoremIndependentLieBracketsNotAR2}, $\E(Y)=\End(Y)$. Since the coefficients $\alpha^{ab}_{cd}$ do not depend on the algebra, its bracket is again $\langle f,g\rangle=fg-gf$. Hence, $\E(Y)$ is the general linear Lie algebra $\mathfrak{gl}(Y)$. Since $\V$ is action representable and $\E(Y)$ is total, \Cref{CorTensorObstruction}(1) shows that $\E(Y)$ is an object of $\V$, that is, $\mathfrak{gl}(Y)$ satisfies every defining identity of $\V$. By the Poincaré--Birkhoff--Witt theorem, the canonical map
$M\to U(M)$ is injective. Moreover, the map
\[
L \colon M\to\End_{\F}(U(M))\colon a\mapsto L_a,
\]
where $L_a$ denotes left multiplication by $a$ in $U(M)$, is a Lie algebra monomorphism, since
\[
[L_a,L_b]=L_{ab-ba}=L_{[a,b]}
\qquad\text{and}\qquad
L_a(1)=a.
\]
Thus, $M$ is a Lie subalgebra of $\mathfrak{gl}(Y)$ and therefore satisfies every defining identity of $\V$. Thus $\Lie\subseteq\V$, and $\V=\Lie$.
\end{proof}

\begin{theorem}\label{TheoremAntiCommNotAR}
Let $\F$ be a field with $\bchar(\F)\neq2$ and let $\V$ be an operadic variety of algebras of reduced type $(0,q)$ with $q\geq2$. Then $\V$ is not action representable.
\end{theorem}

\begin{proof}
By \Cref{LemReduction}, and since action representability and action accessibility are preserved under term equivalences, we may assume that $\V$ satisfies the standing assumptions of this section, with $m=q$ anti-commutative multiplications.

If $\V$ is not action accessible, then it is not action representable.

Now, suppose that $\V$ is action accessible. By \Cref{ThmActionAccessibilityNOperations} it is an Orzech category of interest, so that \Cref{TheoremIndependentLieBracketsNotAR2} and \Cref{RemarkPureTypes} apply. Suppose, towards a contradiction, that the algebra $X$ with trivial multiplications has an actor in $\V$. By \Cref{RemarkPureTypes}, $\V$ satisfies the identities~\eqref{EqForcedMixedIdentity2}. Moreover, since $X$ is infinite-dimensional and $\E(X)=\End(X)^q$ is total, \Cref{CorTensorObstruction}(1) shows that $\E(X)$ is an object of $\V$, and, by~\eqref{EqBracketOnE},
\[
\langle f,g\rangle_a
=
\bigl(F_1G_a-G_1F_a,\ \ldots,\ F_qG_a-G_qF_a\bigr),
\qquad
a=1,\ldots,q .
\]
We choose now $u,v\in\End(X)$ satisfying
\[
u^2=u,\qquad uv=v\neq0,\qquad vu=0 .
\]
For example, take two linearly independent elements $e_1,e_2\in X$, extend them to a basis of $X$, and define $u(e_1)=e_1$, $u(e_2)=0$, $v(e_1)=0$, $v(e_2)=e_1$, both endomorphisms being zero on all the remaining basis elements. Since $q\geq2$, we may consider the elements
\begin{align*}
f&=(F_1,\ldots,F_q)=(0,u,0,\ldots,0),\\
g&=(G_1,\ldots,G_q)=(0,v,0,\ldots,0),\\
h&=(H_1,\ldots,H_q)=(u,0,\ldots,0)
\end{align*}
of $\E(X)=\End(X)^q$. We may check that
\[
\langle f,g\rangle_1=0,
\qquad
\langle f,h\rangle_1=(0,u^2,0,\ldots,0)=f,
\qquad
\langle g,h\rangle_1=(0,vu,0,\ldots,0)=0,
\]
while
\[
\langle f,g\rangle_2=(0,uv-vu,0,\ldots,0)=(0,v,0,\ldots,0)\neq0 .
\]
Since $\E(X)$ is an object of $\V$, it satisfies~\eqref{EqForcedMixedIdentity2}. Evaluating this identity, with $a=1$ and $b=2$, at $(x,y,z)=(f,g,h)$, we obtain
\[
0
=
\bigl\langle\langle f,g\rangle_1,h\bigr\rangle_2
=
\bigl\langle\langle f,h\rangle_1,g\bigr\rangle_2
-
\bigl\langle\langle g,h\rangle_1,f\bigr\rangle_2
=
\langle f,g\rangle_2\neq0,
\]
which is a contradiction. Hence $X$ does not admit an actor in $\V$, and $\V$ is not action representable.
\end{proof}

\begin{remark}\label{RemarkForcedIsMatching2}
Let $q=2$. For $a=1$ and $b=2$, the identity~\eqref{EqForcedMixedIdentity2} can be written as
\[
[[x,y]_1,z]_2=[[x,z]_1,y]_2+[x,[y,z]_1]_2,
\]
whereas the choice of identities for Axiom~2 made in \Cref{SecExternalMatchingLie} is
\[
[[x,y]_1,z]_2=[[x,z]_1,y]_2+[x,[y,z]_2]_1 .
\]
In $\MatchLie$, these two identities are equivalent, because
\[
[[x,y]_i,z]_j=[[x,y]_j,z]_i
\]
is an identity of $\MatchLie$ for any $i,j\in\{1,2\}$: this follows from the matching Jacobi identities by the same computation as in \Cref{DerivedActionMatchingLie}, carried out with all three variables in $X$. The same holds for $a=2$ and $b=1$. Hence, $\MatchLie$ satisfies the identities~\eqref{EqForcedMixedIdentity2}. Moreover, the two brackets of $\MatchLie$ are linearly independent, as one sees by endowing a non-abelian Lie algebra with its bracket and with the trivial bracket, in either order. Thus, $\MatchLie$ has reduced type $(0,2)$, and \Cref{TheoremAntiCommNotAR} gives another proof of \Cref{TheoremMatchingLieNotAR}.
\end{remark}

\begin{theorem}\label{TheoremCommNotAR}
Let $\F$ be a field with $\bchar(\F)\neq2$ and let $\V$ be an operadic variety of algebras of reduced type $(p,0)$ with $p\geq1$. Then $\V$ is not action representable.
\end{theorem}

\begin{proof}
By \Cref{LemReduction}, and since action representability and action accessibility are preserved under term equivalences, we may assume that $\V$ satisfies the standing assumptions of this section, with $m=p$ commutative multiplications.

If $\V$ is not action accessible, then it is not action representable. 

Now, suppose that $\V$ is action accessible. By \Cref{ThmActionAccessibilityNOperations} it is an Orzech category of interest, and we suppose, towards a contradiction, that the algebra $X$ with trivial multiplications has an actor in $\V$. By \Cref{RemarkPureTypes} with $\eta=1$, $\V$ satisfies the identities~\eqref{EqForcedMixedIdentity2}. Since $X$ is infinite-dimensional and $\E(X)=\End(X)^p$ is total, \Cref{CorTensorObstruction}(1) shows that $\E(X)$ is an object of $\V$, so that $\E(X)$ satisfies~\eqref{EqForcedMixedIdentity2} as well. Moreover, by~\eqref{EqBracketOnE},
\[
\bigl(\langle f,g\rangle_a\bigr)_b=-\bigl(F_bG_a+G_bF_a\bigr),
\qquad
1\leq a,b\leq p .
\]
Let $f,g,h\in\E(X)$ and write $u=\langle f,g\rangle_1$, so that $U_1=-(F_1G_1+G_1F_1)$. Then
\[
\bigl(\langle u,h\rangle_1\bigr)_1
=
-\bigl(U_1H_1+H_1U_1\bigr)
=
F_1G_1H_1+G_1F_1H_1+H_1F_1G_1+H_1G_1F_1,
\]
and the first component of~\eqref{EqForcedMixedIdentity2}, with $a=b=1$, evaluated at $(x,y,z)=(f,g,h)$, gives
\[
0
=
2\bigl(F_1G_1H_1+G_1F_1H_1+H_1F_1G_1+H_1G_1F_1+G_1H_1F_1+F_1H_1G_1\bigr),
\]
the second and third rows of the cyclic sum being obtained from the first by the substitutions $(f,g,h)\mapsto(g,h,f)$ and $(f,g,h)\mapsto(h,f,g)$. Choosing now $f,g,h$ with first components $F_1=L_{t_1}$, $G_1=L_{s_1}$, $H_1=L_{r_1}$ and all the remaining components equal to zero, and evaluating at the unit of $\mathcal A$, we obtain
\[
0=2\bigl(t_1s_1r_1+s_1t_1r_1+r_1t_1s_1+r_1s_1t_1+s_1r_1t_1+t_1r_1s_1\bigr),
\]
a sum of six pairwise distinct associative monomials of degree three. Since $\bchar(\F)\neq2$, this is a contradiction. Hence $\V$ is not action representable.
\end{proof}

\begin{remark}\label{RemCommNotAR}
\Cref{TheoremCommNotAR} applies, for instance, to the varieties $\Com$, $\CAssoc$, and $\JJord$.
\end{remark}

When both $p\geq1$ and $q\geq1$, the matrix $E$ is not a scalar matrix, and the conditions~\eqref{EqForcedCoefficients} need not determine $A$. The argument of \Cref{RemarkPureTypes} therefore no longer applies. We settle here the smallest such case, over an infinite field.

\begin{theorem}\cite{Tim}\label{TheoremMixedNotAR}
Let $\F$ be an infinite field with $\bchar(\F)\neq2$ and let $\V$ be an operadic variety of algebras of reduced type $(1,1)$. Then $\V$ is not action representable.
\end{theorem}

\begin{proof}
Write $\odot$ for the commutative and $[-,-]$ for the anti-commutative multiplication of $\V$, and let
\[
xy:=x\odot y+[x,y].
\]
Since $\bchar(\F)\neq2$, we have $x\odot y=\tfrac12(xy+yx)$ and $[x,y]=\tfrac12(xy-yx)$, so that the two maps $(X,\odot,[-,-])\mapsto(X,\cdot\,)$ and $(X,\cdot\,)\mapsto(X,\odot,[-,-])$ are given by terms and are mutually inverse. Hence $\V$ is term equivalent to a variety $\W$ of $\F$-algebras with a single bilinear multiplication.

Suppose that $\V$ is action representable. Then so is $\W$, since action representability is preserved under term equivalences. Since $\F$ is infinite, it follows from~\cite[Theorem~6.1]{Tim} that $\W$ is the variety of Lie algebras, the variety of quasi-Lie algebras, or the variety of abelian algebras, the first two coinciding because $\bchar(\F)\neq2$. In each of these cases $xy+yx=0$ is an identity of $\W$, so that $x\odot y=0$ is an identity of $\V$, contradicting the linear independence of $\odot$ and $[-,-]$. Hence $\V$ is not action representable.
\end{proof}

\begin{remark}\label{RemarkMixedExamples}
\Cref{TheoremMixedNotAR} applies, for instance, to the varieties $\Assoc$ of associative algebras and $\CPois$ of commutative Poisson algebras. In the case of $\CPois$, the term equivalence used in its proof is the polarisation of~\cite{MarklRemm}, and the conclusion holds over an arbitrary field of characteristic different from~$2$ by \Cref{TheoremPoisNotAR}.
\end{remark}

We thus have the following classification theorem.

\begin{theorem}\label{TheoremNLieBracketsAR2}
Let $\F$ be a field with $\bchar(\F)\neq2$, let $n\geq1$, and let $\V$ be an operadic variety of algebras with $n$ bilinear multiplications, of reduced type $(p,q)$. Assume that $p=0$ or $q=0$, or that $(p,q)=(1,1)$ and $\F$ is infinite. Then the following conditions are equivalent:
\begin{tfae}
\item $\V$ is action representable;
\item $\V$ is term equivalent to either $\Lie$ or $\AbAlg$.
\end{tfae}
\end{theorem}

\begin{proof}
By \Cref{LemReduction} we may assume that $\V$ has $m=p+q$ linearly independent multiplications, each commutative or anti-commutative.

$(i)\Rightarrow(ii)$. If $p\geq1$ and $q=0$, then $\V$ is not action representable by \Cref{TheoremCommNotAR}; if $p=0$ and $q\geq2$, by \Cref{TheoremAntiCommNotAR}; and if $(p,q)=(1,1)$ and $\F$ is infinite, by \Cref{TheoremMixedNotAR}. Hence $(p,q)=(0,0)$ or $(p,q)=(0,1)$.

If $(p,q)=(0,0)$, then every multiplication of $\V$ vanishes, so that $\V$ is term equivalent to $\AbAlg$. If instead $(p,q)=(0,1)$, then by \Cref{TheoremAntiCommOne}, $\V$ is term equivalent to $\Lie$.

$(ii)\Rightarrow(i)$. Both $\Lie$ and $\AbAlg$ are action representable, and action representability is preserved under term equivalences.
\end{proof}

\begin{corollary}\label{CorAntiCommClassification}
Let $\F$ be a field with $\bchar(\F)\neq2$ and let $\V$ be an operadic variety of algebras with $n$ anti-commutative bilinear multiplications. Then~$\V$ is action representable if and only if it is term equivalent to either $\Lie$ or $\AbAlg$.
\end{corollary}

\begin{proof}
Since $y\ast_ix=-x\ast_iy$ for every $i$, the reduced type of $\V$ is of the form $(0,q)$ with $q\leq n$. Suppose that $\V$ is action representable. If $q\geq2$, then $\V$ is not action representable by \Cref{TheoremAntiCommNotAR}, a contradiction. If $q=0$, then all the multiplications of $\V$ vanish, so that $\V$ is term equivalent to $\AbAlg$. If $q=1$, then $\V$ is term equivalent to $\Lie$ by \Cref{TheoremAntiCommOne}. The converse follows from \Cref{TheoremNLieBracketsAR2}, whose hypotheses are satisfied because $p=0$.
\end{proof}

\begin{remark}\label{RemarkComparisonSingle}
We collect here the relations between the results of this section and those of~\cite{Tim}, where varieties of algebras with a single multiplication are studied over an infinite field. By \Cref{ExReducedType}, such a variety has reduced type $(0,0)$, $(1,0)$, $(0,1)$ or $(1,1)$. \Cref{LemReduction} plays the role of~\cite[Proposition~1.3]{Tim}, without any hypothesis on the identities of~$\V$. For operadic varieties, \Cref{TheoremCommNotAR}, with $p=1$, and \Cref{TheoremAntiCommOne} recover~\cite[Propositions~4.1, 4.3 and~5.1]{Tim}. The proofs of these propositions are computer-free as well. Ours differ from them in that they rest on the matrix equations~\eqref{EqForcedCoefficients}, treat any number of commutative, or of anti-commutative, multiplications at once, require no assumption on the cardinality of~$\F$, and do not use the multilinearisation theorem~\cite[Theorem~1.5]{Tim}. The only computer-assisted step of that article is the proof of~\cite[Propositions~3.1 and~3.4]{Tim}, which show that, over an infinite field, an action representable variety satisfies a non-trivial identity of degree two. In our terms, this step amounts to excluding the reduced type~$(1,1)$, which is the only case in which we rely on that article, through \Cref{TheoremMixedNotAR}, and the only one in which $\F$ is required to be infinite. Hence, the question raised in~\cite[Section~6.5]{Tim}, whether a less computationally involved proof of this step exists, remains open.
\end{remark}

\begin{remark}\label{RemarkRepresentations}
The obstructions to action representability given by \Cref{CorTensorObstruction}(2) and \Cref{CorCriterion} only involve algebras with trivial multiplications. If $X$ is such an algebra, then the split extensions of an algebra $B$ of $\V$ by $X$ correspond to the representations of $B$ on the vector space $X$, i.e., to the actions of $B$ on the abelian object~$X$ of $\V$. Following~\cite{Tim}, we say that $\V$ has \emph{representable representations} if every algebra of $\V$ with trivial multiplications admits an actor. Every action representable variety has representable representations, and the results above yield the following stronger conclusions.
\begin{enumerate}
\item By \Cref{CorTensorObstruction}(2), the varieties $\Nil_c(\Assoc)$, with $c\geq2$, $\NP^{lr}$, $\Dias$ and $\Match$ do not have representable representations, and the same holds for $\Pois$, $\CPois$ and $\MatchLie$ when $\bchar(\F)\neq2$; see \Cref{RemarkTensorNecessary}, \Cref{TheoremPoisNotAR}, \Cref{TheoremNPlrNotAR}, \Cref{TheoremDiasMatchNotAR} and \Cref{TheoremMatchingLieNotAR}.
\item Let $\bchar(\F)\neq2$. By the proofs of \Cref{TheoremAntiCommNotAR} and \Cref{TheoremCommNotAR}, an action accessible operadic variety of reduced type $(0,q)$, with $q\geq2$, or $(p,0)$, with $p\geq1$, does not have representable representations since, up to the term equivalence of \Cref{LemReduction}, the algebra $X$ of \Cref{TheoremIndependentLieBracketsNotAR2} does not admit an actor. Action accessibility is used there only through \Cref{ThmActionAccessibilityNOperations}, which ensures that $\V$ is an Orzech category of interest, so that identities for Axiom~2 of the form~\eqref{EqMixedAxiomLieN2} can be chosen.
\item Let $\F$ be infinite, with $\bchar(\F)\neq2$. By the proof of \Cref{TheoremMixedNotAR}, an operadic variety of reduced type $(1,1)$ does not have representable representations, since~\cite[Theorem~6.1]{Tim} is stated for varieties with representable representations.
\end{enumerate}
\end{remark}

\begin{remark}\label{RemarkNoCategoricalDefinition}
\Cref{TheoremNLieBracketsAR2} answers negatively, for all the reduced types it covers, the question raised in~\cite[Section~6.6]{Tim}, whether there exist action representable varieties of non-associative algebras with $n\geq2$ bilinear multiplications. In particular, a variety of algebras with $n$ anti-commutative multiplications is action representable if and only if it is term equivalent to either $\Lie$ or $\AbAlg$. Thus, no such variety admits a categorical characterisation through the representability of its actions, unless it is term equivalent to $\Lie$ or to $\AbAlg$. The reduced types not covered by \Cref{TheoremNLieBracketsAR2} are those with $p\geq1$, $q\geq1$ and $p+q\geq3$; the smallest ones are $(1,2)$ and $(2,1)$. \Cref{CorTensorObstruction} nevertheless applies to several varieties of these reduced types. For instance, the varieties $\Pois$, of reduced type $(1,2)$, and $\Dias$ and $\Match$, of reduced type $(2,2)$, are not action representable by \Cref{TheoremPoisNotAR,TheoremDiasMatchNotAR}.

We expect the conclusion of \Cref{TheoremNLieBracketsAR2} to hold for every reduced type, which would settle this question completely. 

It would also be interesting to determine under which additional assumptions the failure of the external weak actor to be an object of the variety provides an obstruction to weak action representability.

Finally, it is shown in~\cite{MPV} that, for every operadic, action accessible, \emph{unit-closed} variety of non-associative algebras $\V$, the \emph{ideally exact} category~\cite{IdeallyExact} $\V_1$ of unital algebras of $\V$ is action representable, the actor of a unital algebra $X$ being $X$ itself; the same holds for unital Poisson algebras. It would be interesting to know whether an analogous result holds for varieties with several bilinear multiplications, for a suitable notion of unit-closed variety, using the external weak actor of \Cref{SecExternalWeakActorN}.
\end{remark}

\section*{Acknowledgements}
The author is grateful to Tim Van der Linden for his valuable comments and for the discussions during the writing of this article. He also thanks Corentin Vienne for some helpful suggestions.

The author is supported by the ``National Group for Algebraic and Geometric Structures and their Applications'' (GNSAGA -- INdAM), and he is a Postdoctoral Researcher of the Fonds de la Recherche Scientifique--FNRS. 


\section*{Use of generative AI tools}
Generative AI tools (Anthropic's Claude Opus~5.5 and OpenAI's GPT-6~Sol) were used to polish the exposition and to help check the correctness of the proofs; they were not used to obtain any mathematical results. The author reviewed all suggestions made by these tools and takes full responsibility for the content of this article.

\vspace{1.5\baselineskip}


\begin{thebibliography}{99}

\bibitem{AlbanoStefanelli}
A.~Albano and P.~Stefanelli, \emph{Action accessibility in the variety of skew braces}, Communications in Mathematics \textbf{35} (2027), no.~2. \doi{10.46298/cm.18819}

\bibitem{Barr}
M.~Barr, \emph{Exact Categories}, Lecture Notes in Mathematics \textbf{236} (1971), 1--120.

\bibitem{BJK}
F.~Borceux, G.~Janelidze and G.~M. Kelly, \emph{Internal object actions}, Commentationes Mathematicae Universitatis Carolinae \textbf{46} (2005), no.~2, 235--255.

\bibitem{BJK2}
F.~Borceux, G.~Janelidze and G.~M. Kelly, \emph{On the representability of actions in a semi-abelian category}, Theory and Applications of Categories \textbf{14} (2005), no.~11, 244--286.

\bibitem{Bourn}
D.~Bourn, \emph{Normalization equivalence, kernel equivalence and affine categories}, Lecture Notes in Mathematics \textbf{1488} (1991), 43--62.

\bibitem{act_accessible}
D.~Bourn and G.~Janelidze, \emph{Centralizers in action accessible categories}, Cahiers de Topologie et G\'eom\'etrie Diff\'erentielle Cat\'egoriques \textbf{50} (2009), no.~3, 211--232.

\bibitem{WRAAlg}
J.~Brox, X.~García-Martínez, M.~Mancini, T.~Van~der Linden and C.~Vienne, \emph{Weak representability of actions of non-associative algebras}, Journal of Algebra \textbf{669} (2025), 401--444. \doi{10.1016/j.jalgebra.2025.02.007}

\bibitem{Burde}
D.~Burde and A.~Fialowski, \emph{Jacobi--Jordan algebras}, Linear Algebra and its Applications \textbf{459} (2014), 586--594.
\doi{10.1016/j.laa.2014.07.034}

\bibitem{CasasDatuashviliNLP}
J.~M.~Casas and T.~Datuashvili, \emph{Noncommutative Leibniz--Poisson algebras}, Communications in Algebra \textbf{34} (2006), no.~7, 2507--2530.
\doi{10.1080/00927870600651091}

\bibitem{CasasDatuashviliLadraNPLR}
J.~M.~Casas, T.~Datuashvili and M.~Ladra, \emph{Left-right noncommutative Poisson algebras}, Central European Journal of Mathematics \textbf{12} (2014), no.~1, 57--78.
\doi{10.2478/s11533-013-0321-x}

\bibitem{Casas}
J.~M.~Casas, T.~Datuashvili and M.~Ladra, \emph{Universal Strict General Actors and Actors in Categories of Interest}, Applied Categorical Structures \textbf{18} (2010), no.~1, 85--114. \doi{10.1007/s10485-008-9166-z}

\bibitem{CasasEtAlTetra}
J.~M.~Casas, R.~Fernández-Casado, X.~García-Martínez and E.~Khmaladze, \emph{Actor of a crossed module of dialgebras via tetramultipliers}, Hacettepe Journal of Mathematics and Statistics \textbf{50} (2021), no.~4, 1063--1078. \doi{\detokenize{10.15672/hujms.701217}}

\bibitem{CigoliManciniMetere}
A.~S.~Cigoli, M.~Mancini and G.~Metere, \emph{On the representability of actions of Leibniz algebras and Poisson algebras}, Proceedings of the Edinburgh Mathematical Society \textbf{66} (2023), no.~4, 998--1021. \doi{10.1017/S0013091523000548}

\bibitem{CigoliMantovaniMetere}
A.~S.~Cigoli, S.~Mantovani and G.~Metere, \emph{A push forward construction and the comprehensive factorization for internal crossed modules}, Applied Categorical Structures \textbf{22} (2014), no.~5--6, 931--960. \doi{10.1007/s10485-013-9348-1}

\bibitem{Datuashvili}
T.~Datuashvili, \emph{Homologically trivial internal categories in categories of groups with operations}, Applied Categorical Structures \textbf{3} (1995), no.~3, 221--237.
\doi{10.1007/BF00878442}

\bibitem{GaoGuoZhang}
X.~Gao, L.~Guo and Y.~Zhang, \emph{Matching Rota--Baxter algebras, matching dendriform algebras and matching pre-Lie algebras}, Journal of Algebra \textbf{552} (2020), 134--170. \doi{10.1016/j.jalgebra.2020.02.011}

\bibitem{XabiMancini}
X.~García-Martínez and M.~Mancini, \emph{Action accessible and weakly action representable varieties of algebras}, Bollettino dell'Unione Matematica Italiana \textbf{19} (2026), no.~2, 367--374. \doi{10.1007/s40574-025-00496-1}

\bibitem{Tim}
X.~García-Martínez, M.~Tsishyn, T.~Van der Linden and C.~Vienne, \emph{Algebras with representable representations}, Proceedings of the Edinburgh Mathematical Society \textbf{64} (2021), no.~3, 555--573. \doi{10.1017/S0013091521000304}

\bibitem{XabiTim}
X.~García-Martínez and T.~Van~der~Linden, \emph{A characterisation of Lie algebras amongst anti-commutative algebras}, Journal of Pure and Applied Algebra \textbf{223} (2019), 4857--4870. \doi{10.1016/j.jpaa.2019.02.018}

\bibitem{Ger63}
M.~Gerstenhaber, \emph{The cohomology structure of an associative ring}, Annals of Mathematics~\textbf{78} (1963), no.~2, 267--288.

\bibitem{Gray}
J.~R.~A. Gray, \emph{A note on the relationship between action accessible and weakly action representable categories}, Theory and Applications of Categories \textbf{44} (2025), no.~8, 272--276.

\bibitem{WRA}
G.~Janelidze, \emph{Central extensions of associative algebras and weakly action representable categories}, Theory and Applications of Categories \textbf{38} (2022), no.~36, 1395--1408.

\bibitem{IdeallyExact}
G.~Janelidze, \emph{Ideally exact categories}, Theory and Applications of Categories \textbf{41} (2024), no.~11, 414--425.

\bibitem{Semi-Ab}
G.~Janelidze, L.~Márki and W.~Tholen, \emph{Semi-abelian categories}, Journal of Pure and Applied Algebra \textbf{168} (2002), no.~2--3, 367--386. \doi{10.1016/S0022-4049(01)00103-7}

\bibitem{LodayDialgebras}
J.-L.~Loday, \emph{Dialgebras}. In: \emph{Dialgebras and Related Operads}, Lecture Notes in Mathematics \textbf{1763}, Springer, Berlin, 2001, 7--66. \doi{\detokenize{10.1007/3-540-45328-8_2}}

\bibitem{loday}
J.-L.~Loday, \emph{Une version non commutative des algèbres de Lie: les algèbres de Leibniz}, L'Enseignement Math\'ematique \textbf{39} (1993), no.~3--4, 269--293.

\bibitem{MacLane58}
S.~Mac~Lane, \emph{Extensions and obstructions for rings}, Illinois Journal of Mathematics \textbf{2} (1958), no.~3, 316--345. \doi{10.1215/ijm/1255454537}

\bibitem{SurveyActions}
M.~Mancini, \emph{Actions in varieties of non-associative algebras: a survey} (2026), submitted, preprint \href{https://arxiv.org/abs/2610.18302}{arXiv:2610.18302}.

\bibitem{MPV}
M.~Mancini, F.~Piazza and C.~Vienne, \emph{On the representability of actions of unital algebras} (2026), submitted, preprint \href{https://arxiv.org/abs/2503.04488}{arXiv:2503.04488}.

\bibitem{MarklRemm}
M.~Markl and E.~Remm, \emph{Algebras with one operation including Poisson and other Lie-admissible algebras}, Journal of Algebra \textbf{299} (2006), no.~1, 171--189. \doi{10.1016/j.jalgebra.2005.09.018}

\bibitem{Montoli}
A.~Montoli, \emph{Action accessibility for categories of interest}, Theory and Applications of Categories \textbf{23} (2010), no.~1, 7--21.

\bibitem{Orzech}
G.~Orzech, \emph{Obstruction theory in algebraic categories {I} and {II}}, Journal of Pure and Applied Algebra \textbf{2} (1972), no.~4, 287--314 and 315--340. \doi{10.1016/0022-4049(72)90008-4} and \doi{10.1016/0022-4049(72)90009-6}

\bibitem{Osborn}
J.~M.~Osborn, \emph{Varieties of algebras}, Advances in Mathematics~\textbf{8} (1972), 163--369. \doi{10.1016/0001-8708(72)90003-5}

\bibitem{Cosmash}
Ü.~Reimaa, T.~Van~der Linden and C.~Vienne, \emph{Associativity and the cosmash product in operadic varieties of algebras}, Illinois Journal of Mathematics \textbf{67} (2023), no.~3, 563--598. \doi{10.1215/00192082-10678862}

\bibitem{Strohmayer}
H.~Strohmayer, \emph{Operads of compatible structures and weighted partitions}, Journal of Pure and Applied Algebra \textbf{212} (2008), no.~11, 2522--2534. \doi{10.1016/j.jpaa.2008.04.009}

\bibitem{VallettePostLie}
B.~Vallette, \emph{Homology of generalized partition posets}, Journal of Pure and Applied Algebra~\textbf{208} (2007), no.~2, 699--725.
\doi{\detokenize{10.1016/j.jpaa.2006.03.012}}

\bibitem{VdL-NAA}
T.~Van~der~Linden, \emph{Non-associative Algebras}. In: M.~M.~Clementino, A.~Facchini, M.~Gran (eds), \emph{New Perspectives in Algebra, Topology and Categories}, Coimbra Mathematical Texts, vol.~1, Springer, Cham, 2021, 225--258.
\doi{\detokenize{10.1007/978-3-030-84319-9_7}}

\bibitem{Coco_Thesis}
C.~Vienne, \emph{Categorical-algebraic conditions in semi-abelian categories}, PhD Thesis (2025), Université catholique de Louvain, Belgium. \url{http://hdl.handle.net/2078.1/301231}

\bibitem{Vin63}
E.~B. Vinberg, \emph{The theory of convex homogeneous cones}, Trudy Moskovskogo Matematicheskogo Obshchestva \textbf{12} (1963), 303--358.

\bibitem{ZhangBaiGuo}
Y.~Zhang, C.~Bai and L.~Guo, \emph{The category and operad of matching dialgebras}, Applied Categorical Structures \textbf{21} (2013), no.~6, 851--865. \doi{10.1007/s10485-012-9290-7}

\end{thebibliography}
\end{document}